\documentclass[12pt]{article}
\usepackage[utf8]{inputenc}
\usepackage[margin=1in]{geometry}
\usepackage{endnotes}

\usepackage{natbib}
 \bibpunct[, ]{(}{)}{,}{a}{}{,}%

\usepackage{soul}
\usepackage{authblk}
\usepackage{booktabs}
\usepackage{array, makecell}
\usepackage{subcaption}
\usepackage{mathrsfs}
\usepackage{multirow}
\usepackage{mathtools}
\usepackage[title]{appendix}
\usepackage{xcolor}
\usepackage{hyperref}
\usepackage{amsthm, amsfonts, amsmath}
\usepackage{bbm}
\hypersetup{
    colorlinks,
    linkcolor={red!60!black},
    citecolor={blue!50!black},
    urlcolor={blue!50!black}
}
\providecommand{\keywords}[1]
{
  \small	
  \textbf{\textit{Keywords---}} #1
}

\newcommand*\mystrut[1]{\vrule width0pt height0pt depth#1\relax}

\title{Data-Driven Optimization for\\Atlanta Police Zone Design}
\author{Shixiang Zhu}
\author{He Wang}
\author{Yao Xie}
\affil{H. Milton Stewart School of Industrial and Systems Engineering, Georgia Institute of Technology, Atlanta, Georgia 30332, USA}
\date{}

\begin{document}

\maketitle
\begin{abstract}
We present a data-driven optimization framework for redesigning police patrol zones in an urban environment. The objectives are
to rebalance police workload among geographical areas and to reduce response time to emergency calls.  
We develop a stochastic model for police emergency response
by integrating multiple data sources, including police incidents reports, demographic surveys, and traffic data. Using this stochastic model, we optimize zone redesign plans using mixed-integer linear programming.
Our proposed design was implemented 
by the Atlanta Police Department in March 2019. By analyzing data before and after the zone redesign, we show that the new design has reduced the response time to high-priority 911 calls by 5.8\% and the imbalance of police workload among different zones by 43\%.
\end{abstract}



%


\section{Introduction}
\label{sec:introduction}

In large urban areas, police departments often organize their patrol forces by dividing the geographical region of a city into multiple patrol areas called \emph{zones} (or \emph{precincts}), and each zone is further divided into smaller areas called \emph{beats} (or \emph{sectors}) \citep{larson1972urban}. 
The design of patrol zones affects both the demand and the capacity for police services in each zone and beat, as well as the travel time of patrol units---together, these factors will determine the police's response time to emergency calls and crime events. 
Therefore, the design of patrol zones has a critical impact on the efficiency of police operations.

\begin{figure}[!htbp]
\centering
\includegraphics[width=.4\linewidth]{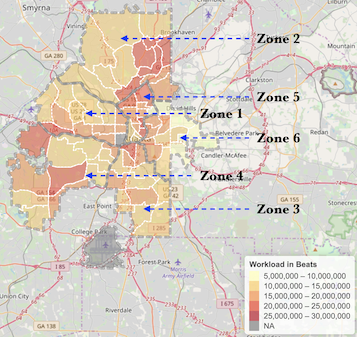}
\caption{Atlanta police zone and beat design (prior to the redesign in March 2019). \emph{Note}: Gray dashed lines represent the  boundaries of the six zones. White lines represent the boundaries of the beats. The color in each beat indicates the average police workload level in 2017 (in seconds).}
\label{fig:workload-beat}
\end{figure}

In this paper, we propose a data-driven framework for designing police zones. The work is developed based on our collaboration with the Atlanta Police Department (APD) for redesigning police zones in Atlanta through a project that lasted from 2017 to 2019.
The APD is the primary police force in metro Atlanta, the ninth-largest metropolitan area in the United States and home to 6 million people.
The APD divides the geographical area of Atlanta into six zones and 81 {beats} (see Figure~\ref{fig:workload-beat}). One police patrol unit is usually assigned to each beat to patrol that area and respond to 911 calls. 
If the response unit is busy handling another incident, available patrol units in other beats of the same zone will be dispatched to answer calls. Therefore, all patrol units in the same zone can back up for each other, but they typically do not travel across the zone boundaries unless there is a major incident. Under this dispatching policy, each zone can be modeled as a spatially distributed queueing system \citep{Larson1974, larson1981urban}. 

The previous zone and beat configuration in Atlanta (Figure~\ref{fig:workload-beat}) was designed in 2011. Since then, Atlanta has experienced significant population growth and urban development. The U.S.\ Census Bureau estimated that the population of Atlanta has increased by 15.8\% from 2010 to 2017 \citep{UScensusbureau}. 
The population growth led to an increase in police workload, which was exacerbated by the difficulty faced by the APD in officer recruitment and retention. In 2018, the average response time to high-priority 911 calls (e.g., violent crimes) in Atlanta was 9.5 minutes, which was above the national average \citep{Fritz2018}. Moreover, changing demographics and traffic patterns in Atlanta created uneven police workloads among different regions. The colors in Figure~\ref{fig:workload-beat} indicate the average workload in each beat, which we compiled from the 911 police incident report data in 2017. A higher workload in a beat often led to longer 911 call response time. For example, Zone 2, which is located in North Atlanta, had a higher than average workload due to its recent commercial development. From the police report data, we found that the average response time for high-priority calls (e.g., carjacks or burglaries) was 13.5 minutes in Zone 2, whereas the average response time citywide for high-priority calls was 9.5 minutes. During the same period, Zone 2, a historically low-crime area, has experienced a rising number of car thefts, burglaries, and armed robberies \citep{HabershamB2019}.  

The challenge faced by the Atlanta Police Department motivated us to develop a rigorous quantitative method using operations research for redesigning police zones. The general methods we proposed can also be applied to other urban regions facing similar problems.
 
\section{Background of Police Operations}
\label{sec:background}

Police departments usually divide the geographical areas of a city into several \emph{zones}, and furthermore each zone is comprised of several \emph{beats}.
For example, Atlanta is currently divided into 6 zones and 81 beats \citep{APD2020}.
New York City, which has the largest police department in the U.S., contains 77 precincts (comparable to zones) and 302 sectors (comparable to beats) \citep{NYPD2020}. Many cities, including Atlanta and New York City, have fixed patrol districts, although in some other cities the districting design may be changed during a 24-hour period if there are significant changes in the pattern of calls throughout the day. We assume the zone districting is fixed over time in this paper.

Typically, one police patrol unit (e.g., a patrol car) is assigned to each beat, where the unit has primary responsibility. When a patrol unit is not busy serving any active calls, it traverses its home beat to perform preventive patrol. If an emergency call is received, the dispatcher will try to send an available unit in the zone to the location of the reported incident. As a result, a patrol unit may be dispatched to any location within its zone, possibly \emph{outside} its home beat. 
If all units in that zone are busy, the call must wait in a queue until a patrol unit becomes available later. We treat each zone as an independent system since in practice, there are very few dispatches across zone boundaries, due to both long travel time and administrative difficulty \citep{larson1972urban}.
In line with practice, we also make several assumptions on the dispatching rule summarized in Appendix~A.

\begin{figure}[!htb]
\centering
\includegraphics[width=.7\linewidth]{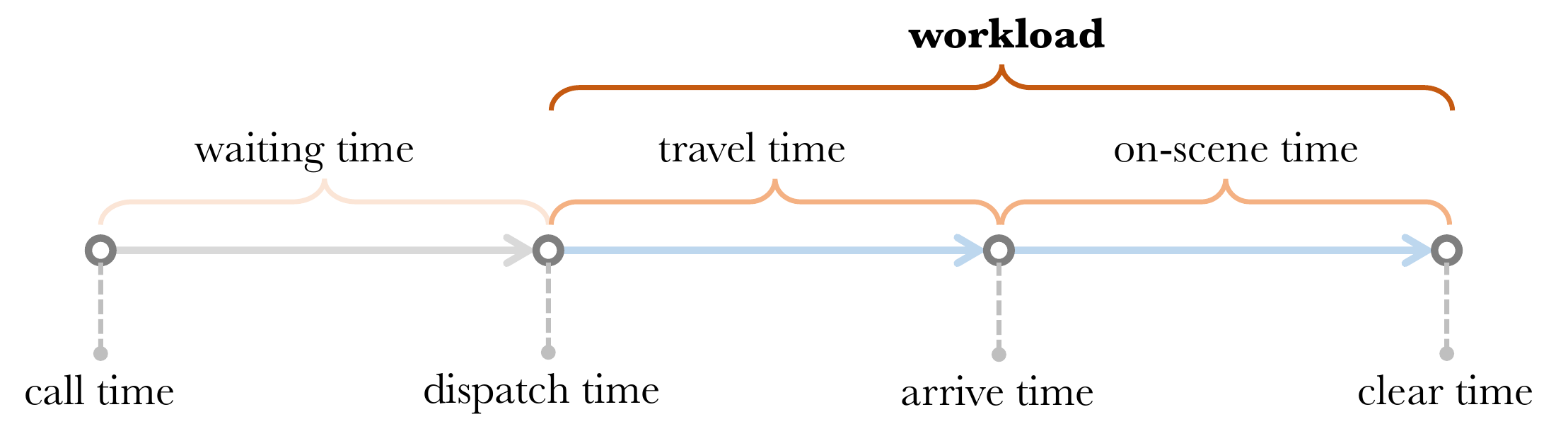}
\caption{The timeline to process a police incident.}
\label{fig:workload-def}
\end{figure}

The typical timeline of a reported incident is as follows (see an illustration in Figure~\ref{fig:workload-def}):
When a 911 call arrives, the dispatcher will try to assign an available patrol unit in the zone to process that incident. If all patrol units are busy, 
the request needs to wait in a queue until at least one unit in that zone becomes available. We refer to the duration between the call arrival time and the dispatch time as the \emph{waiting time}.
Once dispatched, the patrol unit travels from its current location to the location of the reported incident.
The \emph{response time} is thus the sum of the waiting time and the travel time.
After arriving at the scene, the amount of time spent on processing the incident is referred to as the \emph{on-scene time}.
The \emph{service time} for the patrol units to process that incident, which we use to calculate the workload, is equal to the sum of the \emph{travel time} and the \emph{on-scene time}.


\section{Relation to Existing Literature}
\label{sec:related-work}

The study of police zone design and police patrol operations has an extensive history dating back to the 1960s.
The books by \cite{larson1972urban, rios2020optimal} and the survey papers by \cite{chaiken1972methods,green2004anniversary} provide a comprehensive summary of this field.
Below, we review the related literature from three aspects: crime prediction, stochastic modeling of police operations, and optimization for police patrol zones.

\emph{Prediction.}
Big data and analytics techniques have gain popularity in policing over the last two decades.
Several large police departments in the U.S.\ have been experimenting with predictive policing that uses historical crime data to predict future crime activities \citep{perry2013predictive}.
For example, \cite{levine2017new} reported the implementation of a predictive policing system by the New York City Police Department.
However, a downside of predictive policing is that it does not \emph{prescribe} how predictive crime forecasts should be translated into police's operational decisions.
\cite{saunders2016predictions} analyzed the predictive policing program of the Chicago Police Department using a quasi-experiment and found that the effect of predictive policing was not statistically significant.
The practice of predictive policing has also drawn criticisms from civil rights groups for proliferating racial profiling \citep{edwards2016}.

\emph{Stochastic Modeling.}
Police emergency response systems can be viewed as queueing systems where servers and customers are distributed in space, and servers must travel to customer locations.
\cite{Larson1974} proposed a hypercube queueing model for urban emergency services. This hypercube queue model is also studied in \cite{chelst1979estimating,larson1981urban}.
\cite{bammi1975allocation} considered a beat allocation model where the beat design is fixed, but the patrol units can be moved or shared among different beats.

\emph{Optimization.}
\cite{gass1968division} is one of the earliest works that study optimal police beat allocation using integer programming. 
\cite{mitchell1972optimal} studied an optimal beat selection problem.
\cite{bodily1978police} considered fairness issues of police zone design.
\cite{benveniste1985solving} proposed a stochastic optimization model that combines zoning and facility location.
More recently, 
\cite{curtin2005integrating, curtin2010determining} proposed maximal covering models for police patrol area design. 
\cite{cheung2015optimization, chow2015optimization} considered facility location models for police operations. 
\cite{bucarey2015shape} formulated the police districting problem as an enriched $p$-median model.
\cite{camacho2015multi, liberatore2016comparison} introduced a multi-criteria police districting problem that provides a balance between efficiency and workload homogeneity.
\cite{piyadasun2017rationalizing} considered redesigning police patrol beats by clustering methods.
\cite{chen2019designing} proposed a street-level design for urban police districting. Implementations of police zone redesign have been reported for Boston, MA \citep{larson1974illustrative}, Buffalo, NY \citep{sarac1999reconfiguring}, and
Tucson, AZ \citep{kistler2009tucson}.

Unlike most of the existing literature that focuses on one of the three aspects above, our paper presents a comprehensive framework that \emph{integrates} prediction, stochastic modeling, and optimization for the police zone design problem.
We apply this framework to police zone redesign in Atlanta and evaluate the effect of the redesign.

Aside from the police zone design problem, other types of geographical districting problems have also been studied in the operations research literature, including \emph{political districting} \citep{hess1965nonpartisan,GaNe1970} and \emph{school districting} \citep{ferland1990decision}. 
Several papers \citep{WeHe2963, Mi1967, Mo1973, Mo1976, Vi1961, Steven2002} apply meta-heuristics for geographical districting. These studies consider design criteria such as contiguity \citep{Gr1985, Mi1967, GaNe1970, Na1972, Me1972, Vi1961} and compactness \citep{GaNe1970, Ni1990, Yo1988}. In this paper, we use a network flow based formulation by \cite{Shirabe2009} to impose contiguity constraints.
However, the political districting problem and the school districting problem have different objectives than police districting. Political districting and school districting are often formulated as deterministic optimization problems, but police districting must take into under uncertainty associated with crime occurrence and police service time. 
Specifically, our model considers a police patrol system that has multiple moving servers (i.e., patrol units). The workload for each server is stochastic, and the servers' status depends on each other. In addition, changes to the zone configuration will impact the operational dynamics of the zone, including the travel time of servers and the number of 911 calls that need to be served. 

\section{Proposed Method: Data Driven Police Zone Design}

\begin{figure}[!thb]
\centering
\includegraphics[width=.9\linewidth]{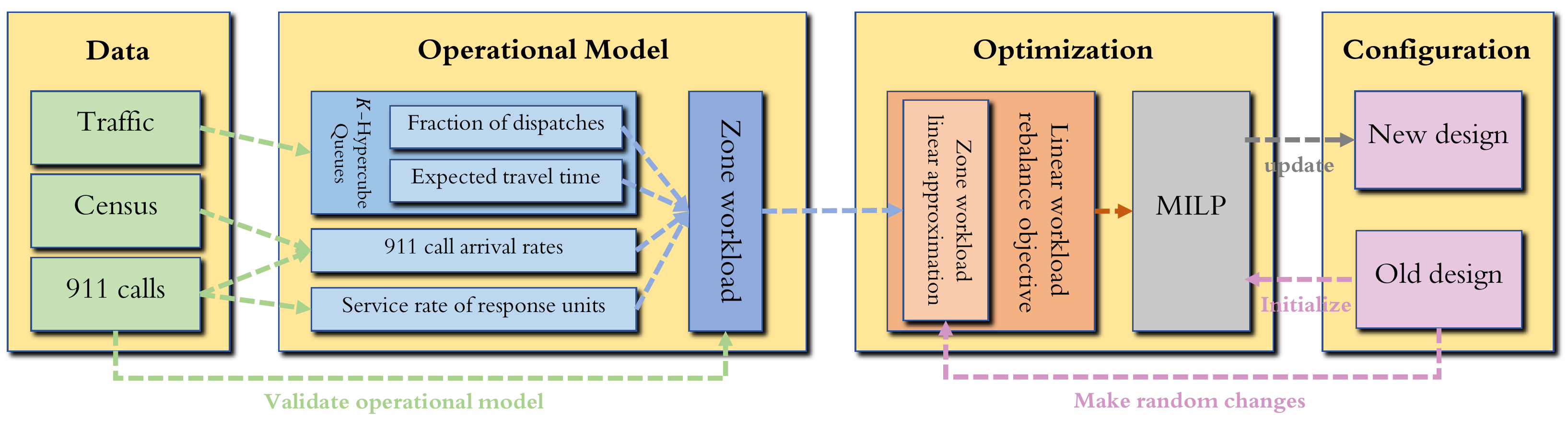}
\caption{An illustration for the data-driven optimization framework of police zone redesign.}
\label{fig:workflow}
\end{figure}

We propose a data-driven modeling and optimization framework for redesigning police patrol zones. This framework integrates statistical prediction for emergency call arrivals, stochastic modeling of police dispatching process, and optimization for zone redesign in an end-to-end design process.
The methodological contribution of this paper is threefold: 
(a) We create a stochastic queueing system for police operations in Atlanta; the system consists of six interdependent queueing models (each model corresponds to one particular zone), which captures the intra-zone dynamics and the inter-zone effects when the zone design is changed.
(b) 
The parameters of the stochastic system are estimated using real police data; the arrival rates of 911 calls across the city are modeled by a spatio-temporal process.
(c) To optimize the zone design using this stochastic model, we overcome the computational challenge by approximating the objective using simulation optimization and finding locally optimal solutions by heuristic search.
A summary of this framework is illustrated in Figure~\ref{fig:workflow}. Next, we describe each component of the framework.


\subsection{Module I. Stochastic Model for Police Patrol and Emergency Response}
The basis of our design process is a stochastic model for the police's response to 911 calls and processing time. 
For a \emph{fixed} zone configuration, the dynamics of patrol units in each zone are independent, since we assume patrol units do not respond to requests from a different zone than they are assigned.
Within each zone, the movement and workload of patrol units are quite complex, as they are determined jointly by the zone design, the dispatch rule, and the beat-to-beat travel time. 

\begin{figure}[!htb]
\centering
\includegraphics[width=0.5\linewidth]{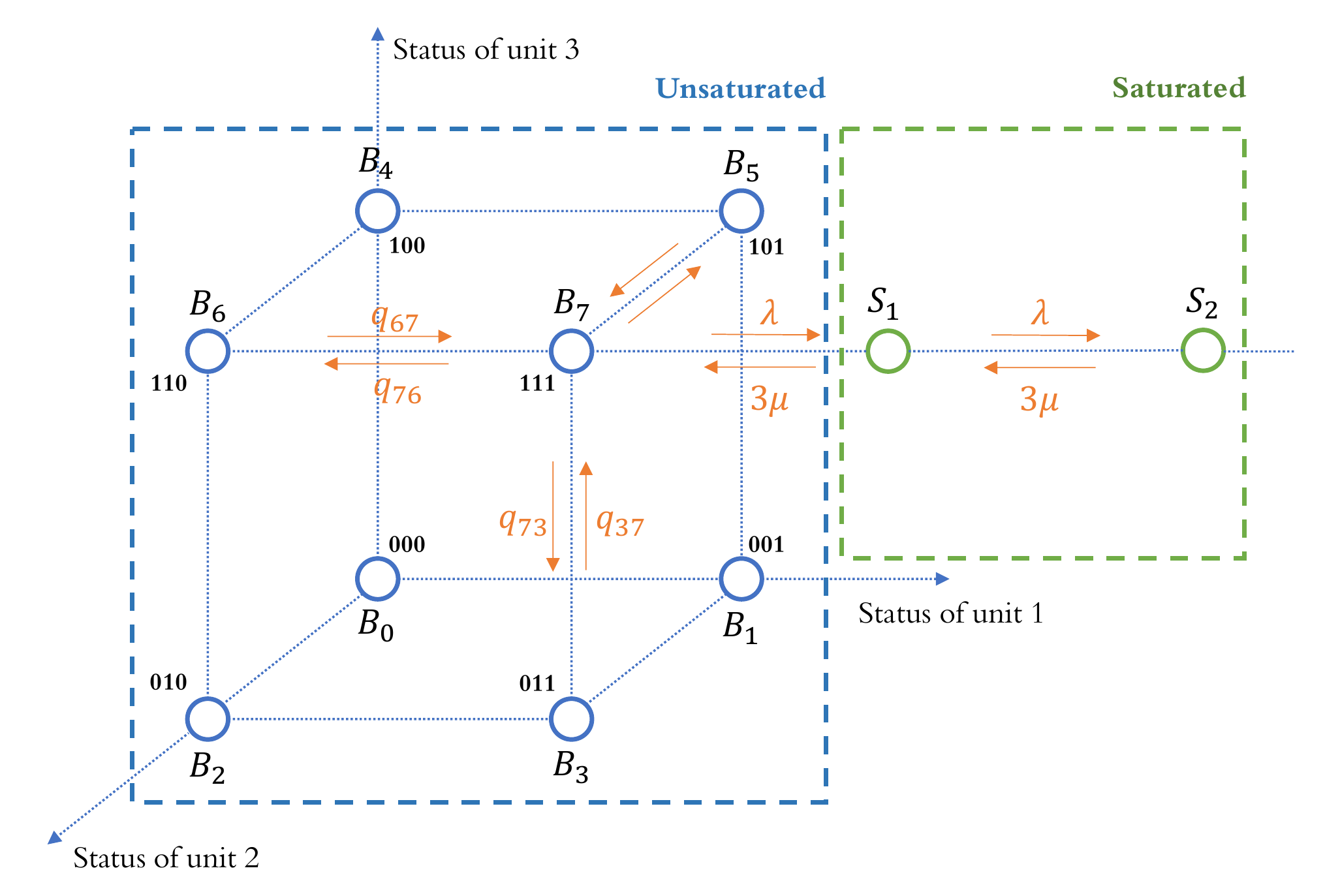}
\caption{The state-space of a hypercube queue model representing a zone with three beats. $B_0, \ldots, B_7$ are unsaturated states and $S_1, S_2\ldots$ are saturated states.
}
\label{fig:hypercube-main}
\end{figure}

We model emergency response operations based on the \emph{hypercube queue model} \citep{Larson1974}; see Appendix~B for details.
The system state depends on the status of all the patrol units in this zone and the number of calls in the queue to be processed. 
Consider the state as a vector of binary numbers, where each element of the vector represents the status of the corresponding unit: the unit is busy processing a call if the corresponding digit of the state equals 1 and otherwise 0. 
Based on whether there is any available unit, the state space can be divided into two parts: (a) \emph{unsaturated states}: the states where the queue for unprocessed calls is empty. 
These states can be represented by a hypercube, where each vertex corresponds to a state. 
(b) \emph{saturated states}: the states where the queue for unprocessed calls is non-empty. By definition, all patrol units are busy under saturated states. 
An example with three patrol units is shown in Figure~\ref{fig:hypercube-main}.

We define the system dynamics using a transition rate matrix based on call arrival rates, travel time, and service rates.
There are two classes of transitions on the hypercube:
\emph{upward transitions} that change a unit's status from idle to busy; 
\emph{downward transitions} that do the reverse. 
The downward transition rate from one state to its adjacent states is always the service rate.
The upward transition rates, however, will depend on the dispatch rule when a call is received and the travel time between beats.
We assume that when a call is ready to be processed, the nearest available unit with the minimum mean travel time to the location of the call is dispatched. 

The steady-state probabilities of this system define the likelihood that each state is occupied in the long run, which can be determined by solving the balance equations. 
We note that although the number of saturated states is infinite, it can be expressed in closed-form as in an M/M/c queue.
The main challenge for computing the steady-state distribution is that the number of \emph{unsaturated} states on the hypercube grows exponentially with the number of beats. 
We developed an efficient computation method to approximate the steady-state distribution by exploiting the fact that the transition rate matrix is sparse.



\subsection{Module II. Model Estimation and Prediction}

We estimate the parameters of the hypercube queue model  using data provided by the Atlanta Police Department (APD) and public data sources. We first describe the datasets that we used for estimation. Then we describe our call arrival rate and workload prediction method. 

\subsubsection{Data Sets.} 
We have utilized two key data sets for police workload prediction: 
\begin{enumerate}
    \item \emph{911 Call Data}: the APD has electronic records of police reports, where each record consists of detailed information about a 911 call, such as the call time, dispatch time, arrive time, clear time, and location. 
    These reports also provided travel information of each dispatch, including a departure time and location, and an arrival time and location.
    Using this data set, we estimate the average travel time between beats, the 911 call arrival rates for each beat, and the service rate of patrol units. 
    \item \emph{American Community Survey}: the U.S.\ Census Bureau provides an annual American Community Survey (ACS), which includes comprehensive information about the population, demographic, and economic status of different areas in Atlanta. 
    We use {population}, {number of housing unit}, {school enrollment}, {median household income}, {median number of rooms}, {media age}, {median house price}, and {average year built} from the ACS as explanatory variables in our predictive model. 
\end{enumerate}

\subsubsection{Call Arrival Rate Prediction.}
We cannot use standard stochastic processes to model
911 call arrivals, because the arrival rates have a significant seasonality pattern and yearly trend; the rate also have a strong spatial correlation among adjacent geographical areas.
Therefore, we propose a spatio-temporal model to predict call arrival rates using spatially lagged endogenous regressors \citep{Rosen1974}.
We assume that the 911 call arrival rate of a beat depends on (a) the arrival rates of adjacent beats, (b) the arrival rate of the same beat in the previous year, and (c) the demographic factors in the past five years. 
We also assume that the noise term has a spatial structure, where the covariance between two beats is determined by an exponential kernel function. 
The parameters of the model are fitted by maximum likelihood estimation. See Appendix~C for details.

\subsubsection{On-scene Time Estimation.}
As we assume the on-scene time of a response unit is independent of its service region, we analyze the distribution of average on-scene time for all units in the same zone. We have observed from data that the on-scene time of response units fits to an exponential distribution shown in Figure~\ref{fig:on-scene}. 

\subsubsection{Travel Time Estimation.}
\begin{figure}
\centering
\begin{minipage}{.45\textwidth}
  \centering
  \includegraphics[width=.8\linewidth]{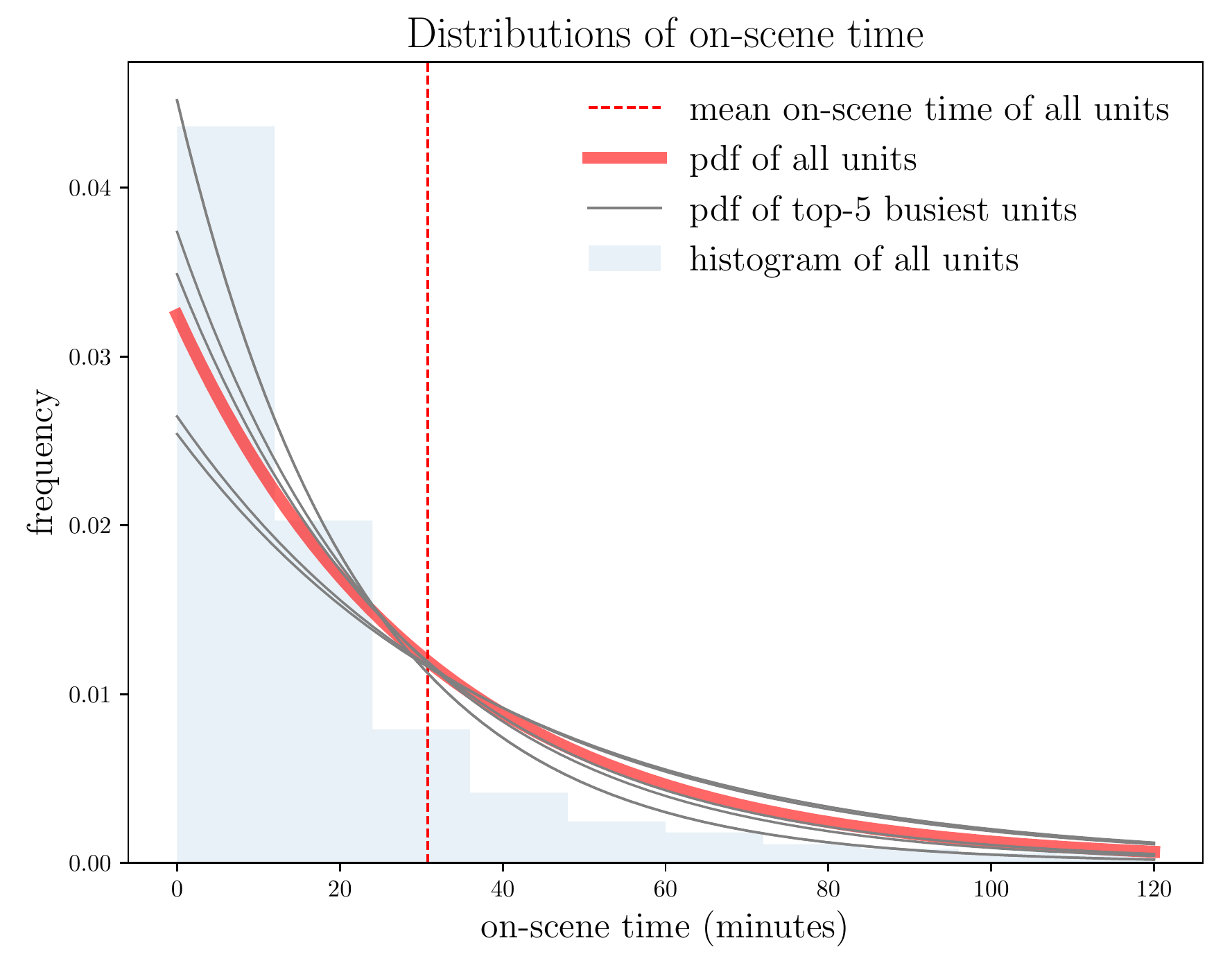}
  \captionof{figure}{Distribution of the on-scene time of all units: The on-scene time data follows an exponential distribution, shown by the red and grey lines. The red dash line indicates the mean on-scene time for all response units.}
  \label{fig:on-scene}
\end{minipage}%
\hspace{.2in}
\begin{minipage}{.45\textwidth}
  \centering
  \includegraphics[width=.71\linewidth]{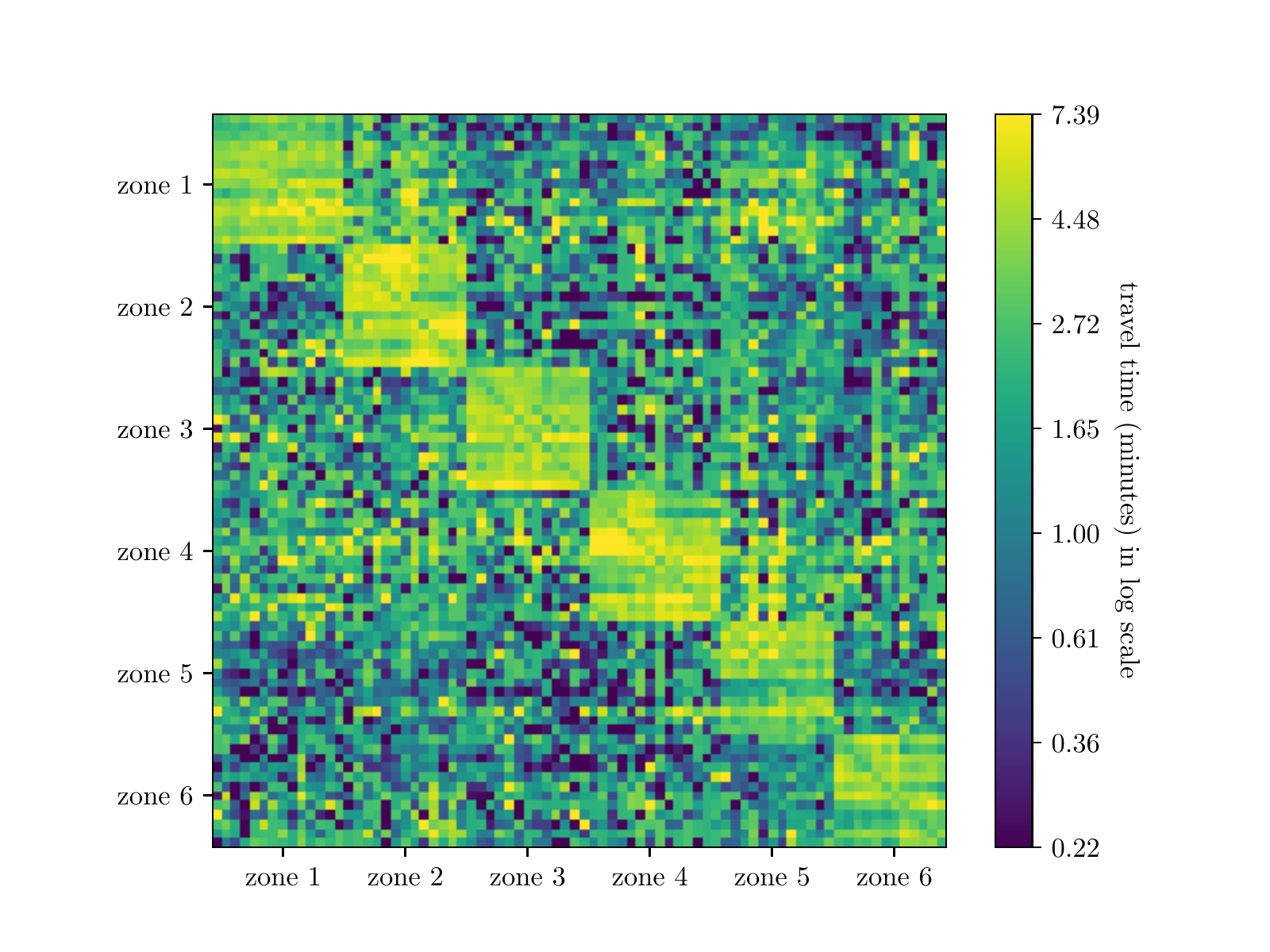}
  \captionof{figure}{The \emph{beat-to-beat} matrix (in seconds) estimated from the real data. Each row and each column corresponds to a beat, and each small square represents average beat-to-beat travel time.
  }
  \label{fig:tau}
\end{minipage}
\vspace{-.15in}
\end{figure}

We estimated the travel time between beats from the 911 call dispatch data. The estimation result is shown in Figure~\ref{fig:tau}.  

\subsubsection{Performance Metrics.}
Using the steady-state distribution of the hypercube queue model with estimated model parameters, we can derive several key performance measures (See Appendix~B for details), which will guide our zone design decisions in the next step: 
\begin{itemize}
\item \emph{Expected Travel Time.} 
The travel time for a dispatch is affected by both the beat-to-beat travel time and the fraction of dispatches. 
The fraction of dispatches critically depends on the utilization factor of the hypercube queue system. If the utilization is low, most calls can be processed by units close to their home beats, resulting in a shorter travel time. However, if the system is congested, calls must wait in the queue, which will be processed by the first unit available.
This will lead to a longer travel time.
\item \emph{Response Time.}
The response time for a call includes both the travel time for the patrol unit to arrive at the scene and the possible waiting time if the call incurs a positive queue delay.
\item \emph{Patrol Unit Workload.}
The workload of one patrol unit is the sum of the steady-state probabilities for all the states where this unit is busy.
The total workload for all patrol units in the zone is simply the aggregation of individual unit workloads.
\end{itemize}

\subsubsection{Model Validation.}

\begin{figure}[!t]
\centering
\begin{subfigure}[h]{0.4\linewidth}
\includegraphics[width=\linewidth]{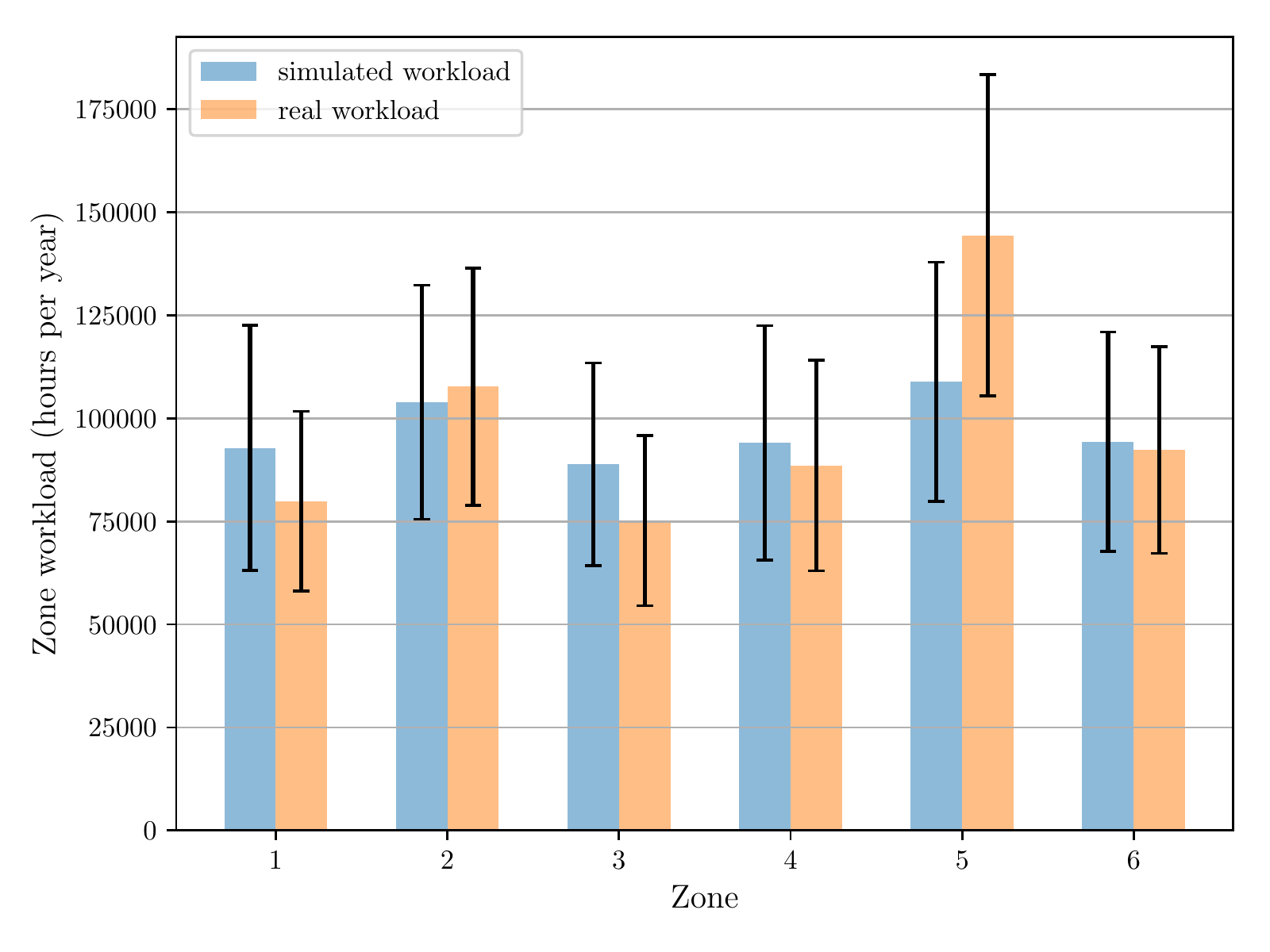}
\caption{Comparisons between zones}
\label{fig:zone-workload-by-zone}
\end{subfigure}
\begin{subfigure}[h]{0.4\linewidth}
\includegraphics[width=\linewidth]{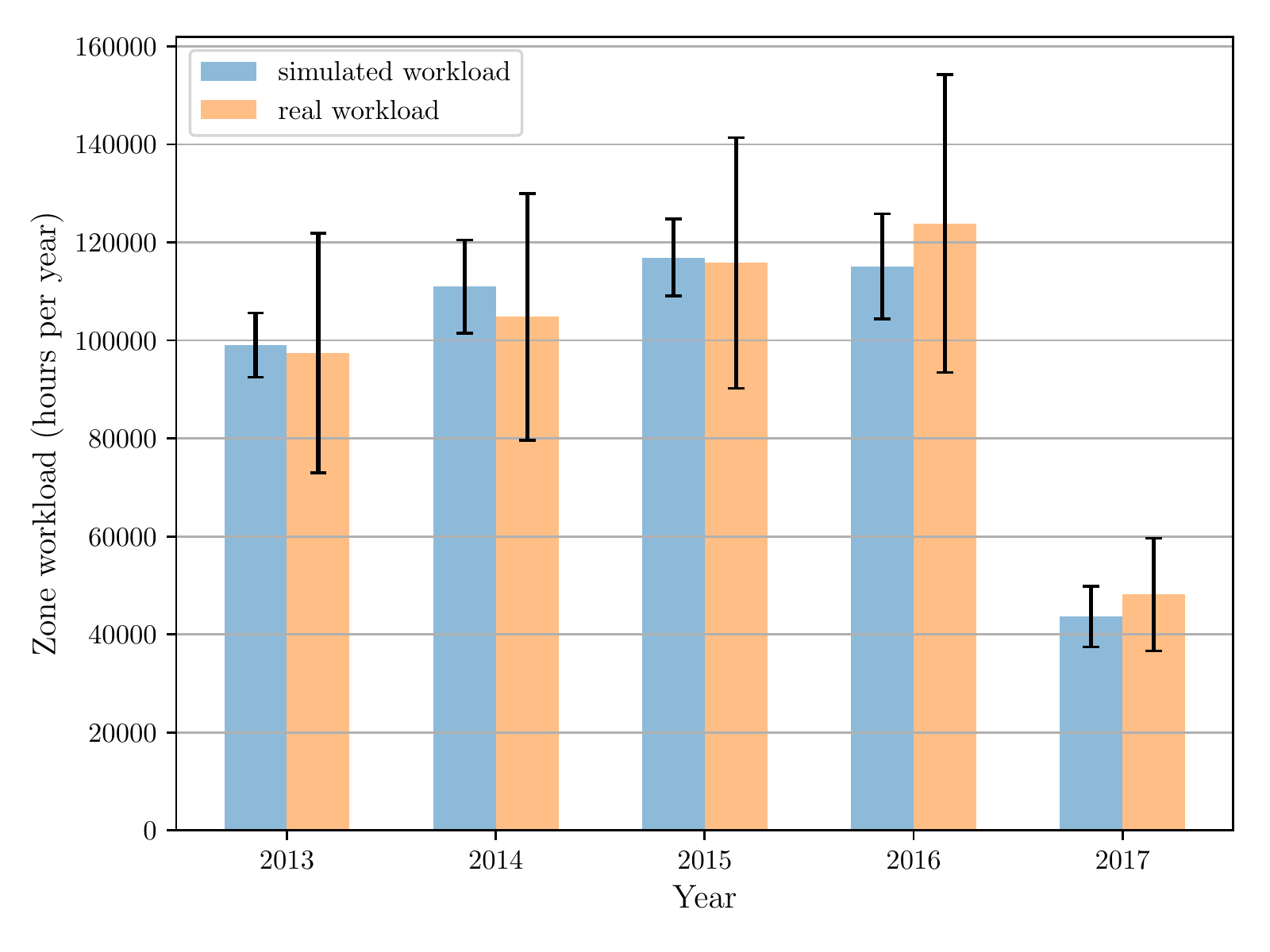}
\caption{Comparisons between years}
\label{fig:zone-workload-by-year}
\end{subfigure}
\caption{Validations for the operations model of annual zone workload (only first quarter data are available in year 2017) (\ref{fig:zone-workload-by-zone}). shows the averaged simulated workload for each zone and (\ref{fig:zone-workload-by-year}) shows the averaged simulated workload for each year.}
\label{fig:zone-workload-validation}
\end{figure}

To evaluate the predictive accuracy of the estimated police zone workload, we compare the real zone workload reported by historical data with the out-of-sample zone workload simulation from our operations model using leave-one-out cross-validation.
Specifically, the operations model 
takes the 911 calls arrival rates, the police service rate, and beat-to-beat travel time in a particular year as input variables and generate the simulated police workload for each zone using the existing zone design.
%
In practice, we estimate these input quantities using the data 
from the past five years. 
Figure~\ref{fig:zone-workload-validation} shows comparisons for the average workload between zones and averaged workload between years. 
To carry out the evaluation for each zone/year in Figure~\ref{fig:zone-workload-validation}, the model estimation is performed using all data except for one zone/year, and a prediction is made for that zone/year.
The result shows that the simulation is consistent with the real data, and the percentage difference is less than 6.7\%.

The operations model allows us to quantify the impact of an arbitrary zone design on the zone workload analytically. However, it cannot be used to search for an improved zone design efficiently due to the randomness and the non-linear structures of zone workloads with respect to decisions. 
In the next section, 
we will develop a way to approximate the zone workload, which can be used in the optimization.

\subsection{Module III. Zone Design Optimization}

We develop an optimization model for zone redesign. We assume that the shape of the beats are fixed, so the zone design problem can be formulated as a graph partition problem, where we allocate the beats to a pre-determined number of zones. 

Different objectives have been studied in the literature for police zone design (see Appendix~D).
After several discussions with the Atlanta Police Department, we learned that their primary concern is the workload imbalance among different geographical areas, as they have observed a trend of police overwork in some zones. 
In particular, longer shifts and police officer fatigue may lead to a host of public safety issues and affect the police department's budget (in the form of bonuses and overtime).
Therefore, in this paper, we follow the zone workload definition introduced by \cite{mitchell1972optimal} 
and choose the objective function as the workload variance among different zones, which quantitatively measures the police workload imbalance among zones.
This metric was also used internally by APD to measure previous zone redistricting plans \citep{egbert2016}. 
The optimization model aims to minimize the workload imbalance across zones subject to some shape constraints for each zone, such as contiguity and compactness.
The full optimization model is included in Appendix~E.



\subsubsection{Zone-level Workload Approximation.}
A key challenge for solving the optimization problem is characterizing the complex dependence between the zone design and the zone-level workload. In theory, the expected workload in each zone can be computed using the hypercube queue model given any zone design.
However, it is impractical to carry out such computation to search for an improved zone design since the decision space is extremely large. For example, with the six zones and 81 beats in Atlanta, this will generate more than $2.4 \times 10^{59}$ possible zone designs, where each zone design corresponds to a separate queueing model.

To tackle this challenge, we consider a local linear approximation of the workload function and
then use the approximated value in the objective function of the optimization model,
which allows the objective function of the optimization problem to be expressed explicitly without repeatedly computing the hypercube queue model.
This linear approximation represents the first-order Taylor expansion of the workload change when the current zone design is modified slightly by adding a few beats or removing a few beats.
In Appendix~E, we discuss the detailed process of generating this linear approximation and validating the effectiveness of our approximation method 
using randomly selected {out-of-sample} designs. 

\vspace{-3pt}
\subsubsection{Contiguity and Compactness Constraints.}
In addition to balancing the police workload, the shapes of zones should be contiguous and compact. 
Since patrol units almost always travel within a zone,
a narrow or snakelike shape increases the travel time and reduces the efficiency of patrol operations and emergency response.
Therefore, we formulate the \emph{contiguity} and \emph{compactness} criteria using binary variables and a set of linear constraints,
which requires that beats in the same zone to be geographically connected and close to each other; see Appendix~E for details. 

\vspace{-3pt}
\subsubsection{MILP Formulation and Local Search.}
With the linear approximation of the zone workload, the objective function measuring the workload variance becomes a \emph{quadratic} function of the binary decision variables.
Using the McCormick envelopes technique in integer programming \citep{mccormick1976computability},
this quadratic function can be expressed equivalently as a linear form.
The reformulated mixed-integer linear programming (MILP) model contains more than 240,000 variables and has a huge solution space with more than $2.4\times 10^{59}$ possible solutions. 

Based on the feedback from APD, we also consider a practical constraint that beat changes should be as few as possible. A drastic zone design change is undesirable for two reasons: First, a large-scale operational change will result in high implementation costs; second, a radical design change will usually face significant uncertainties and unpredictable risks in the future operation. 
For this reason,  we only consider plans with no more than six beat changes.  
Moreover, two specific design requirements were suggested by APD and included as constraints in our model: (a) beat 305, 111 should be staying in zone 3 and 1, respectively; (b) zone 3 will not take any existing beat from zone 4.

In order to find a near optimal solution to the MILP efficiently and avoid limit the number of beat changes from the existing design, we use a local search method based on simulated annealing.
The optimization process takes about 30 to 50 minutes 
on a standard laptop with a quad-core 4.7 GHz processor.



\section{Implementation Results}
\label{sec:results}

During a one-year collaboration between Georgia Tech and APD, we used our data-driven zone design framework as a decision support tool and implemented a  
zone reconfiguration plan. 
In the following, we first discuss the police workload imbalance analysis using real and predicted results; then we present our zone redistricting plan and the major changes we made to the original zone configuration; we also discuss the implementation's impact; lastly, we investigate the effectiveness of the proposed plan by analyzing the post-implementation data.

\subsection{Workload Imbalance Analysis}
\label{sec:workload-analysis}

Under the previous zone design in Atlanta, the workload imbalance across zones had worsened in the past few years. This imbalance is partly due to the growth and redistribution of population, as well as changing (mostly worsening) traffic conditions. Here, we analyze how zone configuration affects the workload imbalance over the years, and predict the future trend if there is no change in the configuration.
To quantify the imbalance of workload, we use the variance of workload across zones as the metric.

\begin{figure}[!h]
  \centering
  \begin{subfigure}[b]{0.24\textwidth}
      \centering
      \includegraphics[width=\textwidth]{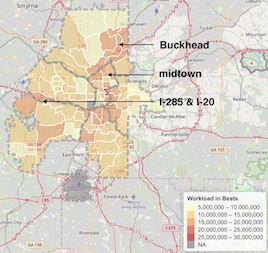}
      \caption{2014}
      \label{fig:beat-2014}
  \end{subfigure}
  \begin{subfigure}[b]{0.24\textwidth}
      \centering
      \includegraphics[width=\textwidth]{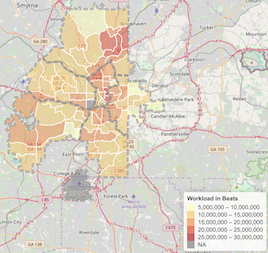}
      \caption{2015}
      \label{fig:beat-2015}
  \end{subfigure}
  \begin{subfigure}[b]{0.24\textwidth}  
      \centering 
      \includegraphics[width=\textwidth]{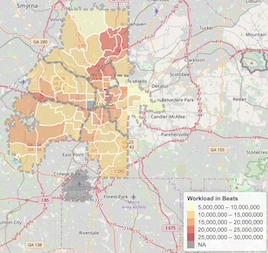}
      \caption{2016}
      \label{fig:beat-2016}
  \end{subfigure}
\vfill
  \begin{subfigure}[b]{0.24\textwidth}   
      \centering 
      \includegraphics[width=\textwidth]{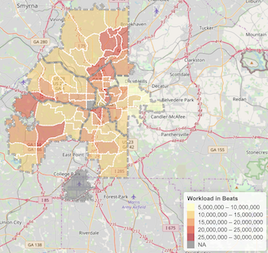}
      \caption{2017}
      \label{fig:beat-2017}
  \end{subfigure}
  \begin{subfigure}[b]{0.24\textwidth}
      \centering
      \includegraphics[width=\textwidth]{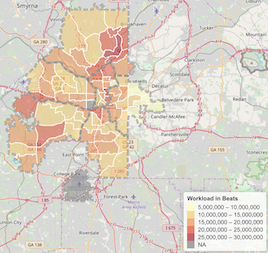}
      \caption{2018 *}
      \label{fig:beat-2018}
  \end{subfigure}
  \begin{subfigure}[b]{0.24\textwidth}   
      \centering 
      \includegraphics[width=\textwidth]{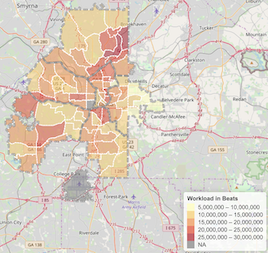}
      \caption{2019 *}
      \label{fig:beat-2019}
  \end{subfigure}
  \caption{The beat level workload (in seconds) in Atlanta from the year 2014 to 2019 (2018 and 2019 levels are predicted by the hypercube queue model).  Dashed grey lines represent the boundaries of the six zones. Darker color means a higher workload. (The gray area at the bottom corresponds to the Atlanta International Airport, which we exclude from our analysis.)}  
  \label{fig:beat-change}
\end{figure}

We compiled the actual workload from the year 2014 to 2017 (the most recent year at the time of our redesign process) from the police report data, and computed the predicted workload for the year 2018 and 2019 using the stochastic queue model we developed earlier.
The workload in each beat is displayed on the city map in Figure~\ref{fig:beat-change}. One can observe a clear trend that both the mean and the variance of workload were growing up to 2017, and would continue to grow according to our model prediction.
Moreover, the workload was exceptionally high in several concentrated areas, such as Buckhead in Zone 2, Midtown in Zone 5, and the intersection of highway I-285 and highway I-20 in Zone 4. 

\begin{figure}[!ht]
\centering
\includegraphics[width=.5\linewidth]{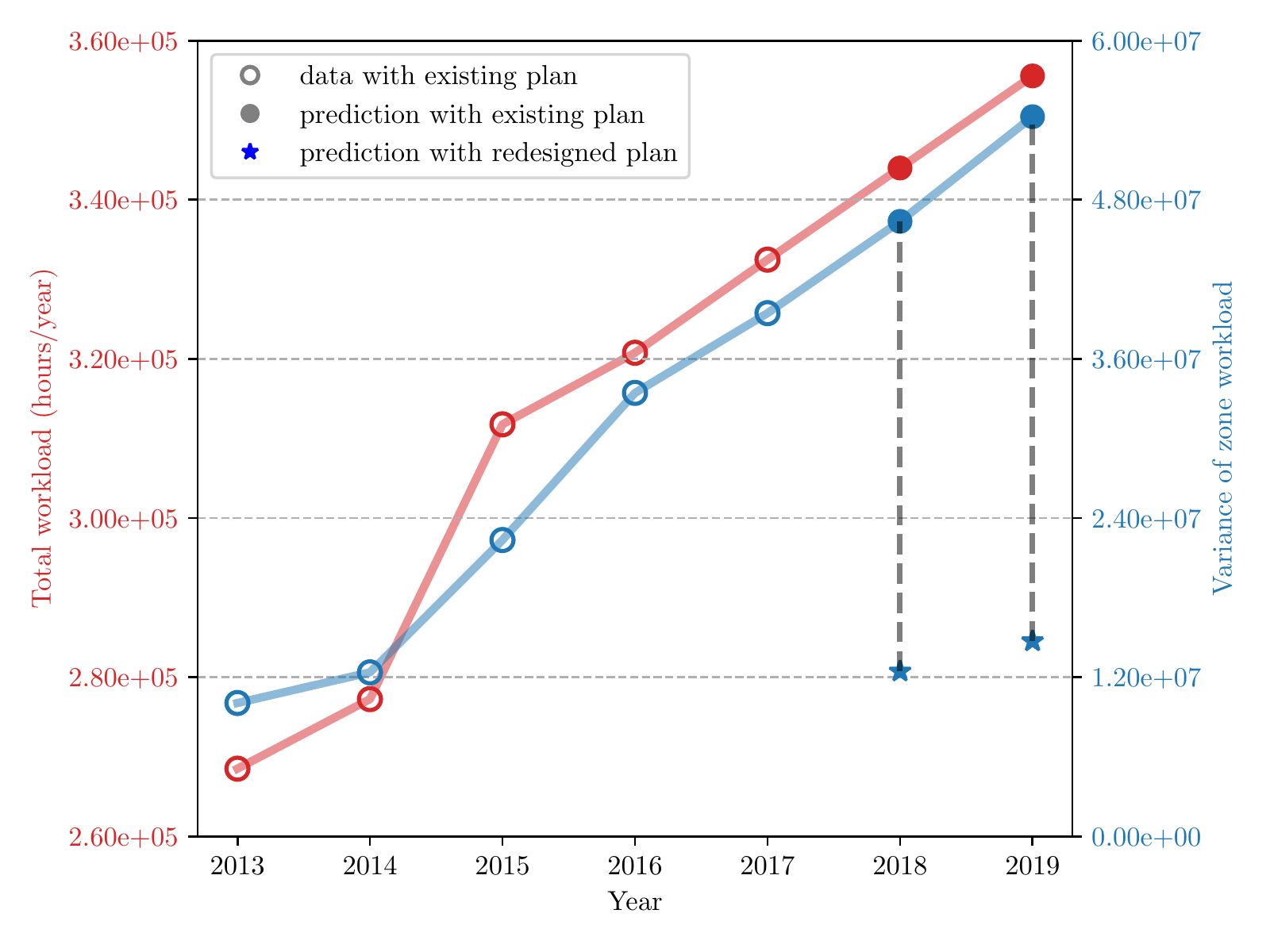}
\caption{The total workload of all zones (red) and variance of workload across zones (blue) from the year 2013 to 2019: The values for 2018 and 2019 are predicted levels assuming the configuration is not changed, and the stars for the year 2018 and 2019 are predicted levels under our proposed redesigned configuration.}
\label{fig:workload-var}
\end{figure}

In Figure~\ref{fig:workload-var}, we plot the total workload of the entire city (red curve) and the workload imbalance across the six zones (blue curve).
As before, the results for 2018 and 2019 are based on predictions assuming the same zone design. The figure shows that both the total workload and the workload variance would have increased without zone reconfiguration.

\subsection{Zone Reconfiguration}
\label{sec:redesign-analysis}

\begin{figure}[!h]
  \centering
  \begin{subfigure}[b]{0.24\textwidth}
      \centering
      \includegraphics[width=\textwidth]{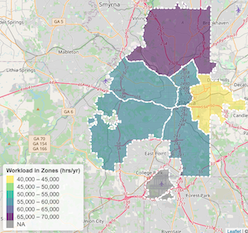}
      \caption{existing 2018}
      \label{fig:old-2018}
  \end{subfigure}
  \begin{subfigure}[b]{0.24\textwidth}  
      \centering 
      \includegraphics[width=\textwidth]{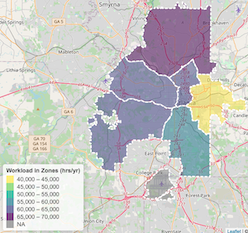}
      \caption{existing 2019}
      \label{fig:new-2018}
  \end{subfigure}
  \begin{subfigure}[b]{0.24\textwidth}   
      \centering 
      \includegraphics[width=\textwidth]{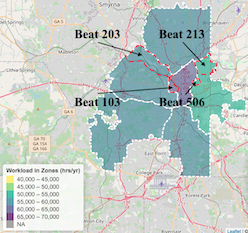}
      \caption{redesigned 2018}
      \label{fig:old-2019}
  \end{subfigure}
  \begin{subfigure}[b]{0.24\textwidth}   
      \centering 
      \includegraphics[width=\textwidth]{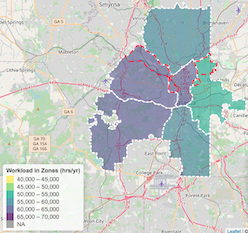}
      \caption{redesigned 2019}
      \label{fig:new-2019}
  \end{subfigure}
  \caption{Comparison of zone workload distribution between the existing plan and the redesigned plan: for 2018 and 2019 (based on prediction), respectively. The color depth of each zone represents the level of its annual workload. Beat shifts in the redesigned plans have been highlighted by red dash lines. As we can see, for the predicted year 2018 and 2019, the redesigned plan achieves a more balanced workload across zones comparing to the existing plan.} 
  \label{fig:design-comparison}
\end{figure}

Using the queueing and optimization models developed earlier in this paper,
we proposed a new police zone plan for Atlanta. 
Our analysis report in 2018 contained beat-wise workload prediction for the next two years (2018 and 2019) and proposed three candidate designs with similar beat shifts that all attain improved workload balance. 
In Table~\ref{tab:config-comparison}, we list the predicted annual workload in each zone, total workload, and workload variance. 
After our proposed plan was reported to the police, several follow-up meetings were held to 
elicit feedback from the APD and
deliberate various design trade-offs. The APD Deputy Chief and other senior police officers also participated in these discussions and provided feedback.
After several rounds of discussions, the police chose one of our proposed designs as the new police zone configuration.
The new design and the previous design are shown in Figure~\ref{fig:design-comparison}.
This particular design was selected for three major reasons: (a) this plan would only change four beats, so it has fewer beat shifts in comparison with other candidate plans; (b) Beat 203 and 213 are two areas with relatively low police workload in zone 2, so moving them out of zone 2 may help focus their police resources on curbing the crime surge in the Buckhead area (the center of zone 2); (c) Zone 3 and Zone 4 have maintained a good operational balance for the recent years; keeping these two zones intact can also help reduce the implementation cost of the plan. 
For both predictions in 2018 and 2019, the redesigned zone plan has a more balanced workload distribution over the six zones (the colors in Figures \ref{fig:design-comparison} (c) and \ref{fig:design-comparison} (d) are more uniform). 
The new design is projected to reduce the workload variance significantly by 62.79\% and 55.97\% in 2018 and 2019, respectively. 




\begin{table}[!h]
\centering
\caption{Comparison of the predicted workload under the previous zone design and the new zone design.}
\label{tab:config-comparison}
\resizebox{\textwidth}{!}{%
\begin{tabular}{lccccccccc}
  \toprule[1pt]\midrule[0.3pt]
 & \multicolumn{7}{c}{Workload ($\times10^4$ hrs/yr)} & \multirow{2}{*}{Variance ($\times10^{7}$)} & \multirow{2}{*}{\begin{tabular}[c]{@{}c@{}}Variance increase ratio\\compare with 2016 (\%)\end{tabular}} \\ 
 & Zone 1 & Zone 2 & Zone 3 & Zone 4 & Zone 5 & Zone 6 & Total &  &  \\ \hline
Real 2016 & 5.3744 & 6.0082 & 5.3891 & 5.6164 & 5.5479 & 4.1413 & 32.0776 & 3.3441 & N/A\\
Real 2017 & 5.6648 & 6.2981 & 5.5154 & 5.7919 & 5.7319 & 4.2454 & 33.2479 & 3.9454 & $+17.98$ \\
Predicted 2018 & 5.9558 & 6.5833 & 5.6366 & 5.9606 & 5.9143 & 4.3488 & 34.3996 & 4.6375 & $+38.67$ \\
Predicted 2019 & 6.2487 & 6.8650 & 5.7642 & 6.1292 & 6.0976 & 4.4518 & 35.5568 & 5.4267 & $+62.27$ \\
Predicted 2018 with redesign & 5.9018 & 5.6961 & 5.6366 & 5.9606 & 6.1595 & 5.0448 & 34.3996 & 1.2442 & $-62.79$ \\
Predicted 2019 with redesign & 6.1833 & 5.9294 & 5.7642 & 6.1292 & 6.3686 & 5.1819 & 35.5568 & 1.4721 & $-55.97$  \\
\midrule[0.3pt]\bottomrule[1pt]
\end{tabular}%
}
\end{table}

\subsection{Implementation and Impact}
\label{sec:implementation}

In January 2019, we submitted the final report for zone redesign to the Atlanta Police Department, which the Police Chief and her team reviewed. Our report analyzed the annual growth trend of police workload and proposed a detailed redistricting plan. 
Our redistricting plan involves changing four zones: Zone 6 in East Atlanta will increase by four square miles. Zone 1 in Northwest Atlanta will increase by two square miles. Zone 2, which covers Northeast Atlanta and Buckhead, decreases by seven square miles. Zone 5 also has some minor changes. The detailed plan is shown in Figure~\ref{fig:design-comparison}. Overall, the redistricting would rebalance the police workload between the four zones and reduced the average response time in Zone 2 \citep{HabershamA2019}. 

In February 2019, Atlanta's City Council voted to approve our proposed redistricting plan. The Atlanta Police Department officially implemented the plan on March 17, 2019. In a statement about this redesign \citep{APD2019}, the APD Deputy Chief Jeff Glazier stated: ``It is important that we examine our officer workload periodically, and with the help of Georgia Tech, we were able to do so in a data-driven manner. We are confident these changes will balance the workload in all zones.'' The new zone design was also covered by several media outlets, including the Atlanta Journal-Constitution, the largest daily newspaper in the metro Atlanta area. On social media, many residents praised the change and thanked the APD and the Georgia Tech team for the new plan \citep{socialmedia}.

\subsection{Post-Implementation Analysis}
\label{sec:post-implementation}

In August 2020, we received a new set of 911 call data provided by the Atlanta Police Department thanks to Major John Quigley's help. 
This dataset contains records from the pre-implementation period (March 2017 $\sim$ March 2019) to the post-implementation period (March 2019 $\sim$ March 2020). Note that we excluded the data after March 2020 because police operations had been drastically affected by the COVID-19 pandemic and the citywide lockdown. 
Also, due to the inconsistency between the new data format and the previous data format, we normalized the new data set by bringing the workload of each zone in the overlapping period (i.e., Year 2017) to the same level.
Using this dataset, we evaluate the zone-level workload distribution and variance before and after the implementation.
Apart from the workload analysis, we also focus on three other metrics as we described in the background section,
namely, waiting time, travel time, and response time.
(We do not study on-scene time since it is exogenous to the zone design.) In particular, response time is important for the perception of police operations efficiency by the general public and directly reflects the effectiveness of the districting plan. 

\begin{table}[!htb]
\centering
\caption{Comparison of the real workload before and after the redesign. 
}
\label{tab:redesign-comparison}
\resizebox{\textwidth}{!}{%
\begin{tabular}{ccccccccccc}
  \toprule[1pt]\midrule[0.3pt]
  \multirow{2}{*}{Time period} & \multicolumn{7}{c}{Workload ($\times10^4$ hrs/yr)} 
  & \multirow{2}{*}{Variance ($\times10^{7}$)} & \multirow{2}{*}{\begin{tabular}[c]{@{}c@{}}Variance increase\\year-on-year ratio (\%)\end{tabular}} \\ 
 & Zone 1 & Zone 2 & Zone 3 & Zone 4 & Zone 5 & Zone 6 & Total 
& &  \\ \hline
Mar 2018 $\sim$ Mar 2019 (before) & 6.0568 & 6.8737 & 5.5450 &  5.7040 & 6.1607 & 4.7360 & 35.0762 
&  4.2375 & $+7.40$\\
Mar 2019 $\sim$ Mar 2020 (after) &
5.5142 & 6.3753 &  5.7157 &  5.0090 & 5.9462 & 5.0129 &
33.5733 
& 2.3925 & $-43.54$\\
\midrule[0.3pt]\bottomrule[1pt]
\end{tabular}%
}
\end{table}

\begin{figure}[!htb]
  \centering
  \begin{subfigure}[b]{0.24\textwidth}  
      \centering 
      \includegraphics[width=\textwidth]{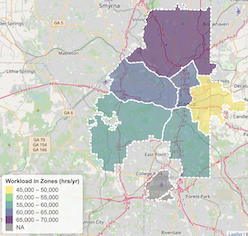}
      \caption{before redesign}
  \end{subfigure}
  \quad
  \begin{subfigure}[b]{0.24\textwidth}   
      \centering 
      \includegraphics[width=\textwidth]{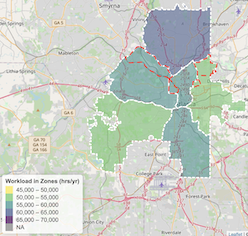}
      \caption{after redesign}
  \end{subfigure}
  \caption{
  Zone-level average workloads before and after the redesign.
  The beats changed in the redesigned plans are highlighted by red dash lines. 
  }
  \label{fig:design-comparison-newdata}
\end{figure}

In Table~\ref{tab:redesign-comparison} and Figure~\ref{fig:design-comparison-newdata}, we show the average workload of each zone and their variances before and after the redesign in March 2019. 
It is evident that our redesign successfully mitigated the workload increase rate and reduced the deterioration of workload imbalance: the workload variance of 2019 decreased by 43.54\% comparing to 2018.
For Zone 2, a region with rapidly increasing crime rates and a focus of this redesign, the analysis shows that our new plan has reduced its police workload by 7.25\% compared to 2018.

\begin{figure}[!h]
  \centering
  \begin{subfigure}[b]{0.25\textwidth}
      \centering
      \includegraphics[width=\textwidth]{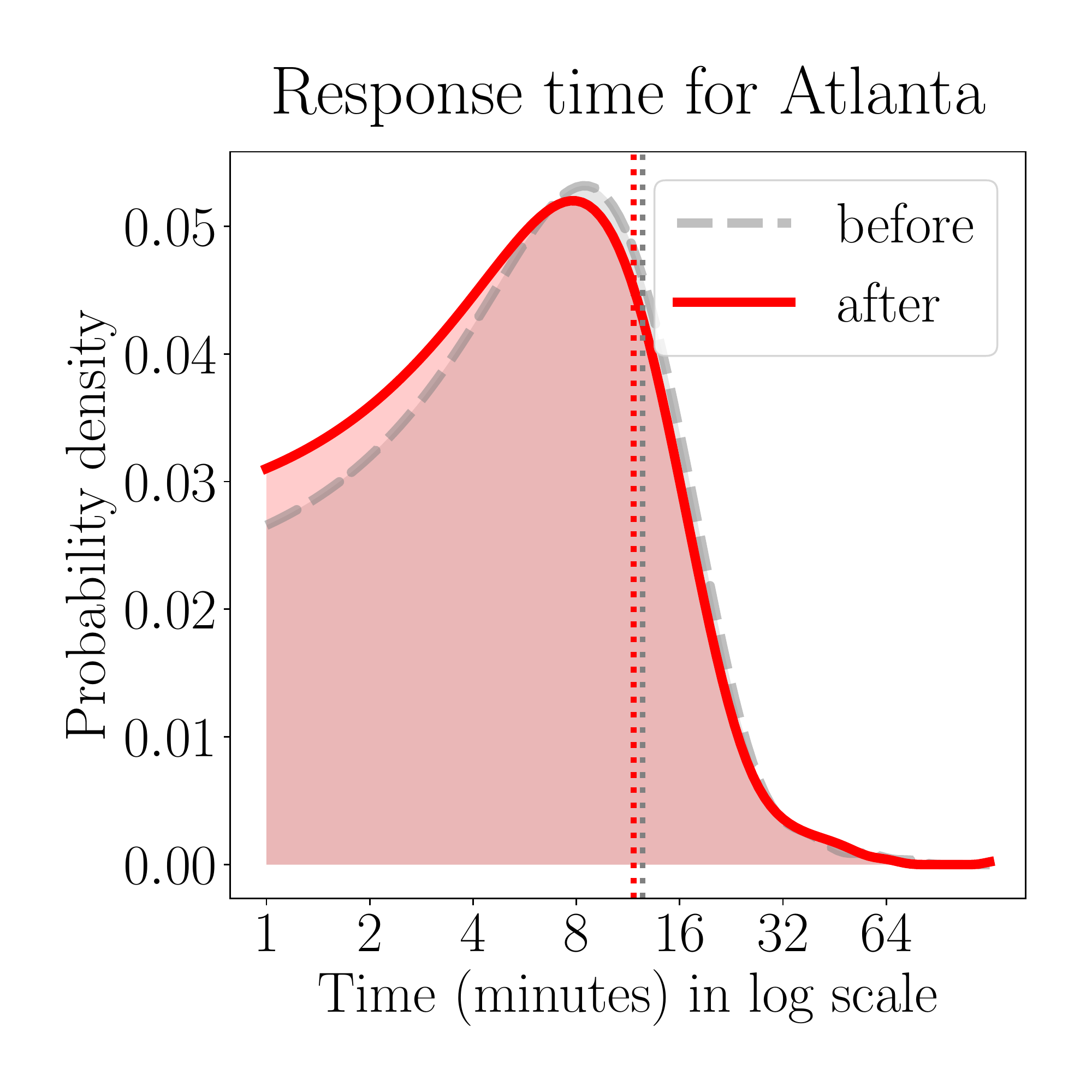}
  \end{subfigure}
  \begin{subfigure}[b]{0.25\textwidth}  
      \centering 
      \includegraphics[width=\textwidth]{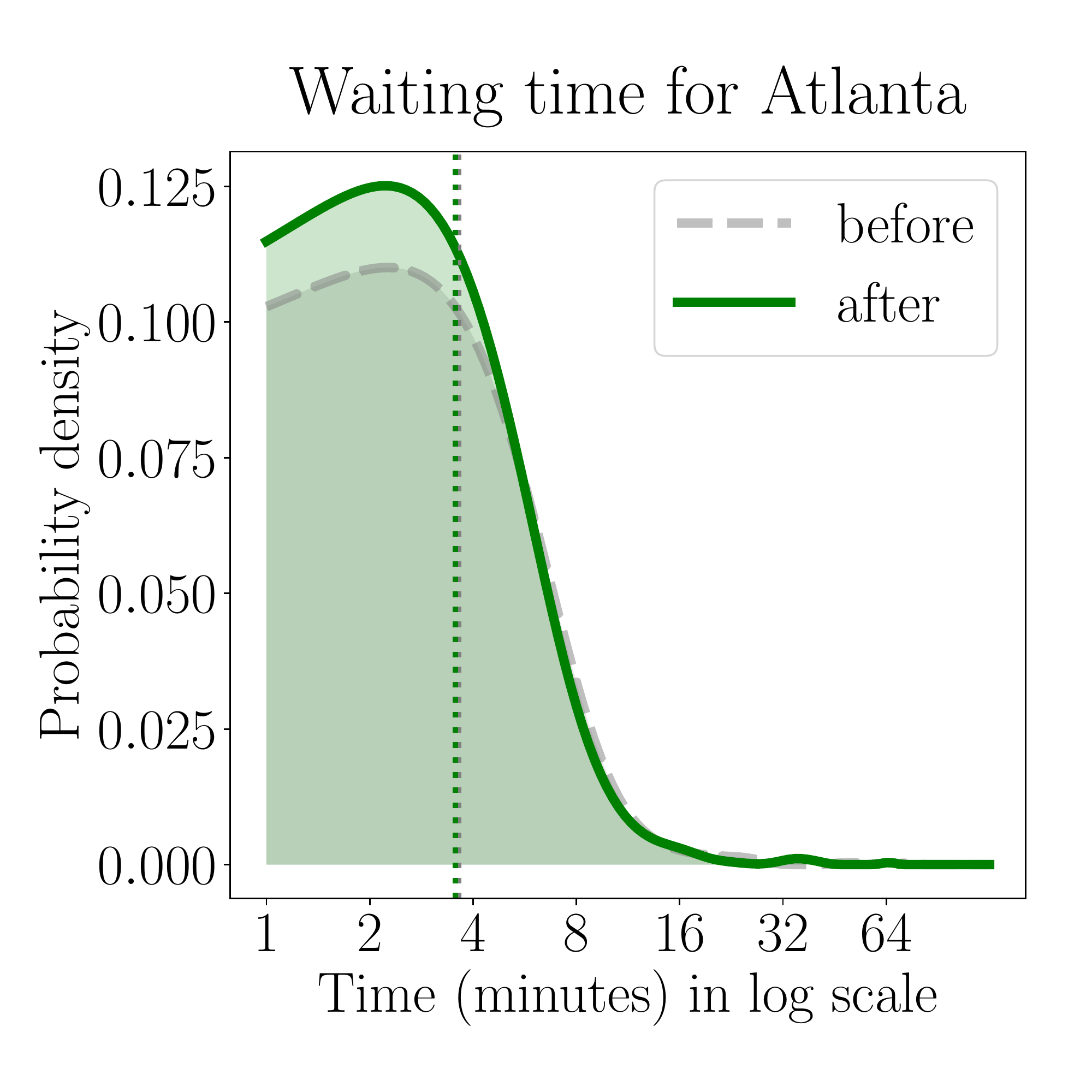}
  \end{subfigure}
  \begin{subfigure}[b]{0.25\textwidth}   
      \centering 
      \includegraphics[width=\textwidth]{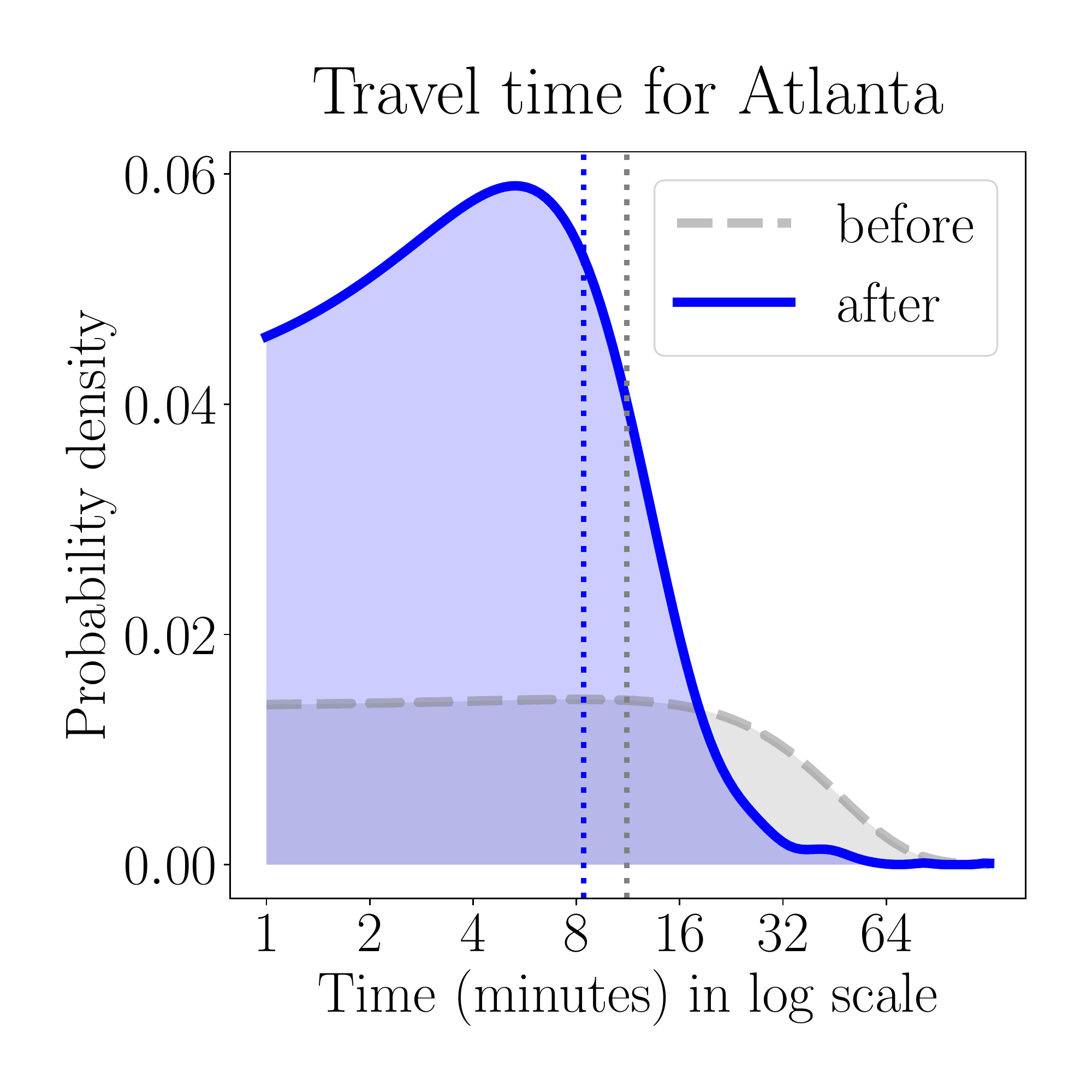}
  \end{subfigure}
  \caption{
Distribution of response, waiting, and travel time per call for high-priority calls before and after the redesign, calculated by kernel density estimation. Vertical dotted lines indicate the mean of each distribution. 
  } 
  \label{fig:time-dist-kde}
\end{figure}

We also report the average response time, waiting time, and travel time for all 911 calls and high-priority calls (e.g., serious crimes)
in Tables~\ref{tab:time-comparison} and \ref{tab:time-comparison-category86}. 
All three time metrics have worsened from March 2018 to March 2019, the year before the zone redesign, as a result of economic development and population growth.
After the new zone design was implemented,
all three metrics have seen decreased growth rates. For high-priority calls, in particular, these time metrics have decreased significantly.
In Figure~\ref{fig:time-dist-kde}, we also show the distributions of the three metrics for high-priority calls; the redesign has resulted in the probability density shifting left toward zero.

\begin{table}[!h]
\centering
\caption{Average response, waiting, and travel time per call for high-priority calls in six zones.}
\label{tab:time-comparison-category86}
\resizebox{\textwidth}{!}{%
\begin{tabular}{ccccccccccc}
  \toprule[1pt]\midrule[0.3pt]
\multirow{2}{*}{Category} & \multirow{2}{*}{Time period} & \multicolumn{7}{c}{Time per call (minutes)} & \multirow{2}{*}{\begin{tabular}[c]{@{}c@{}}Time increase\\year-on-year ratio (\%)\end{tabular}} \\ 
& & Zone 1 & Zone 2 & Zone 3 & Zone 4 & Zone 5 & Zone 6 & Citywide &  \\ 
\hline
\multirow{2}{*}{\begin{tabular}[c]{@{}c@{}}Response\end{tabular}} 
& Mar 2018 $\sim$ Mar 2019 (before) & 11.94 & 11.61 & 10.91 & 12.99 & 13.57 & 9.81 & 11.90 & $+10.90$ \\
& Mar 2019 $\sim$ Mar 2020 (after) & 13.92 & 10.59 & 12.06 & 9.36 & 10.34 & 9.80 & 11.21 & $-5.80$\\ 
\hline
\multirow{2}{*}{\begin{tabular}[c]{@{}c@{}}Waiting\end{tabular}}
& Mar 2018 $\sim$ Mar 2019 (before) & 3.49 & 3.16 & 4.69 & 3.22 & 3.39 & 2.46 & 3.55 & $+36.02$ \\
& Mar 2019 $\sim$ Mar 2020 (after) & 4.91 & 2.38 & 3.41 & 3.31 & 2.95 & 3.29 & 3.47 & $-2.25$\\ 
\hline
\multirow{2}{*}{\begin{tabular}[c]{@{}c@{}}Travel\end{tabular}} 
& Mar 2018 $\sim$ Mar 2019 (before) & 8.48 & 9.55 & 17.73 & 9.45 & 10.05 & 7.08 & 10.72 & $-0.74$ \\
& Mar 2019 $\sim$ Mar 2020 (after) & 9.59 & 8.64 & 8.28 & 7.73 & 6.59 & 7.36 & 8.07 & $-24.72$\\ 
\midrule[0.3pt]\bottomrule[1pt]
\end{tabular}%
}
\end{table}

\begin{table}[!h]
\centering
\caption{Average response, waiting, and travel time per call for all cases in six zones.}
\label{tab:time-comparison}
\resizebox{\textwidth}{!}{%
\begin{tabular}{ccccccccccc}
  \toprule[1pt]\midrule[0.3pt]
\multirow{2}{*}{Category} & \multirow{2}{*}{Time period} & \multicolumn{7}{c}{Time per call (minutes)} & \multirow{2}{*}{\begin{tabular}[c]{@{}c@{}}Time increase\\year-on-year ratio (\%)\end{tabular}} \\ 
& & Zone 1 & Zone 2 & Zone 3 & Zone 4 & Zone 5 & Zone 6 & Citywide &  \\ 
\hline
\multirow{2}{*}{\begin{tabular}[c]{@{}c@{}}Response\end{tabular}} 
& Mar 2018 $\sim$ Mar 2019 (before) & 30.35 & 32.21 & 27.81 & 30.70 & 30.42 & 26.07 & 29.69 & $+11.70$ \\
& Mar 2019 $\sim$ Mar 2020 (after) &28.72 & 24.77 & 24.70 & 28.59 & 36.20 & 33.25 & 31.26 & $+5.29$\\ 
\hline
\multirow{2}{*}{\begin{tabular}[c]{@{}c@{}}Waiting\end{tabular}}
& Mar 2018 $\sim$ Mar 2019 (before) & 18.41 & 17.75& 16.61 & 17.87 & 17.97 & 14.21 & 17.24 & $+12.83$ \\
& Mar 2019 $\sim$ Mar 2020 (after) & 18.00 & 13.96 & 14.63 & 16.56 & 21.74 & 19.63 & 19.08 & $+10.67$\\ 
\hline
\multirow{2}{*}{\begin{tabular}[c]{@{}c@{}}Travel\end{tabular}} 
& Mar 2018 $\sim$ Mar 2019 (before) & 8.20 & 7.88 & 7.49 & 8.60 & 6.80 & 6.27 & 7.50 & $+28.21$ \\
& Mar 2019 $\sim$ Mar 2020 (after) & 7.56 & 6.02 & 7.28 & 7.51 & 5.98 & 6.27 & 7.43 & $-0.93$\\ 
\midrule[0.3pt]\bottomrule[1pt]
\end{tabular}%
}
\end{table}

\begin{figure}[!h]
  \centering
  \begin{subfigure}[b]{0.25\textwidth}
      \centering
      \includegraphics[width=\textwidth]{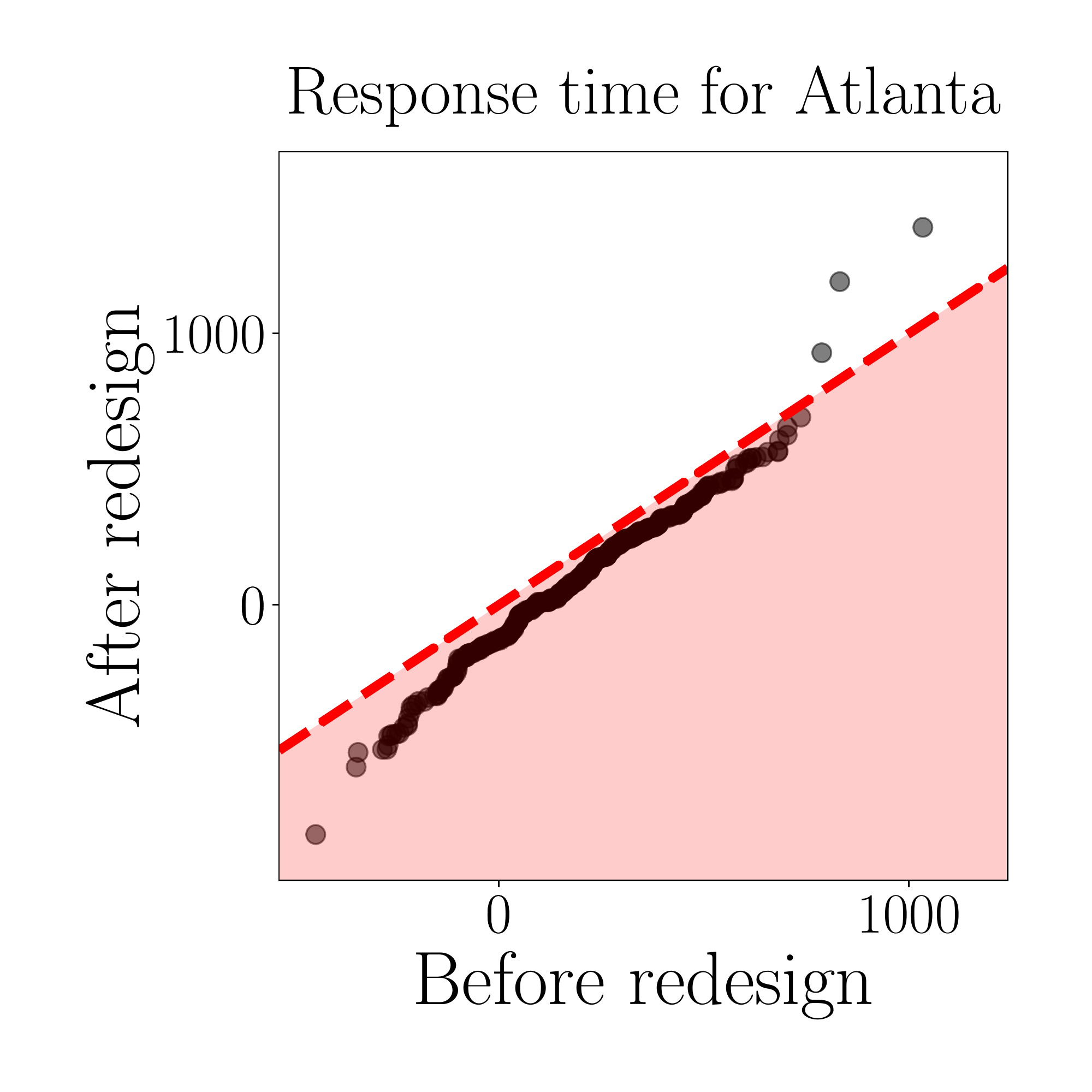}
  \end{subfigure}
  \begin{subfigure}[b]{0.25\textwidth}  
      \centering 
      \includegraphics[width=\textwidth]{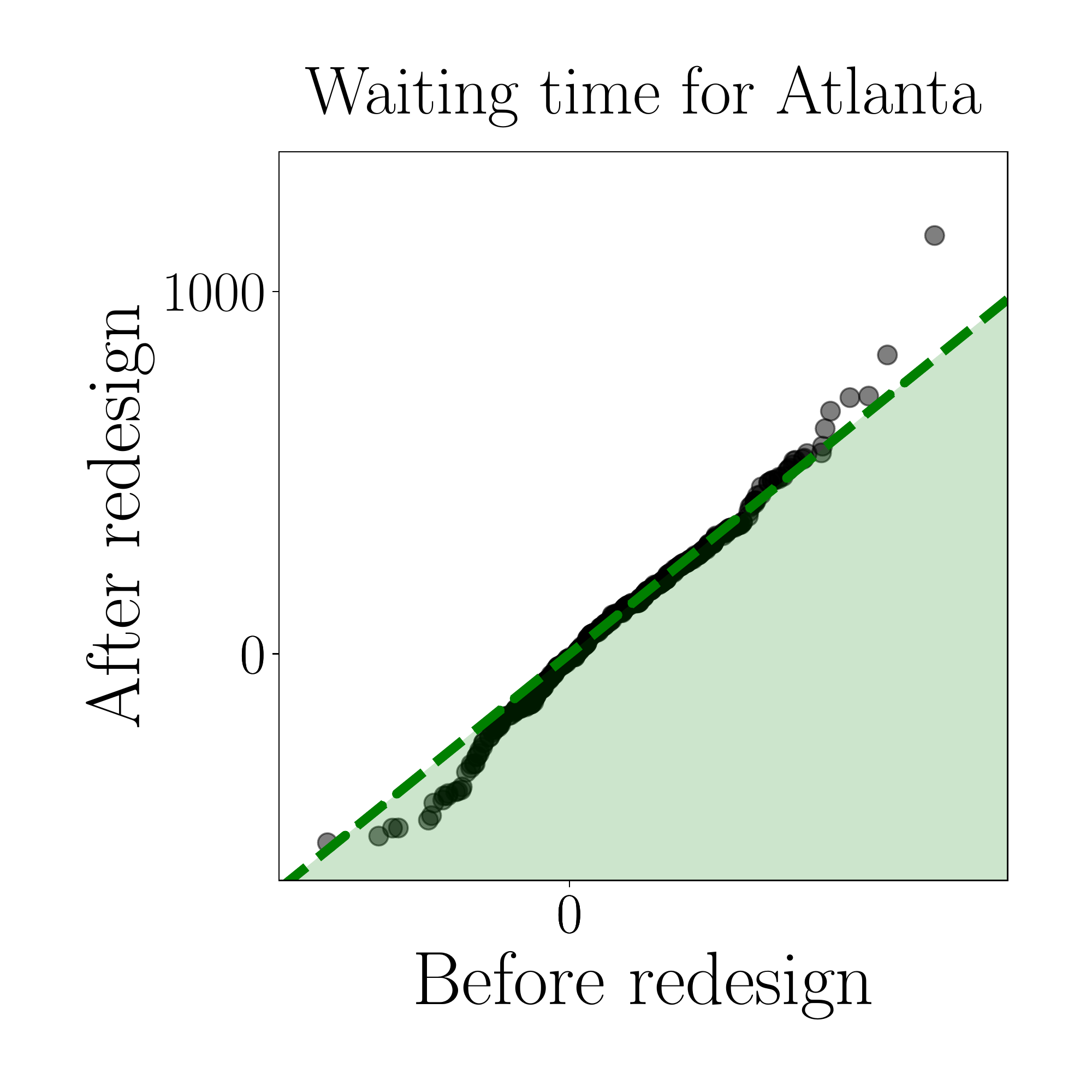}
  \end{subfigure}
  \begin{subfigure}[b]{0.25\textwidth}   
      \centering 
      \includegraphics[width=\textwidth]{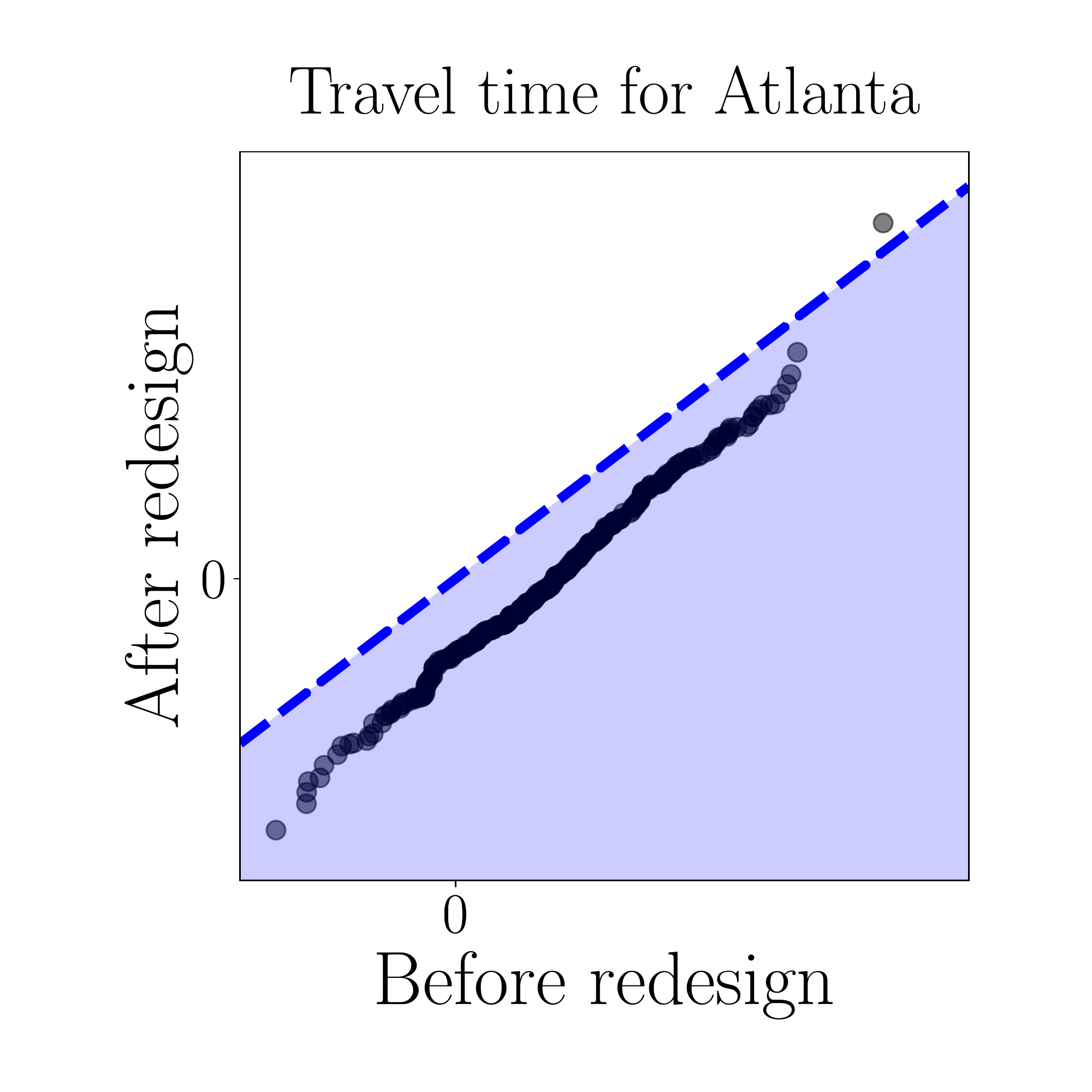}
  \end{subfigure}
  \caption{
  The difference in difference (DID) analysis of variance for response time (red), waiting time (green), and travel time (blue). Each point represents one day of the year. The x-axis corresponds to the variance difference for the two years prior to the redesign, and the y-axis corresponds to the variance difference for the two years before and after the redesign. The 45-degree dashed line passes through the origin. 
  Points below the dashed line indicate in a reduction in variance increase.
  } 
  \label{fig:time-diff-summary}
\end{figure}

\begin{figure}[!h]
  \centering
  \begin{subfigure}[b]{0.15\textwidth}
      \centering
      \includegraphics[width=\textwidth]{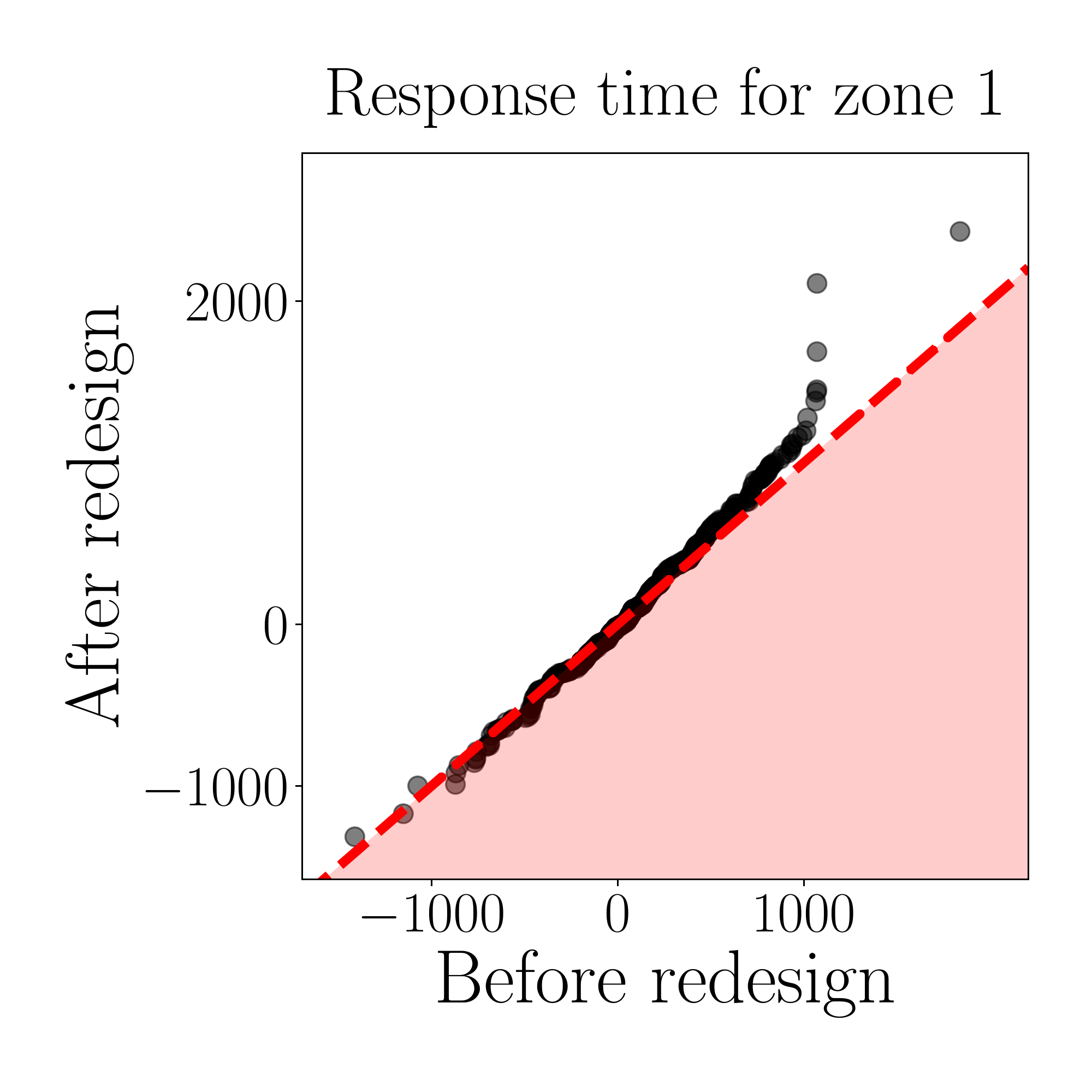}
  \end{subfigure}
  \begin{subfigure}[b]{0.15\textwidth}  
      \centering 
      \includegraphics[width=\textwidth]{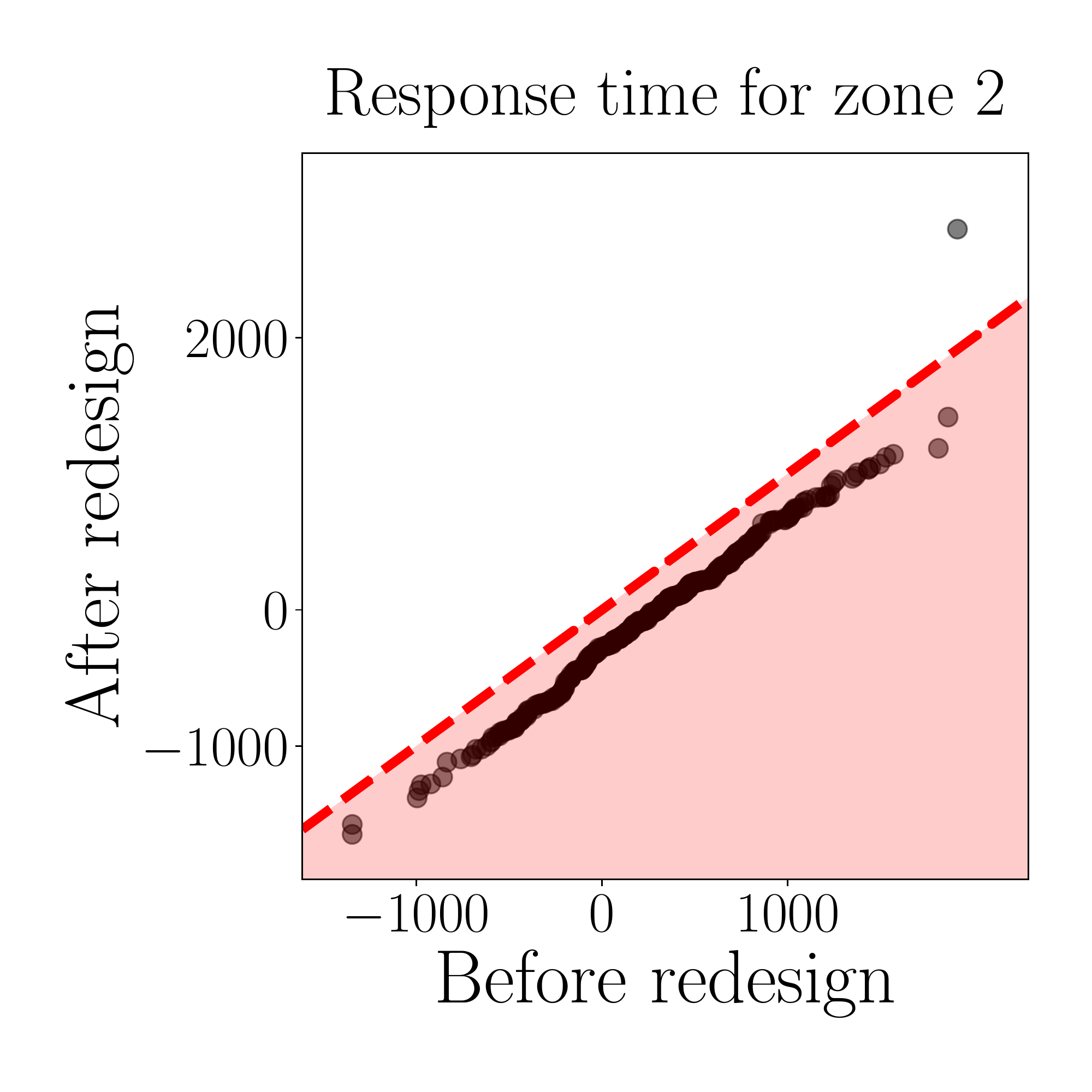}
  \end{subfigure}
  \begin{subfigure}[b]{0.15\textwidth}   
      \centering 
      \includegraphics[width=\textwidth]{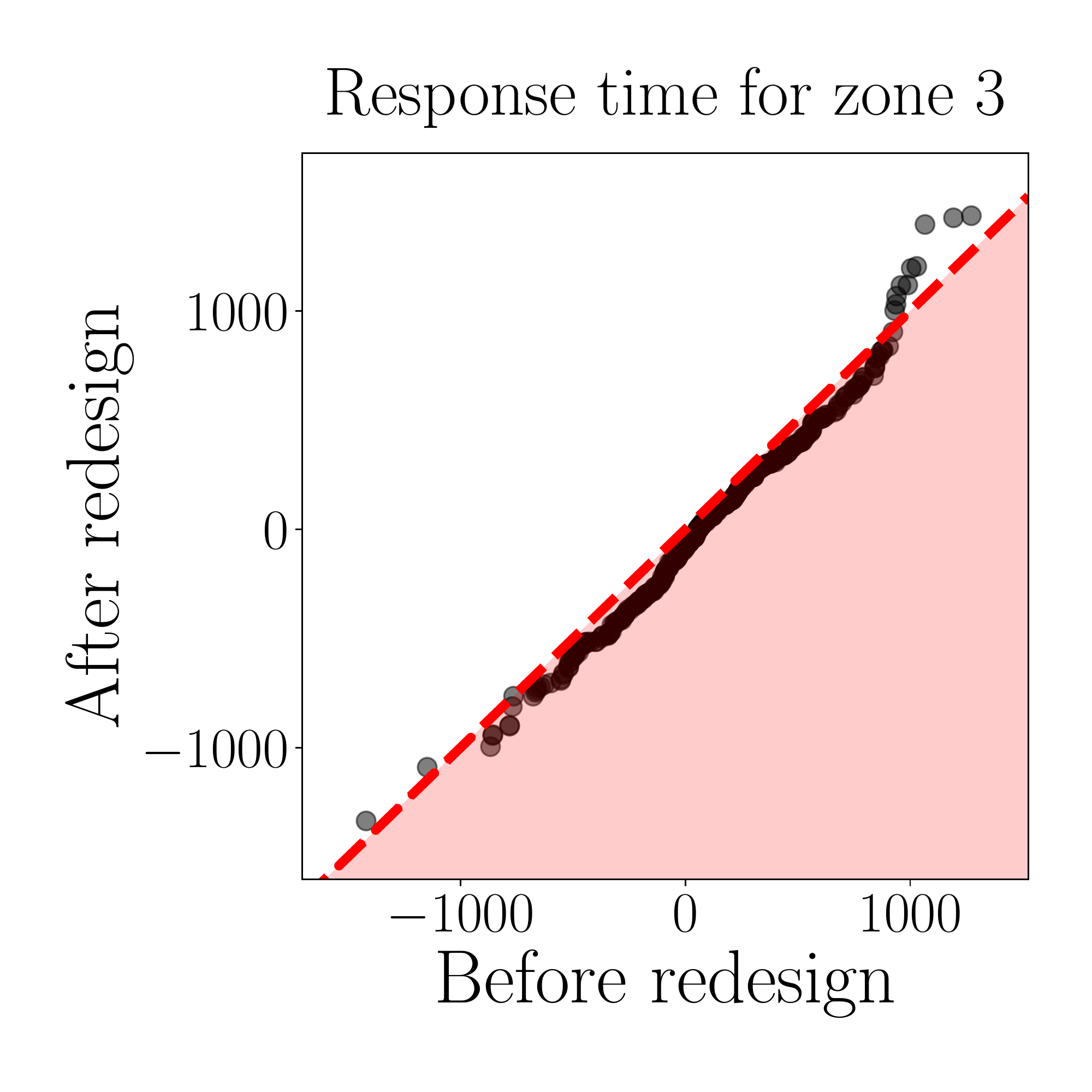}
  \end{subfigure}
  \begin{subfigure}[b]{0.15\textwidth}   
      \centering 
      \includegraphics[width=\textwidth]{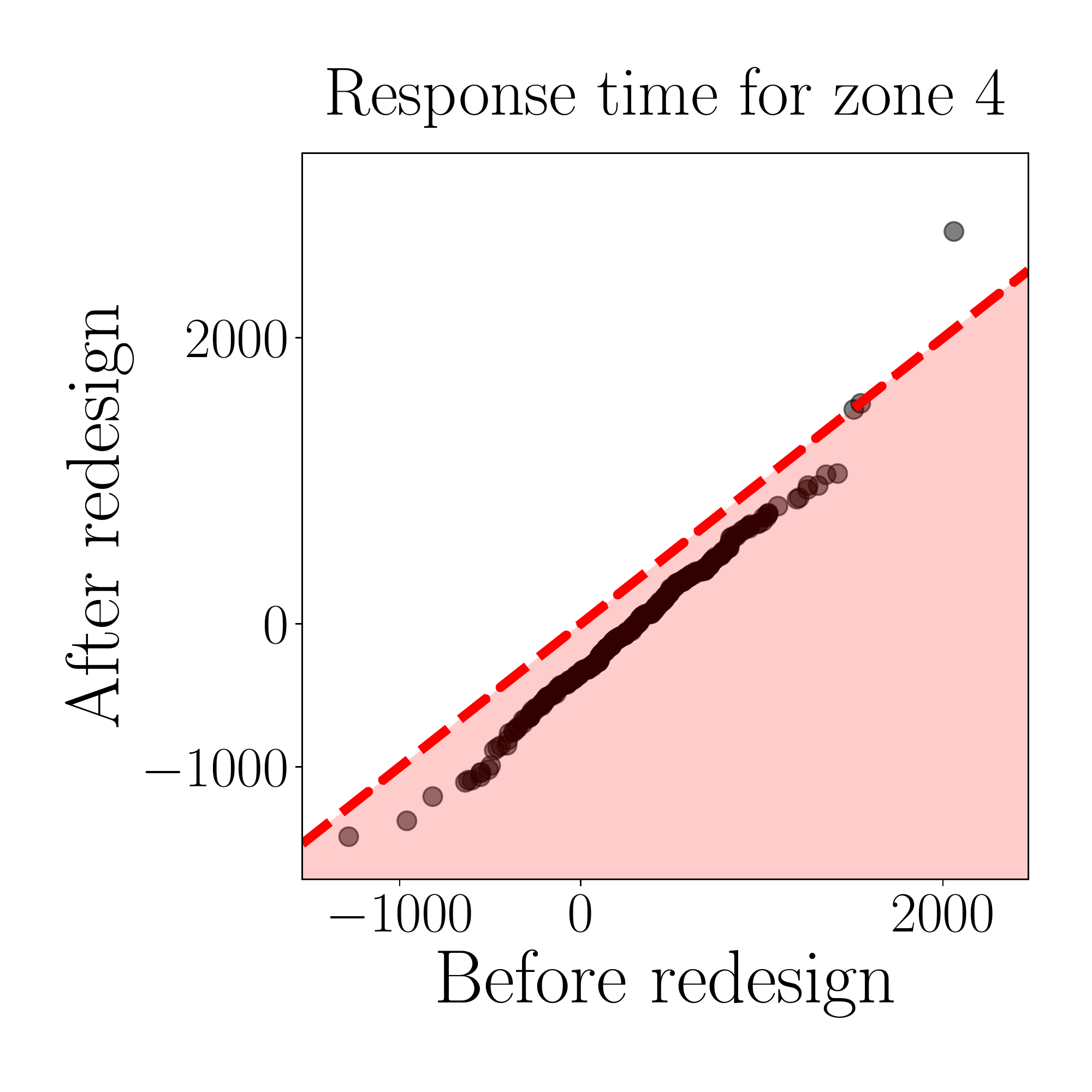}
  \end{subfigure}
  \begin{subfigure}[b]{0.15\textwidth}   
      \centering 
      \includegraphics[width=\textwidth]{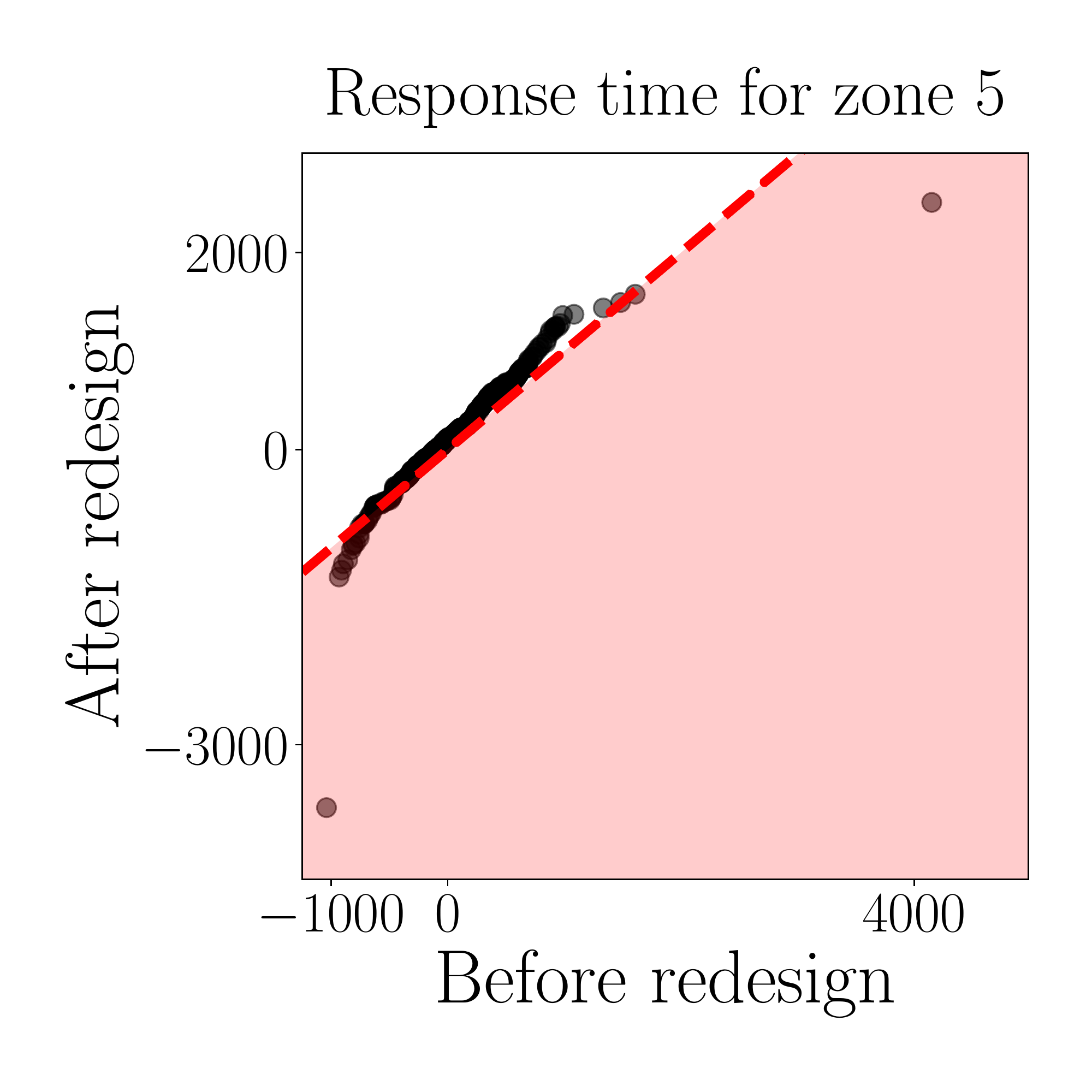}
  \end{subfigure}
  \begin{subfigure}[b]{0.15\textwidth}   
      \centering 
      \includegraphics[width=\textwidth]{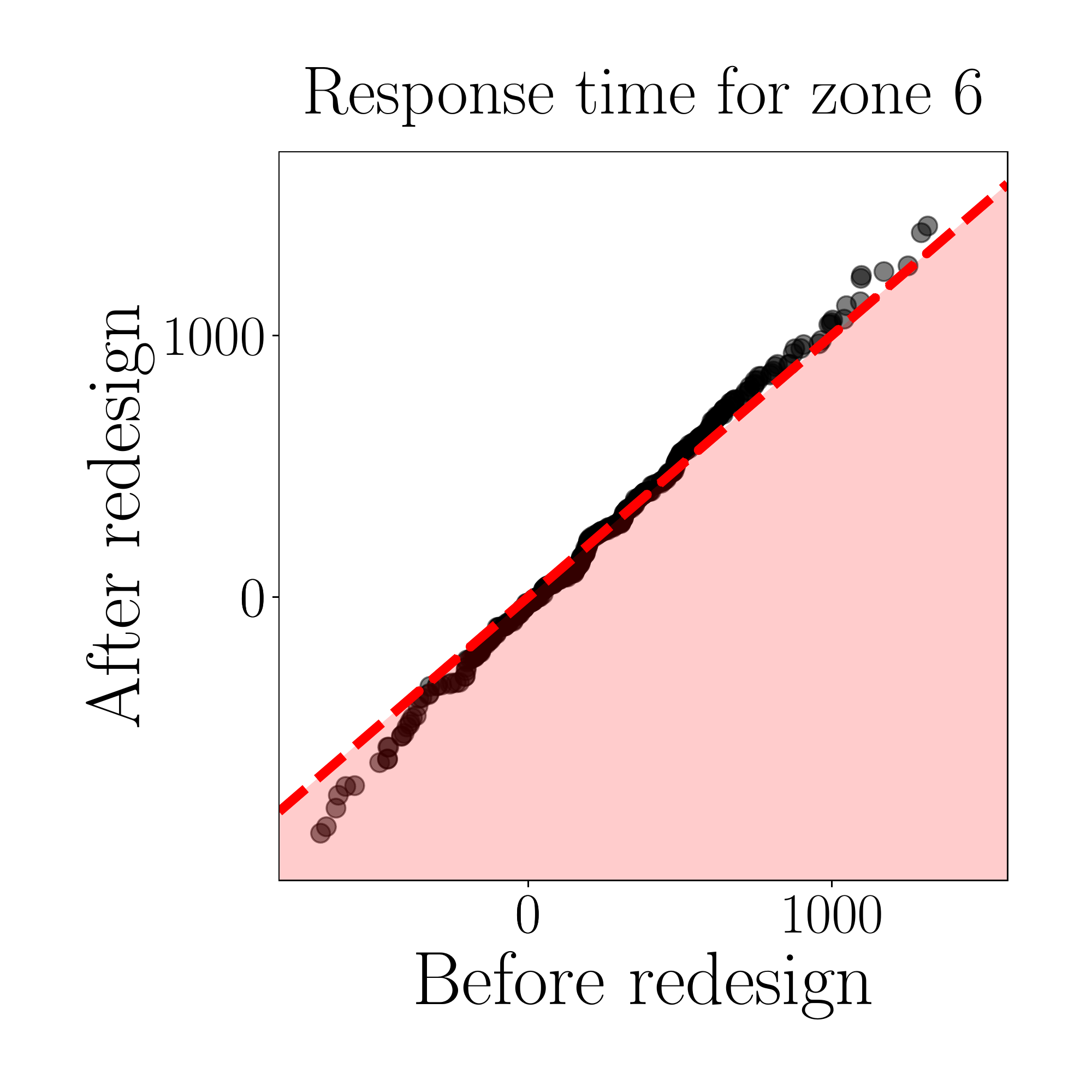}
  \end{subfigure}
  \vfill
  \begin{subfigure}[b]{0.15\textwidth}
      \centering
      \includegraphics[width=\textwidth]{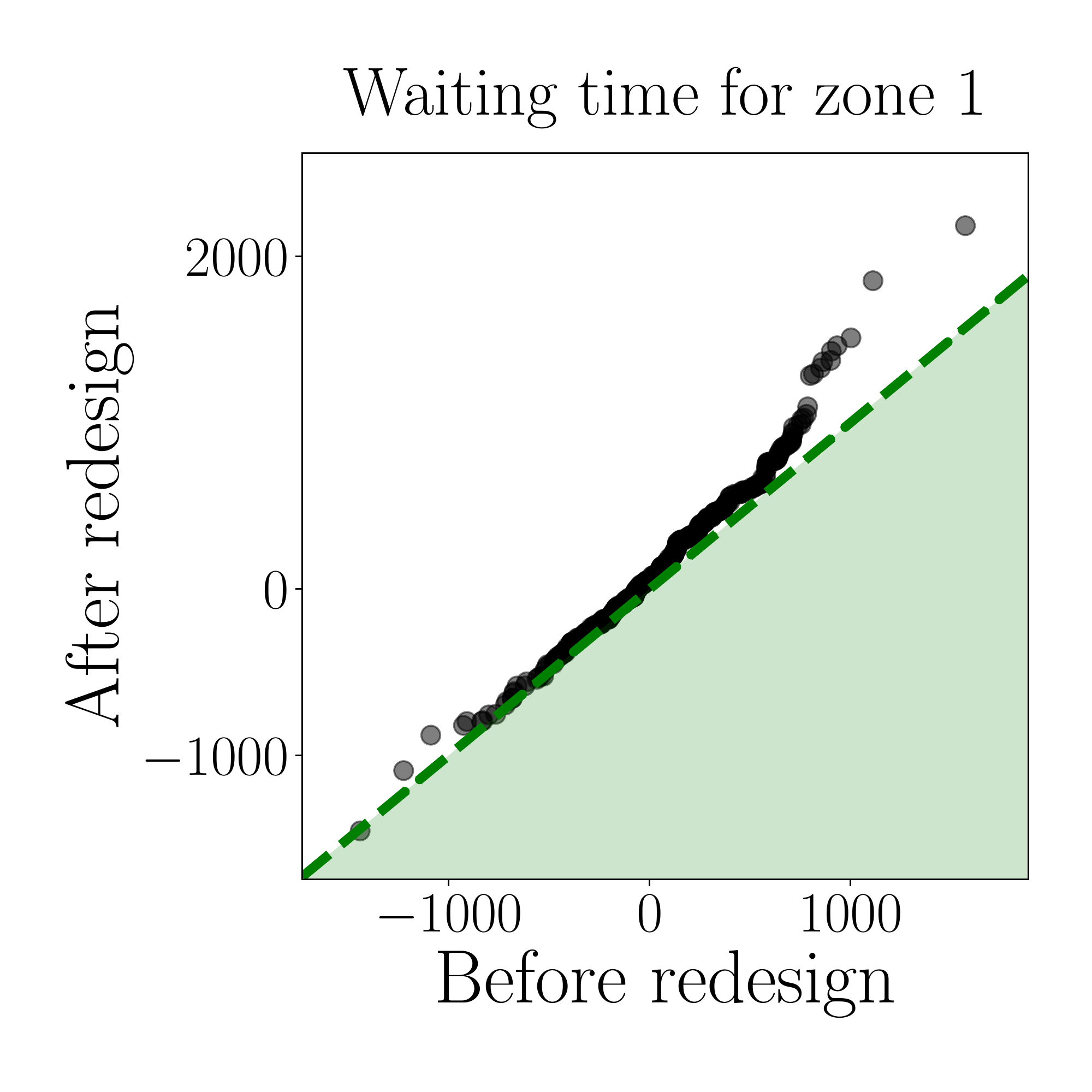}
  \end{subfigure}
  \begin{subfigure}[b]{0.15\textwidth}  
      \centering 
      \includegraphics[width=\textwidth]{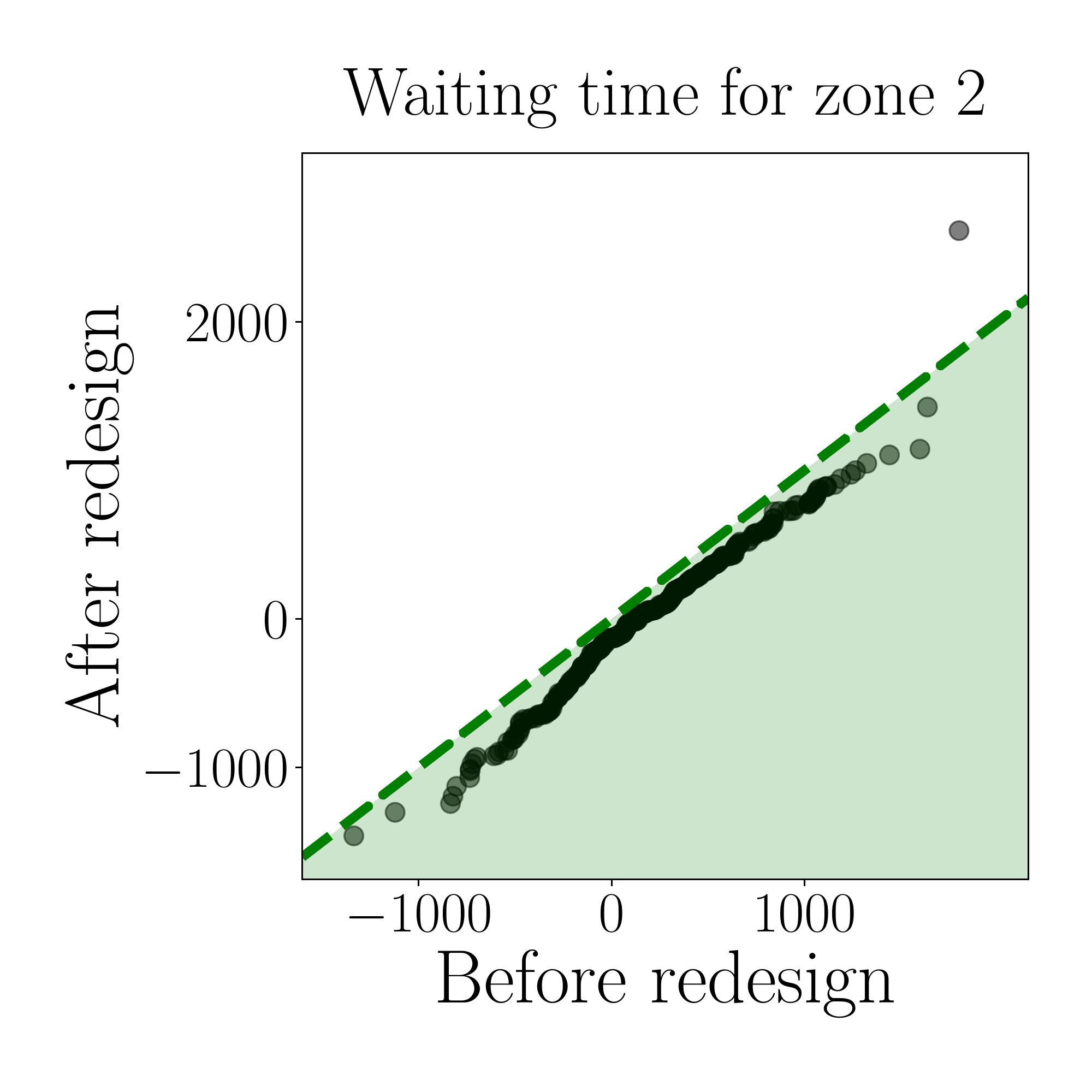}
  \end{subfigure}
  \begin{subfigure}[b]{0.15\textwidth}   
      \centering 
      \includegraphics[width=\textwidth]{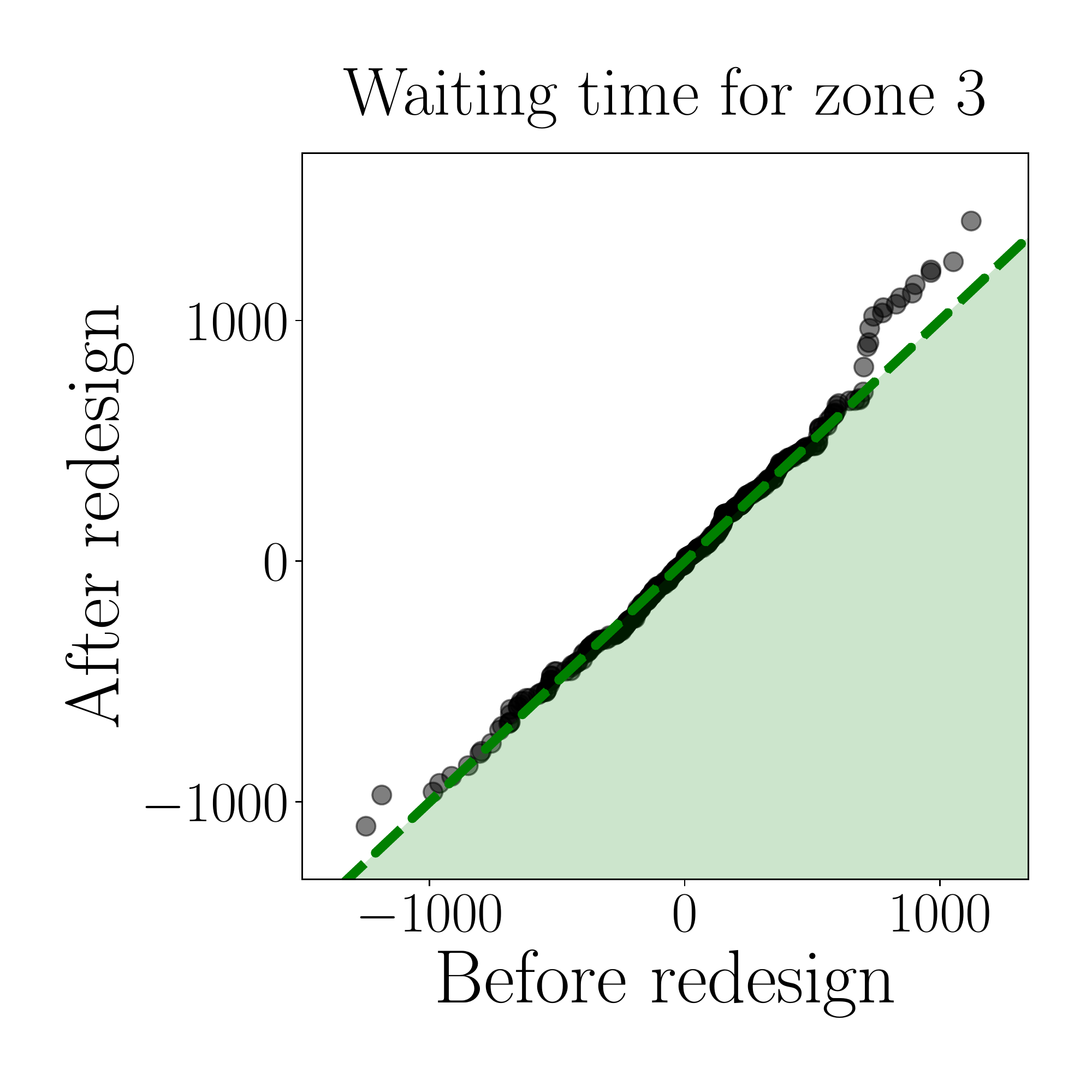}
  \end{subfigure}
  \begin{subfigure}[b]{0.15\textwidth}   
      \centering 
      \includegraphics[width=\textwidth]{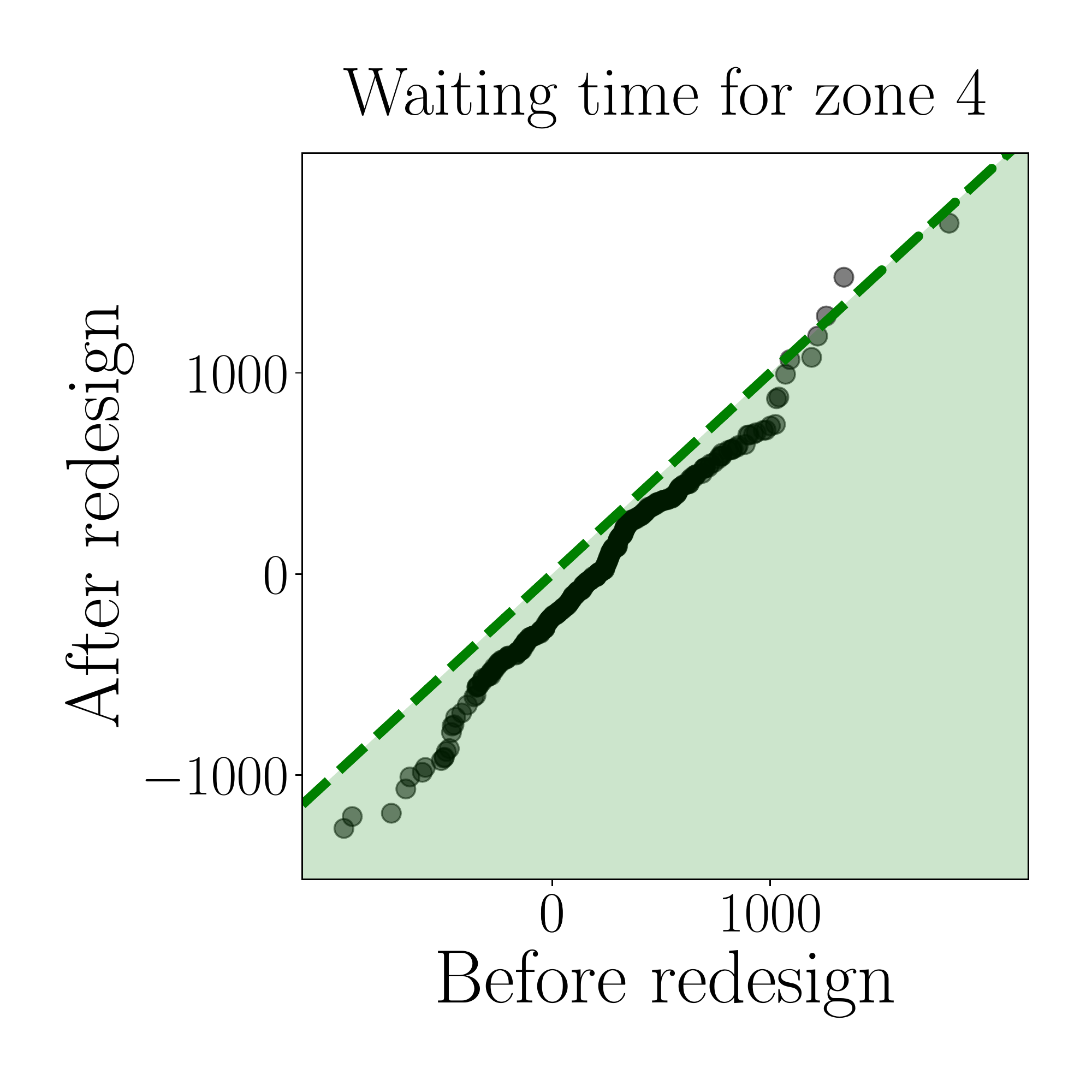}
  \end{subfigure}
  \begin{subfigure}[b]{0.15\textwidth}   
      \centering 
      \includegraphics[width=\textwidth]{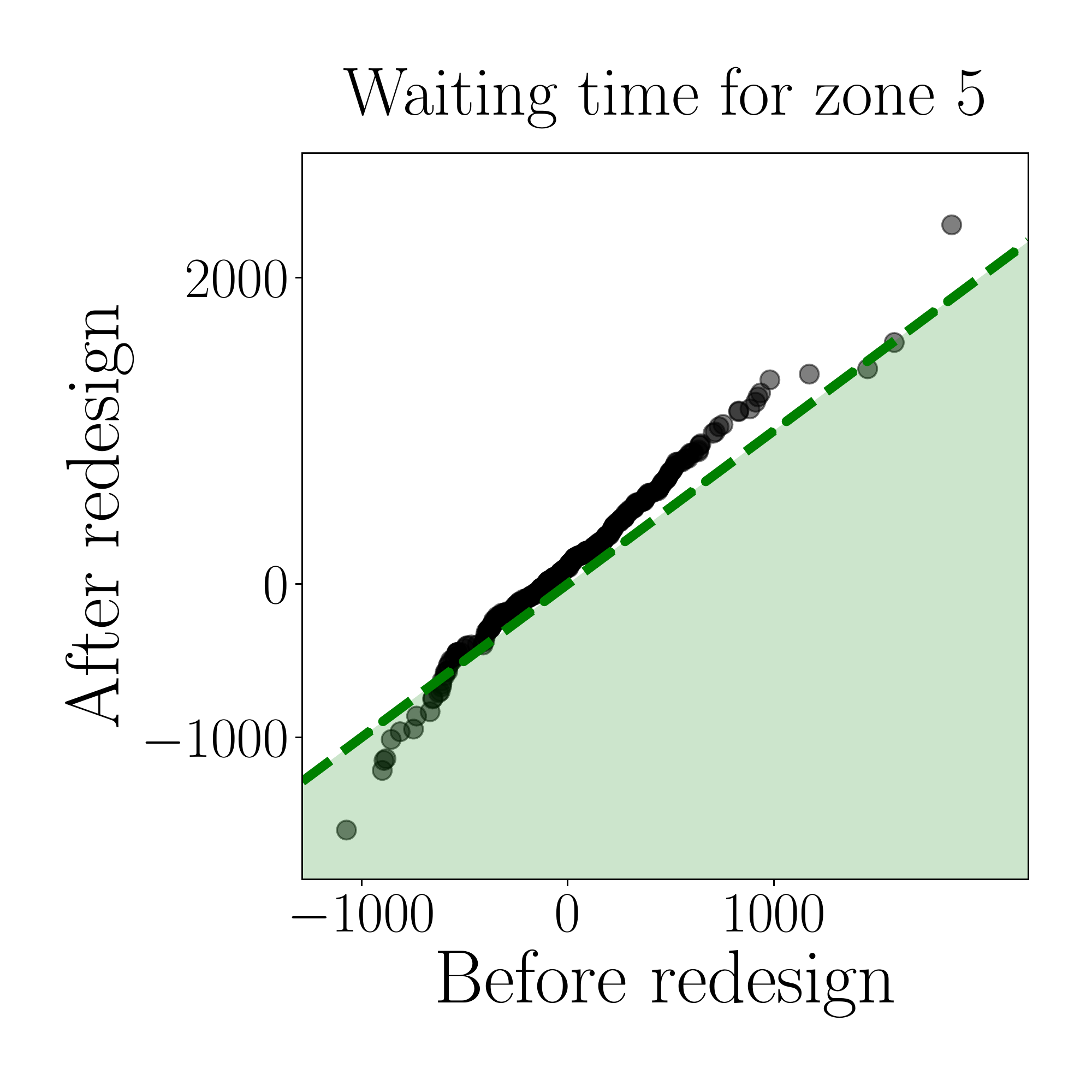}
  \end{subfigure}
  \begin{subfigure}[b]{0.15\textwidth}   
      \centering 
      \includegraphics[width=\textwidth]{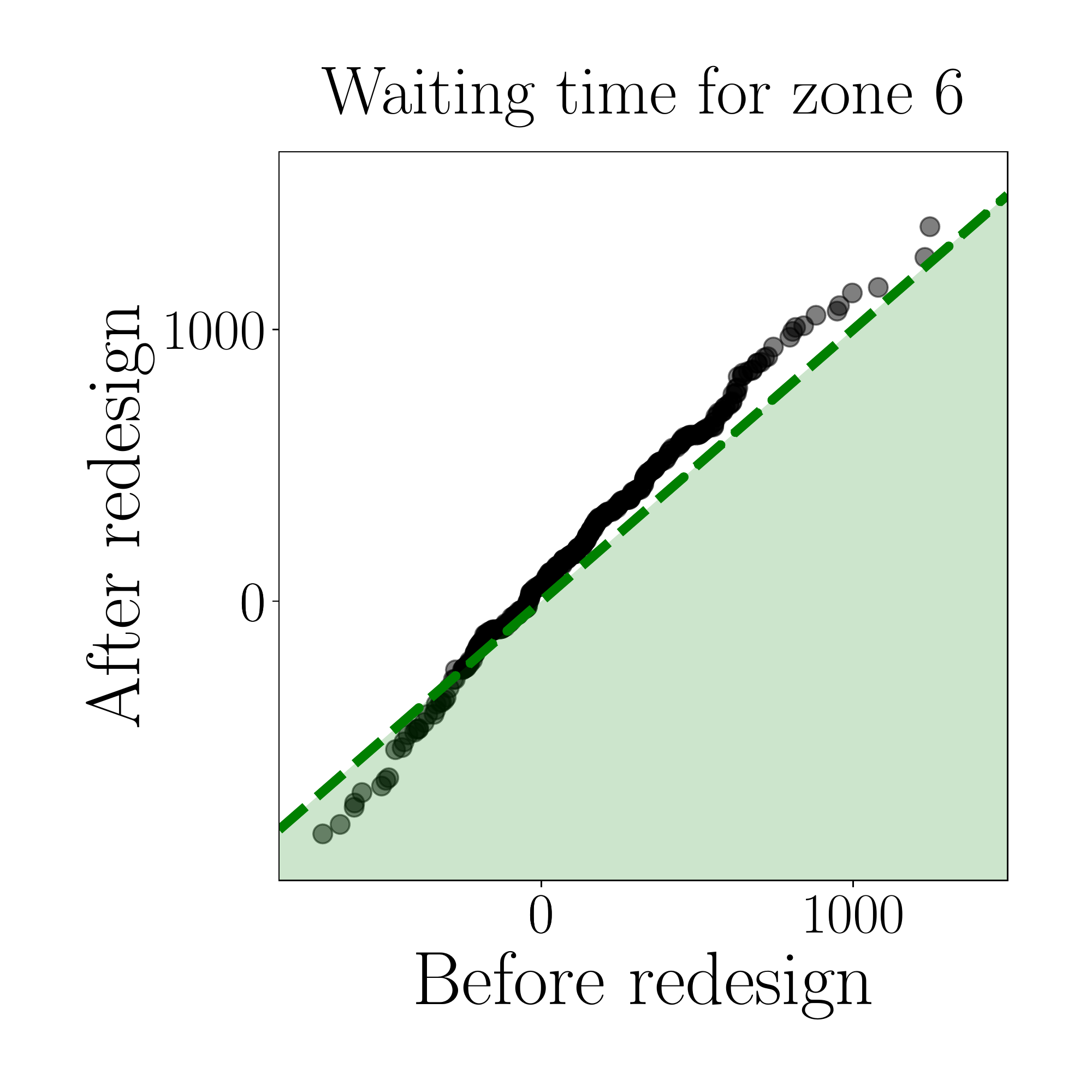}
  \end{subfigure}
  \vfill
  \begin{subfigure}[b]{0.15\textwidth}
      \centering
      \includegraphics[width=\textwidth]{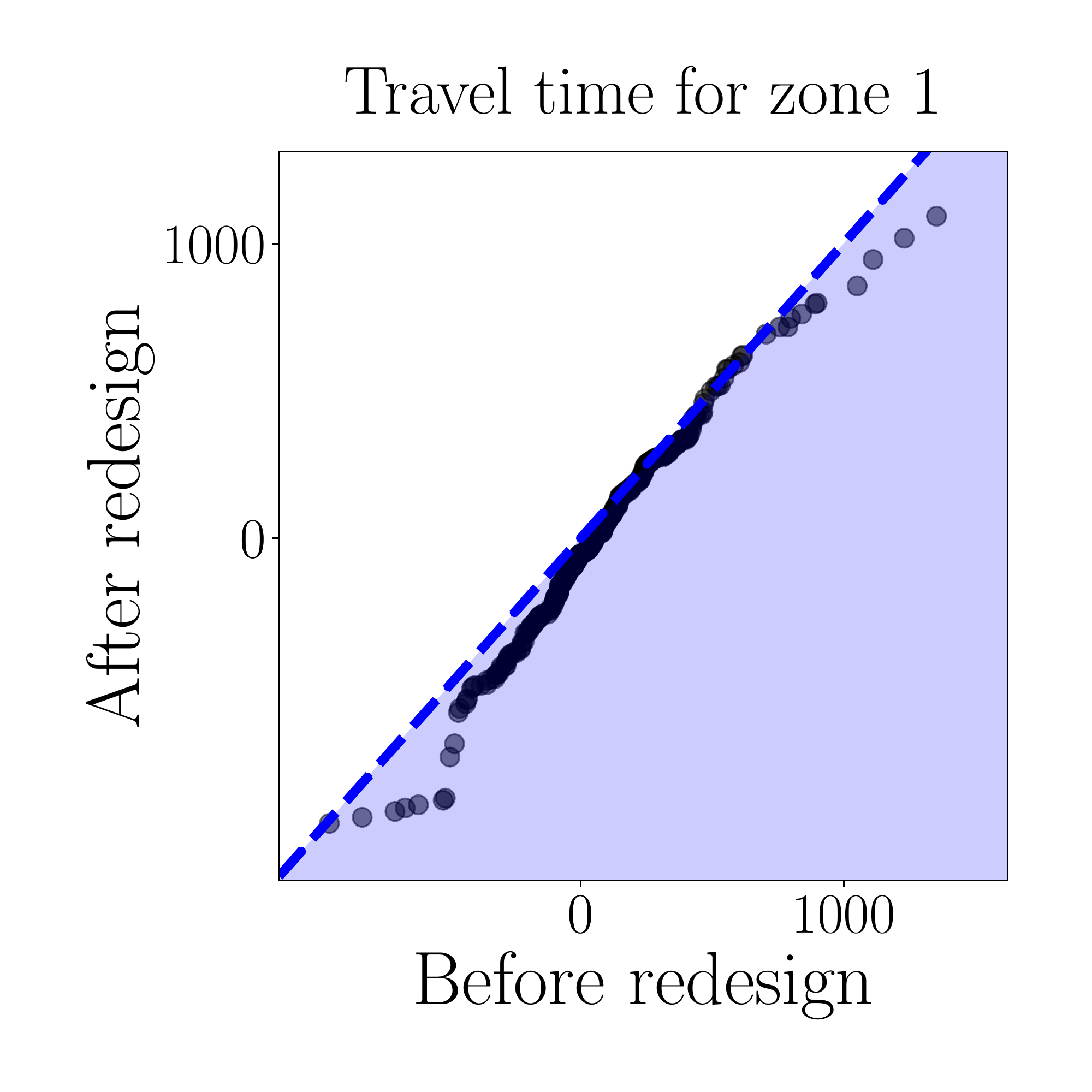}
  \end{subfigure}
  \begin{subfigure}[b]{0.15\textwidth}  
      \centering 
      \includegraphics[width=\textwidth]{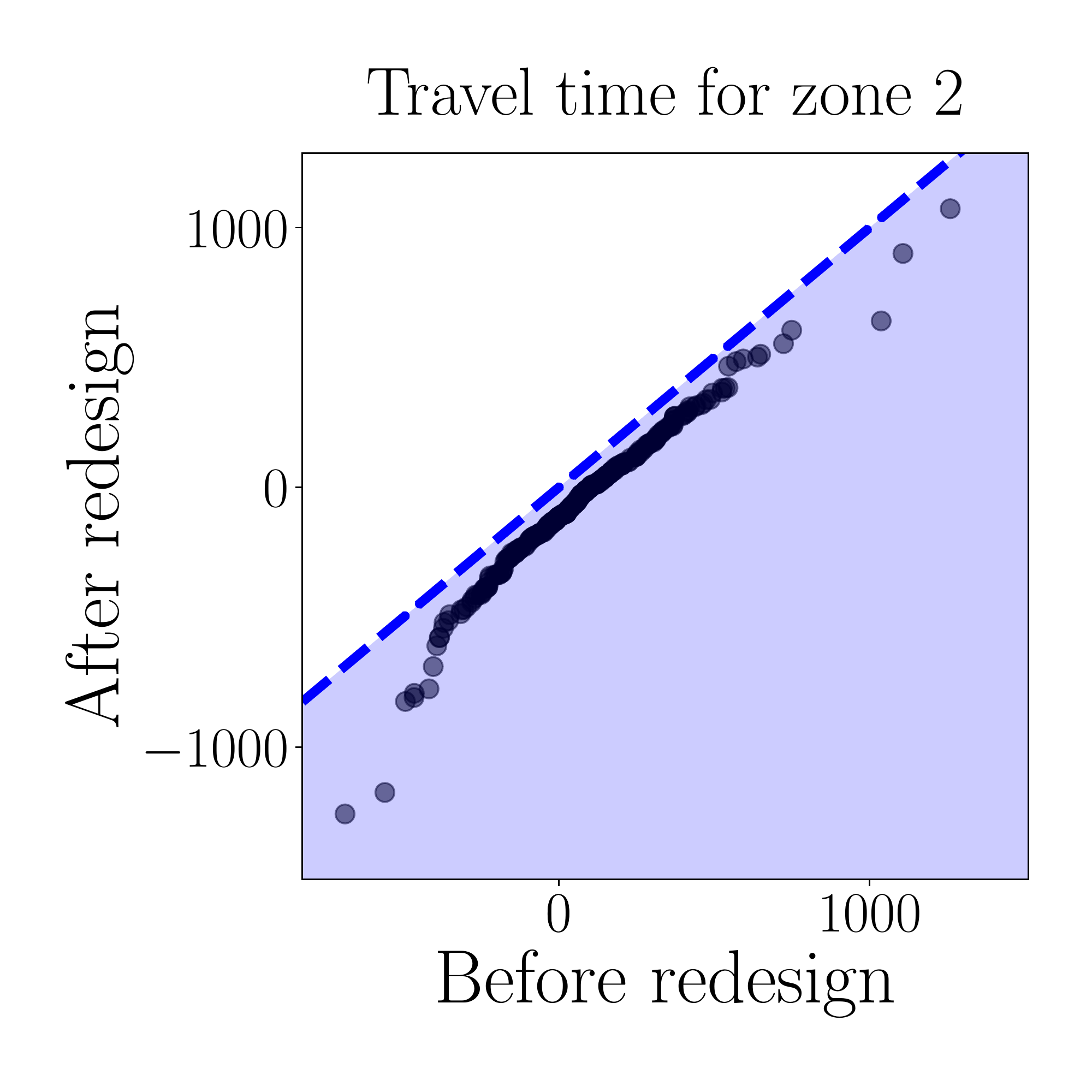}
  \end{subfigure}
  \begin{subfigure}[b]{0.15\textwidth}   
      \centering 
      \includegraphics[width=\textwidth]{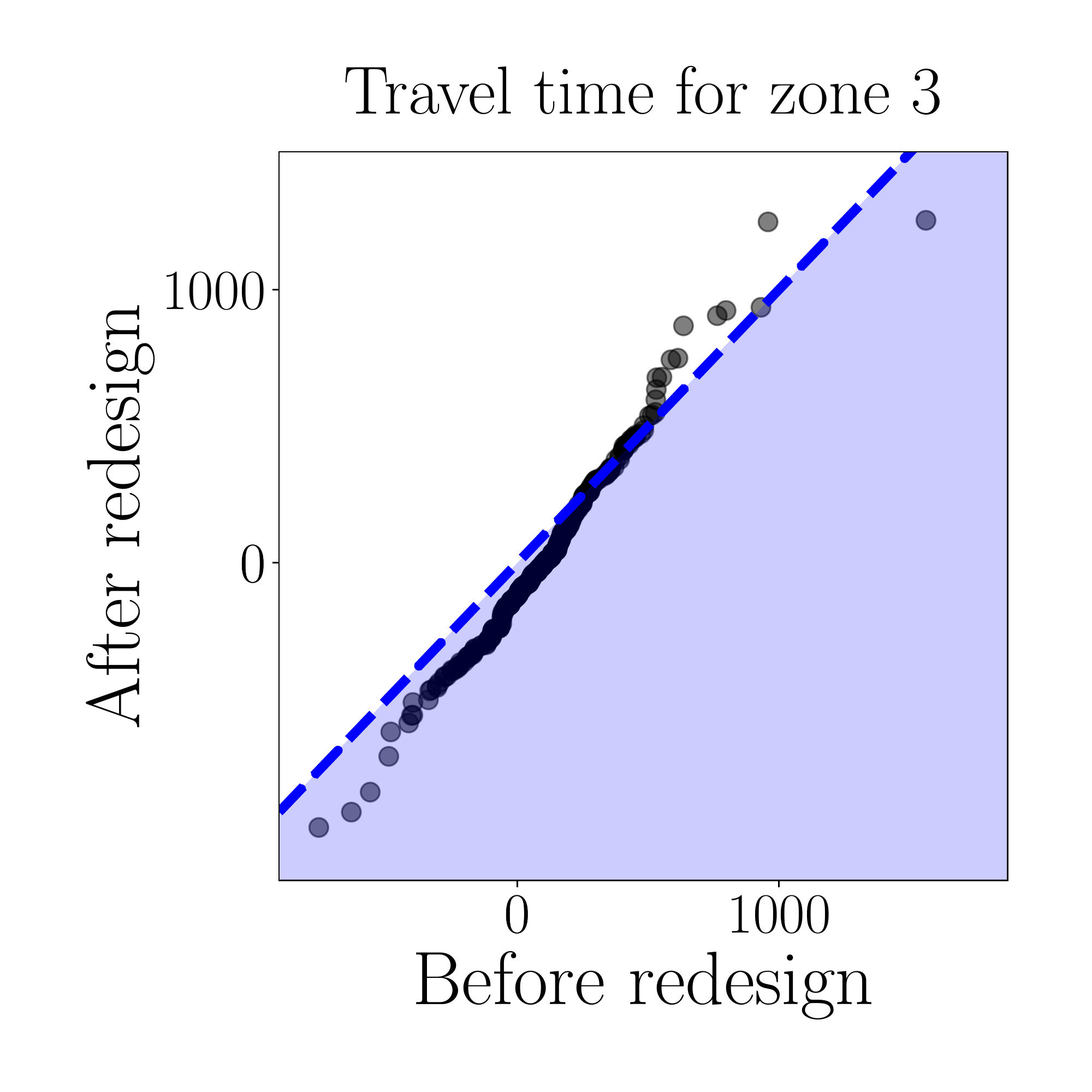}
  \end{subfigure}
  \begin{subfigure}[b]{0.15\textwidth}   
      \centering 
      \includegraphics[width=\textwidth]{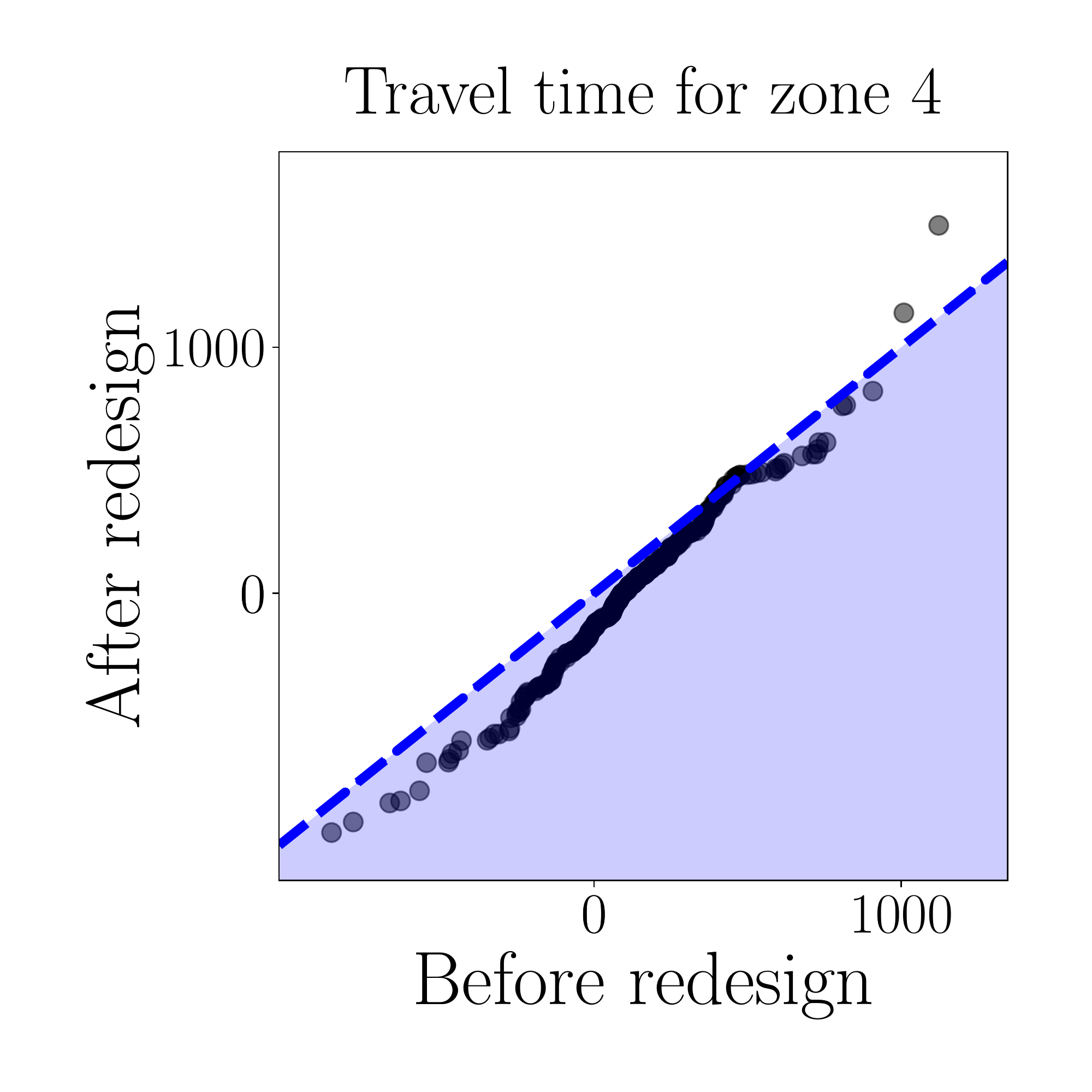}
  \end{subfigure}
  \begin{subfigure}[b]{0.15\textwidth}   
      \centering 
      \includegraphics[width=\textwidth]{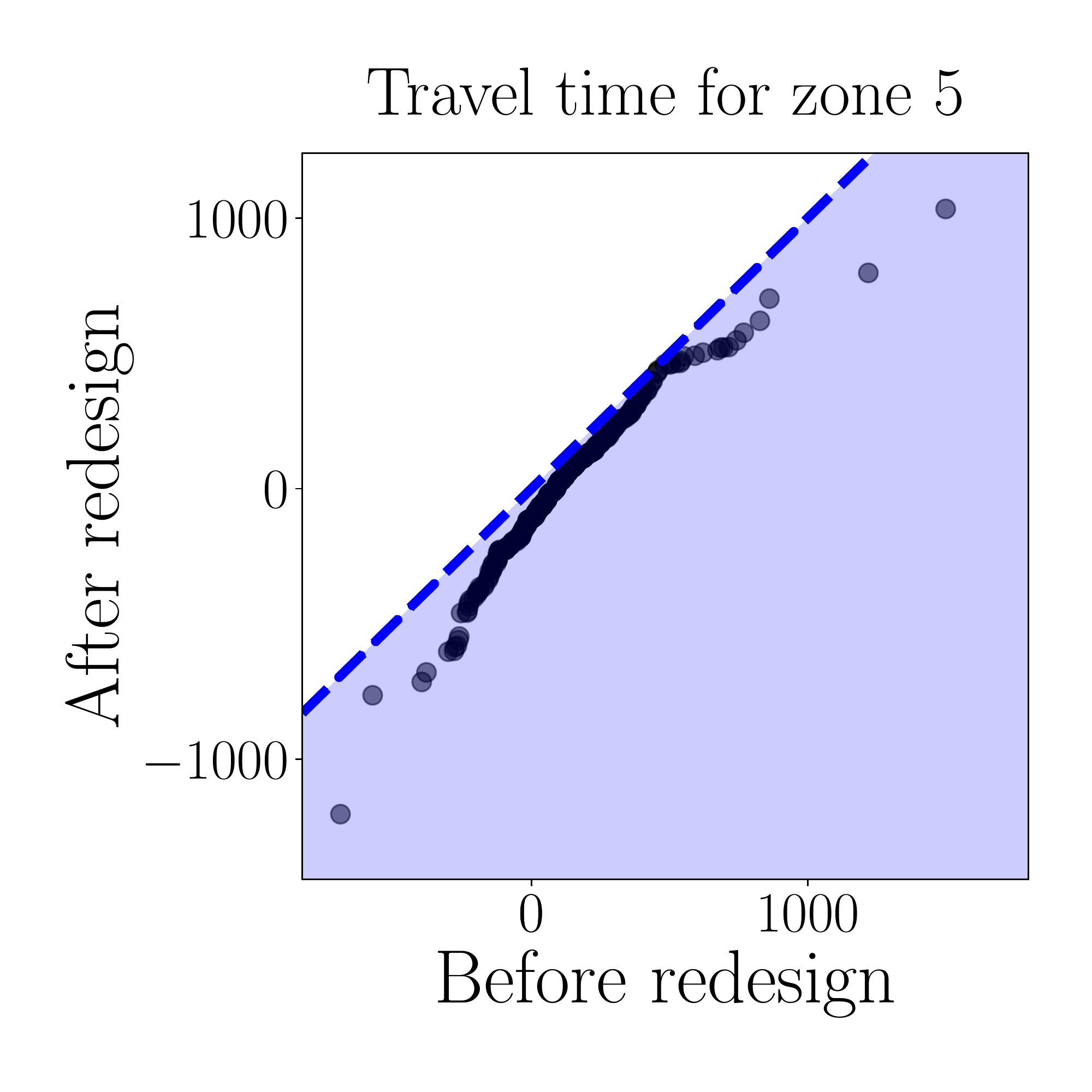}
  \end{subfigure}
  \begin{subfigure}[b]{0.15\textwidth}   
      \centering 
      \includegraphics[width=\textwidth]{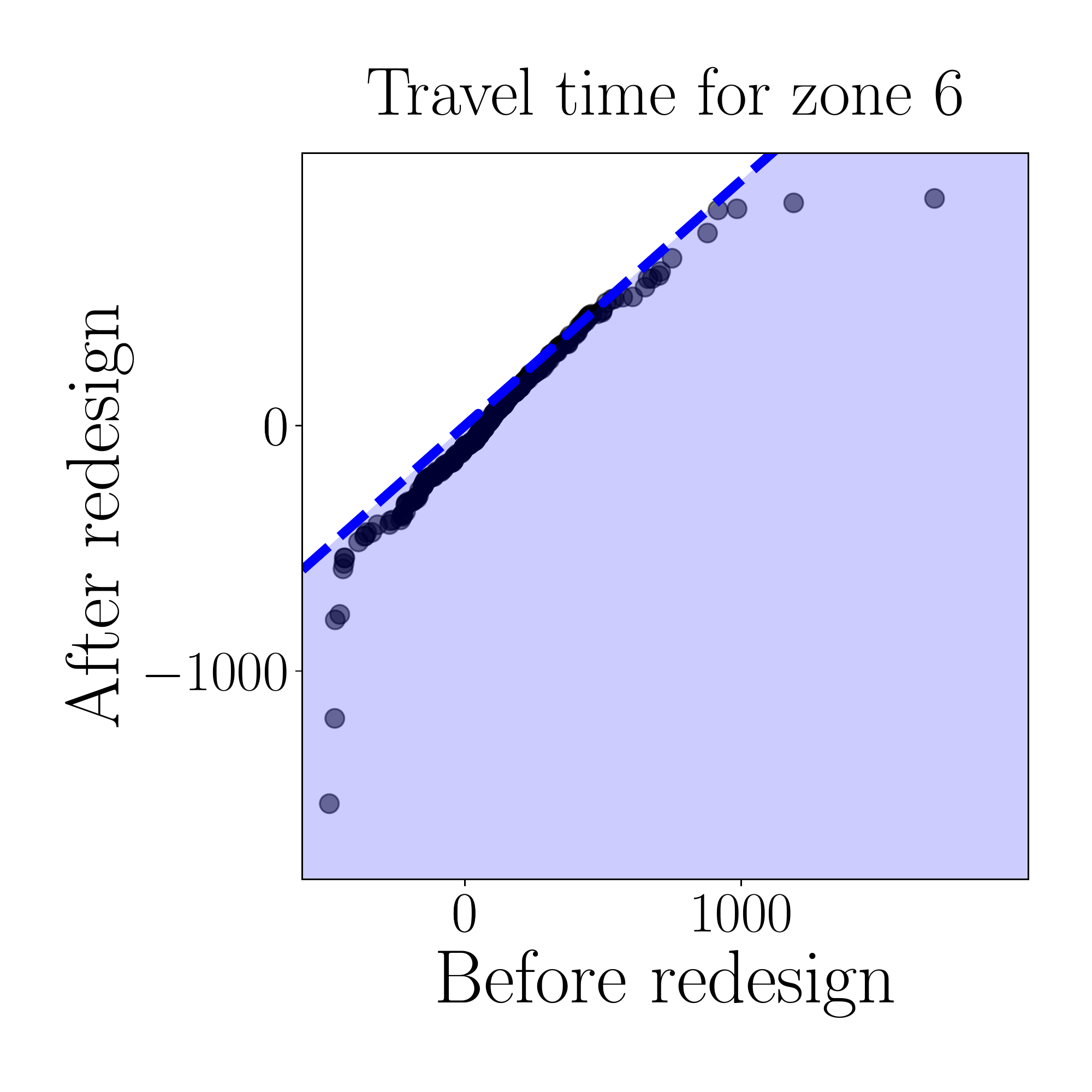}
  \end{subfigure}
  \caption{
  Difference in difference (DID) of variance for response time (1st row), waiting time (2nd row), and travel time (3rd row) in each zone.
  } 
  \label{fig:time-diff-comparison}
\end{figure}

To control for the time series effect, we also perform a \emph{difference in difference} (DID) analysis for the yearly variance change of the three time metrics.
We calculate the variance of these metrics on each day of the year and then compute the variance difference between two years. We then compare the variance difference before and after the redesign.
In Figure~\ref{fig:time-diff-summary}, we report this comparison for the entire Atlanta. For travel time, the result shows that the variance difference after the redesign is significantly lower than the difference before the redesign, as most points are below the 45-degree line. For response and waiting time, the figure also shows a slight decrease after the redesign.
In Figure~\ref{fig:time-diff-comparison}, we show the same comparison for each zone separately. For Zones 2, 3, and 4, which have experienced upticks in crime rates and police workloads (see Figure~\ref{fig:design-comparison-newdata}), the variance differences for all three metrics have decreased.



\section{Final Remarks}
\label{sec:conclusion}

In this paper, we presented a data-driven framework for police zone redesign.
The framework includes three main components:
a spatial queueing model for police patrol and emergency response operations, a statistical model to predict emergency call arrivals and police service rates, and a mixed-integer linear programming model for zone districting optimization.
We applied this framework to redesigning police zones in Atlanta. The Atlanta Police Department officially adopted our proposed redesign plan in March 2019. By analyzing data before and after the zone redesign, we show that the new design has reduced the response time to high-priority 911 calls by 5.8\% and the imbalance of police workload among different zones by 43\%.


We note several potential limitations in our method that can be improved in future research:
(a) our model assumes that the occurrence of crime events and 911 calls are exogenous to the police zone design; however, we observe in the post-implementation data that the call rates changed after the zone redesign, which imply that these factors may be endogenous and affected by the zone redesign;
(b) our method assumes that beat designs are fixed, and hence the zone design problem can be formulated as an assignment problem. For a complete redesign, beat designs and zone designs can be optimized jointly.

\section*{Acknowledgement}
We are grateful to Chief Erika Shields, Major John Quigley, Captain Jacquelyn Gwinn-Villaroel, and other police officers at the Atlanta Police Department for their valuable inputs and assists throughout this project. We would like to thank grateful to Mr. Song Wei for his help to prepare the data and the post-implementation results. This research is supported in part by the Atlanta Police Foundation and NSF (CMMI-2015787).

\bibliographystyle{informs2014}
\bibliography{refs}

\begin{thebibliography}{56}
\providecommand{\natexlab}[1]{#1}
\providecommand{\url}[1]{\texttt{#1}}
\providecommand{\urlprefix}{URL }

\bibitem[{{Atlanta Police Department}(2019)}]{APD2019}
{Atlanta Police Department} (2019) Atlanta police make changes to four zone
  boundaries as part of beat redesign.
  \url{http://www.atlantapd.org/Home/Components/News/News/190/} (accessed May
  1, 2020).

\bibitem[{{Atlanta Police Department}(2020)}]{APD2020}
{Atlanta Police Department} (2020) Atlanta police department zones.
  \url{https://www.atlantapd.org/community/apd-zones} (accessed May 1, 2020).

\bibitem[{Bammi(1975)}]{bammi1975allocation}
Bammi D (1975) Allocation of police beats to patrol units to minimize response
  time to calls for service. \emph{Computers \& Operations Research}
  2(1):1--12.

\bibitem[{Benveniste(1985)}]{benveniste1985solving}
Benveniste R (1985) Solving the combined zoning and location problem for
  several emergency units. \emph{Journal of the Operational Research Society}
  36(5):433--450.

\bibitem[{Bodily(1978)}]{bodily1978police}
Bodily SE (1978) Police sector design incorporating preferences of interest
  groups for equality and efficiency. \emph{Management Science}
  24(12):1301--1313.

\bibitem[{Bucarey et~al.(2015)Bucarey, Ord{\'o}{\~n}ez, \protect\BIBand{}
  Bassaletti}]{bucarey2015shape}
Bucarey V, Ord{\'o}{\~n}ez F, Bassaletti E (2015) Shape and balance in police
  districting. \emph{Applications of Location Analysis}, 329--347 (Springer).

\bibitem[{Camacho-Collados et~al.(2015)Camacho-Collados, Liberatore,
  \protect\BIBand{} Angulo}]{camacho2015multi}
Camacho-Collados M, Liberatore F, Angulo JM (2015) A multi-criteria police
  districting problem for the efficient and effective design of patrol sector.
  \emph{European journal of operational research} 246(2):674--684.

\bibitem[{Chaiken \protect\BIBand{} Larson(1972)}]{chaiken1972methods}
Chaiken JM, Larson RC (1972) Methods for allocating urban emergency units: a
  survey. \emph{Management Science} 19(4-part-2):110--130.

\bibitem[{Chelst \protect\BIBand{} Jarvis(1979)}]{chelst1979estimating}
Chelst K, Jarvis JP (1979) Estimating the probability distribution of travel
  times for urban emergency service systems. \emph{Operations Research}
  27(1):199--204.

\bibitem[{Chen et~al.(2019)Chen, Cheng, \protect\BIBand{}
  Ye}]{chen2019designing}
Chen H, Cheng T, Ye X (2019) Designing efficient and balanced police patrol
  districts on an urban street network. \emph{International Journal of
  Geographical Information Science} 33(2):269--290.

\bibitem[{Cheung et~al.(2015)Cheung, Yoon, \protect\BIBand{}
  Chow}]{cheung2015optimization}
Cheung CY, Yoon H, Chow AH (2015) Optimization of police facility deployment
  with a case study in greater london area. \emph{Journal of Facilities
  Management} .

\bibitem[{Chow et~al.(2015)Chow, Cheung, \protect\BIBand{}
  Yoon}]{chow2015optimization}
Chow AH, Cheung C, Yoon H (2015) Optimization of police facility locationing.
  \emph{Transportation research record} 2528(1):60--68.

\bibitem[{Curtin et~al.(2010)Curtin, Hayslett-McCall, \protect\BIBand{}
  Qiu}]{curtin2010determining}
Curtin KM, Hayslett-McCall K, Qiu F (2010) Determining optimal police patrol
  areas with maximal covering and backup covering location models.
  \emph{Networks and Spatial Economics} 10(1):125--145.

\bibitem[{Curtin et~al.(2005)Curtin, Qui, Hayslett-McCall, \protect\BIBand{}
  Bray}]{curtin2005integrating}
Curtin KM, Qui F, Hayslett-McCall K, Bray TM (2005) Integrating gis and maximal
  covering models to determine optimal police patrol areas. \emph{Geographic
  information systems and crime analysis}, 214--235 (IGI Global).

\bibitem[{D'Amico et~al.(2002)D'Amico, Wang, Batta, \protect\BIBand{}
  Rump}]{Steven2002}
D'Amico SJ, Wang SJ, Batta R, Rump CM (2002) A simulated annealing approach to
  police district design. \emph{Computers \& Operations Research}
  29(6):667--684.

\bibitem[{Edwards(2016)}]{edwards2016}
Edwards E (2016) Predictive policing software is more accurate at predicting
  policing than predicting crime.
  Https://www.aclu.org/blog/criminal-law-reform/reforming-police/predictive-policing-software-more-accurate-predicting.

\bibitem[{Egbert(2016)}]{egbert2016}
Egbert J (2016) Beat redistricting proposal 2016. Atlanta Police Department
  (internal report).

\bibitem[{Ferland \protect\BIBand{} Gu{\'e}nette(1990)}]{ferland1990decision}
Ferland JA, Gu{\'e}nette G (1990) Decision support system for the school
  districting problem. \emph{Operations Research} 38(1):15--21.

\bibitem[{Fritz(2018)}]{Fritz2018}
Fritz J (2018) What is the average police response time in the {U.S.}?
  (accessed {A}pril 30, 2019).
  Https://www.safesmartliving.com/home-security/average-police-response-time/.

\bibitem[{Gaetan \protect\BIBand{} Guyon(2010)}]{gaetan2010spatial}
Gaetan C, Guyon X (2010) \emph{Spatial statistics and modeling}, volume~90
  (Springer).

\bibitem[{Garfinkel \protect\BIBand{} Nemhauser(1970)}]{GaNe1970}
Garfinkel RS, Nemhauser GL (1970) Optimal political districting by implicit
  enumeration techniques. \emph{Management Science} 16(8):B--495--B--508.

\bibitem[{Gass(1968)}]{gass1968division}
Gass SI (1968) On the division of police districts into patrol beats.
  \emph{Proceedings of the 1968 23rd ACM National Conference}, 459--473 (ACM).

\bibitem[{Green \protect\BIBand{} Kolesar(2004)}]{green2004anniversary}
Green LV, Kolesar PJ (2004) Anniversary article: Improving emergency
  responsiveness with management science. \emph{Management Science}
  50(8):1001--1014.

\bibitem[{Grofman(1985)}]{Gr1985}
Grofman B (1985) Criteria for redistricting: A social science perspective.
  \emph{UCLA Law Review} 33:77--184.

\bibitem[{Habersham(2019)}]{HabershamA2019}
Habersham R (2019) Atlanta police hope changes to patrol zones shorten response
  times. \emph{The Atlanta Journal-Constitution}
  \url{https://www.ajc.com/news/local/atlanta-police-change-zoning-boundaries-help-with-response-times/R4hjqmYEmaveOgkW6NfRUO/}.

\bibitem[{Habersham \protect\BIBand{} Deere(2019)}]{HabershamB2019}
Habersham R, Deere S (2019) Buckhead residents confront mayor, police chief
  about crime. \emph{The Atlanta Journal-Constitution}
  \url{https://www.ajc.com/news/local/buckhead-residents-confront-mayor-police-chief-about-crime/5VYluSvFkIAmqUloAIIRSK/}.

\bibitem[{Hess et~al.(1965)Hess, Weaver, Siegfeldt, Whelan, \protect\BIBand{}
  Zitlau}]{hess1965nonpartisan}
Hess SW, Weaver J, Siegfeldt H, Whelan J, Zitlau P (1965) Nonpartisan political
  redistricting by computer. \emph{Operations Research} 13(6):998--1006.

\bibitem[{Kistler(2009)}]{kistler2009tucson}
Kistler A (2009) Tucson police officers redraw division boundaries to balance
  their workload. \emph{Geography \& Public Safety} 1(4):3--5.

\bibitem[{Larson(1972)}]{larson1972urban}
Larson RC (1972) \emph{Urban police patrol analysis}, volume~28 (MIT Press
  Cambridge, MA).

\bibitem[{Larson(1974{\natexlab{a}})}]{Larson1974}
Larson RC (1974{\natexlab{a}}) A hypercube queuing model for facility location
  and redistricting in urban emergency services. \emph{Computers \& Operations
  Research} 1(1):67 -- 95.

\bibitem[{Larson(1974{\natexlab{b}})}]{larson1974illustrative}
Larson RC (1974{\natexlab{b}}) Illustrative police sector redesign in district
  4 in boston. \emph{Urban Analysis} 2(1):51--91.

\bibitem[{Larson \protect\BIBand{} Odoni(1981)}]{larson1981urban}
Larson RC, Odoni AR (1981) \emph{Urban Operations Research} (Prentice Hall,
  NJ).

\bibitem[{Lee et~al.(2014)Lee, Goldsman, Kim, \protect\BIBand{}
  Tsui}]{lee2014spatiotemporal}
Lee ML, Goldsman D, Kim SH, Tsui KL (2014) Spatiotemporal biosurveillance with
  spatial clusters: control limit approximation and impact of spatial
  correlation. \emph{IIE Transactions} 46(8):813--827.

\bibitem[{Levine et~al.(2017)Levine, Tisch, Tasso, \protect\BIBand{}
  Joy}]{levine2017new}
Levine E, Tisch J, Tasso A, Joy M (2017) The {New York City} police
  department's domain awareness system. \emph{Interfaces} 47(1):70--84.

\bibitem[{Liberatore \protect\BIBand{}
  Camacho-Collados(2016)}]{liberatore2016comparison}
Liberatore F, Camacho-Collados M (2016) A comparison of local search methods
  for the multicriteria police districting problem on graph. \emph{Mathematical
  Problems in Engineering} 2016.

\bibitem[{McCormick(1976)}]{mccormick1976computability}
McCormick GP (1976) Computability of global solutions to factorable nonconvex
  programs: Part i—convex underestimating problems. \emph{Mathematical
  programming} 10(1):147--175.

\bibitem[{Mehrotra et~al.(1998)Mehrotra, Johnson, \protect\BIBand{}
  Nemhauser}]{Me1972}
Mehrotra A, Johnson EL, Nemhauser GL (1998) An optimization based heuristic for
  political districting. \emph{Management Science} 44(8):1100--1114.

\bibitem[{Mills(1967)}]{Mi1967}
Mills G (1967) The determination of local government electoral boundaries.
  \emph{Operations Research Quarterly} 18:243--255.

\bibitem[{Mitchell(1972)}]{mitchell1972optimal}
Mitchell PS (1972) Optimal selection of police patrol beats. \emph{Journal of
  Criminal Law, Criminology \& Police Science} 63:577.

\bibitem[{Morrill(1973)}]{Mo1973}
Morrill RL (1973) Ideal and reality in reapportionment. \emph{Annals of the
  Association of American Geographers} 63:463--477.

\bibitem[{Morrill(1976)}]{Mo1976}
Morrill RL (1976) Redistricting revisited. \emph{Annals of the Association of
  American Geographers} 66:548--556.

\bibitem[{Nagel(1972)}]{Na1972}
Nagel SS (1972) Computers and the law and policitics of redistricting.
  \emph{Polity} 5:77--93.

\bibitem[{{New York City Government}(2020)}]{NYPD2020}
{New York City Government} (2020) {NYPD} precincts and sectors.
  \url{https://data.cityofnewyork.us/Public-Safety/NYPD-Sectors/eizi-ujye}.

\bibitem[{Niemi(1990)}]{Ni1990}
Niemi ea R~G (1990) Measuring compactness and the role of a compactness
  standard in a test for partisan and racial gerrymandering. \emph{The Journal
  of Politics} 52(4):1155--1181.

\bibitem[{Perry et~al.(2013)Perry, McInnis, Price, Smith, \protect\BIBand{}
  Hollywood}]{perry2013predictive}
Perry WL, McInnis B, Price CC, Smith SC, Hollywood JS (2013) \emph{Predictive
  policing: The role of crime forecasting in law enforcement operations} (Rand
  Corporation).

\bibitem[{Piyadasun et~al.(2017)Piyadasun, Kalansuriya, Gangananda, Malshan,
  Bandara, \protect\BIBand{} Marru}]{piyadasun2017rationalizing}
Piyadasun T, Kalansuriya B, Gangananda M, Malshan M, Bandara HD, Marru S (2017)
  Rationalizing police patrol beats using heuristic-based clustering.
  \emph{2017 Moratuwa Engineering Research Conference (MERCon)}, 431--436
  (IEEE).

\bibitem[{R{\'\i}os-Mercado(2020)}]{rios2020optimal}
R{\'\i}os-Mercado RZ (2020) \emph{Optimal Districting and Territory Design},
  volume 284 (Springer).

\bibitem[{Rosen(1974)}]{Rosen1974}
Rosen S (1974) Hedonic prices and implicit markets: Product differentiation in
  pure competition. \emph{Journal of Political Economy} 82(1):34--55.

\bibitem[{Rosete(2019)}]{socialmedia}
Rosete A (2019) People's response to the beat redesign on {N}extdoor.
  \url{https://nextdoor.com/agency-post/ga/atlanta/atlanta-police-department/atlanta-police-make-changes-to-three-zone-boundaries-as-part-of-beat-redesign-104899950/}.

\bibitem[{Sarac et~al.(1999)Sarac, Batta, Bhadury, \protect\BIBand{}
  Rump}]{sarac1999reconfiguring}
Sarac A, Batta R, Bhadury J, Rump C (1999) Reconfiguring police reporting
  districts in the city of buffalo. \emph{{OR} Insight} 12(3):16--24.

\bibitem[{Saunders et~al.(2016)Saunders, Hunt, \protect\BIBand{}
  Hollywood}]{saunders2016predictions}
Saunders J, Hunt P, Hollywood JS (2016) Predictions put into practice: a
  quasi-experimental evaluation of chicago’s predictive policing pilot.
  \emph{Journal of Experimental Criminology} 12(3):347--371.

\bibitem[{Shirabe(2009)}]{Shirabe2009}
Shirabe T (2009) Districting modeling with exact contiguity constraints.
  \emph{Environment and Planning B: Planning and Design} 36(6):1053--1066.

\bibitem[{{United States Census Bureau}(2018)}]{UScensusbureau}
{United States Census Bureau} (2018) 2018 population estimates program.
  \url{https://factfinder.census.gov}.

\bibitem[{Vickrey(1961)}]{Vi1961}
Vickrey W (1961) On the prevention of gerrymandering. \emph{Political Science
  Quarterly} 76:105-- 110.

\bibitem[{Weaver \protect\BIBand{} Hess.(1963)}]{WeHe2963}
Weaver JB, Hess SW (1963) Districting modeling with exact contiguity
  constraints. \emph{The Yale Law Journal} 72:288--308.

\bibitem[{Yong(1988)}]{Yo1988}
Yong HP (1988) Measuring the compactness of legislative districts.
  \emph{Legislative Studies Quarterly} 13:105--115.

\end{thebibliography}

%

\newpage

\appendix

\vspace{.2in}
\section{Key Assumptions on Dispatching Rule}
We make the following assumptions on the dispatching rule:
(a) The number of patrol units allocated to each zone is equal to the number of beats; for simplicity, we do not consider backup units.
(b) In response to each call for police service, the nearest available patrol unit is dispatched to the scene of the request; a major incident may require multiple patrol units to respond, but we ignore such possibility since these incidents are rare.
(c) If no unit is available in the zone of the request, the call enters a queue with other backlogged calls, which is processed according to the first-come-first-served (FIFO) rule.
(d) Upon completion of the service, a patrol unit is either dispatched to a call waiting in the queue, or it immediately returns to its home beat if the queue is empty. In the formulation of the stochastic model below, we assume returning to the home beat is instantaneous.

\vspace{.2in}
\section{Police Patrol and Emergency Response Model}
\label{sec:simulation-model}

In this section, we present a stochastic model for the police patrol and emergency responses.
First we model the patrol units as a queueing system; then we derive expressions for several performance measures. The notations defined in this section are summarized in Table~\ref{tab:notation}.

\begin{table}[!ht]
\caption{Summary of Notations}
\label{tab:notation}
\resizebox{\textwidth}{!}{%
\begin{tabular}{lll}
\toprule[1pt]\midrule[0.3pt]
\bf Section & \bf Notation & \bf Description \\ \hline
Setup
& $\mathscr{I} = \{i=1,\dots, I\}$ & Set of all beats \\
 & $\mathscr{K} = \{k=1,\dots, K\}$ & Set of all zones \\
 & $D=(d_{ik})_{i\in\mathscr{I}, k\in\mathscr{K}}$ & Districting design matrix indicating if beat $i$ is assigned to zone $k$, $i \in \mathscr{I},~k \in \mathscr{K}$  \\
 & $\mathscr{I}^{(k)}$ & Set of beats assigned to zone $k \in \mathscr{K}$ \\
 & $N^{(k)}$ & Total number of beats assigned to zone $k \in \mathscr{K}$ \\ 
 & $\lambda_i$ & 911 calls-for-service arrival rate in beat $i \in \mathscr{I}$ \\
 & $\mu$ & Service rate of each response unit \\ \hline
Model
& 
$B^{(k)}, S^{(k)}$ & The unsaturated and saturated states for zone $k \in \mathscr{K}$ \\ 
 & $H^{(k)} = (\eta^{(k)}_{nj})$ & Optimal unit dispatched to beat $j \in \mathscr{I}^{(k)}$ under unsaturated states $B^{(k)}_n$ \\
 & $Q^{(k)}=(q^{(k)}_{nm})$ & Transition rate matrix associated with unsaturated states \\ 
 & $T=(\tau_{ij})$ & Average travel time from beat $i$ to beat $j$, $\forall i, j \in \mathscr{I}$ \\ \hline
Metrics 
& $E^{(k)}_{ij}$ & Set of states in which unit $i \in \mathscr{I}^{(k)}$ is an optimal unit to assign to a call from beat $j \in \mathscr{I}^{(k)}$ \\
 & $\rho^{(k)}_{ij}$ & Fraction of all calls that send unit $i \in \mathscr{I}^{(k)}$ to beat $j \in \mathscr{I}^{(k)}$ in zone $k$ \\
 & $\xi^{(k)}_i$ & Expected travel time of unit $i \in \mathscr{I}^{(k)}$ \\ 
 & $w^{(k)}(D)$ & Total workload of zone $k \in \mathscr{K}$ given districting $D$ \\ \midrule[0.3pt]\bottomrule[1pt]
\end{tabular}%
}
\end{table}

\subsection{Model Formulation}
\label{sec:model-desc}

\noindent\emph{Zone Districting.}
Consider the police districting problem with $K$ zones and $I$ beats ($K<I$). 
Let $k \in \mathscr{K} = \{1,...,K\}$ denote the set of zones. 
Let $i \in \mathscr{I} = \{1, \dots,I\}$ denote both the set of beats, as well as the set of police patrol units which are assigned to these beats.
We assume that the shape of each beat is given.
A particular zone districting design is thus a graph partition of $I$ beats into $K$ zones.
Let the binary decision variable $d_{ik} \in \{0, 1\}$ denote whether beat $i$ is assigned to zone $k$. 
The districting decision is represented by a matrix $D=(d_{ik}) \in \{0,1\}^{I \times K}$. For each $i \in \mathscr{I}$, we have $\sum_{k = 1}^{K} d_{ik} = 1$. Given a zone design $D$, the set of beats assigned to zone $k$ is denoted by $\mathscr{I}^{(k)} = \{i~|~d_{ik} = 1\} \subseteq \mathscr{I}$ and the number of beats in zone $k$ is denoted by $N^{(k)} = \sum_{i \in \mathscr{I}} d_{ik}$. 
Throughout this section, we assume the zone districting matrix $D$ is fixed. 

\vspace{.1in}
\noindent\emph{Hypercube Queue Model.}
For a given zone $k \in \mathscr{K}$, recall that the set of patrol units in this zone is denoted by $\mathscr{I}^{(k)}$. 
Let $N^{(k)} = |\mathscr{I}^{(k)}|$ be the total number of patrol units, which is equal to the number of beats in zone $k$.
We assume that the arrivals of 911 calls in beat $i\in \mathscr{I}$
follow a time-homogeneous Poisson process with rate $\lambda_i$. 
Let $\lambda^{(k)} = \sum_{i\in \mathscr{I}^{(k)}} \lambda_i$ be the aggregate call arrival rate in  zone $k$.
Denote the average travel time between two arbitrary beats as a matrix $T = (\tau_{ij}) \in \mathbb{R}_+^{I \times I}$, where $\tau_{ij}$ is the average travel time from beat $i$ to beat $j$.
The service time of each patrol unit is i.i.d.\ and follows an exponential distribution with mean $1/\mu$.
We require $\lambda^{(k)} < \mu N^{(k)}$ for the stability of this queueing system.

The system state depends on the status of all the units in this zone, and the number of calls in queue to be processed. Based on whether there is any available unit, the state space (denoted by $C^{(k)}$) can be divided into two parts: (1) \emph{unsaturated states}: the states where the queue for unprocessed calls is empty. These states can be represented by a hypercube in dimension $N^{(k)}$ (hence the name ``hypercube queue model'').
Each vertex of the hypercube corresponds to a state $B^{(k)} = (b_i)_{i \in \mathscr{I}^{(k)}}$, where unit $i$ is busy processing a call if $b_i=1$ and idle if $b_i=0$. 
Considering the state $B^{(k)}$ as a vector of binary numbers, we index the state by its corresponding decimal form $B^{(k)}_n$ where $n = \sum_{j = 1}^{N^{(k)}} 2^{j-1} b_j$. 
(2) \emph{saturated states}: the states where the queue for unprocessed calls is non-empty. By definition, all patrol units are busy under saturated states.
We denote these states by $\{S^{(k)}_n\}_{n\geq 1}$; the state $S^{(k)}_n$ indicates that there are exactly $n$ calls waiting in queue for zone $k$.


We now define the state transition rates.
For the saturated states, the transition rate from $S^{(k)}_{n}$ to $S^{(k)}_{n+1}$, as well as from the transition rate from $B^{(k)}_{2^{N^{(k)}}-1}$ to $S^{(k)}_{1}$, is $\lambda^{(k)}$. The transition rate from $S^{(k)}_{n+1}$ to $S^{(k)}_{n}$, as well as the transition rate from $S^{(k)}_{1}$ to $B^{(k)}_{2^{N^{(k)}}-1}$, is $\mu N^{(k)}$.

For the unsaturated states, we define the transition rate matrix $Q^{(k)} = (q^{(k)}_{nm})$  for zone $k$, where $q^{(k)}_{nm}$ (for $n \neq m$) denotes the rate departing from the $n$-th unsaturated state and arriving in the $m$-th unsaturated state. Diagonal entries $q^{(k)}_{nn}$ are defined such that $q^{(k)}_{nn} = -\sum_{m: m \neq n} q^{(k)}_{nm}$ and therefore the rows of the matrix sum to zero. 

\begin{figure}[!th]
\centering
\begin{subfigure}[h]{0.38\linewidth}
\includegraphics[width=\linewidth]{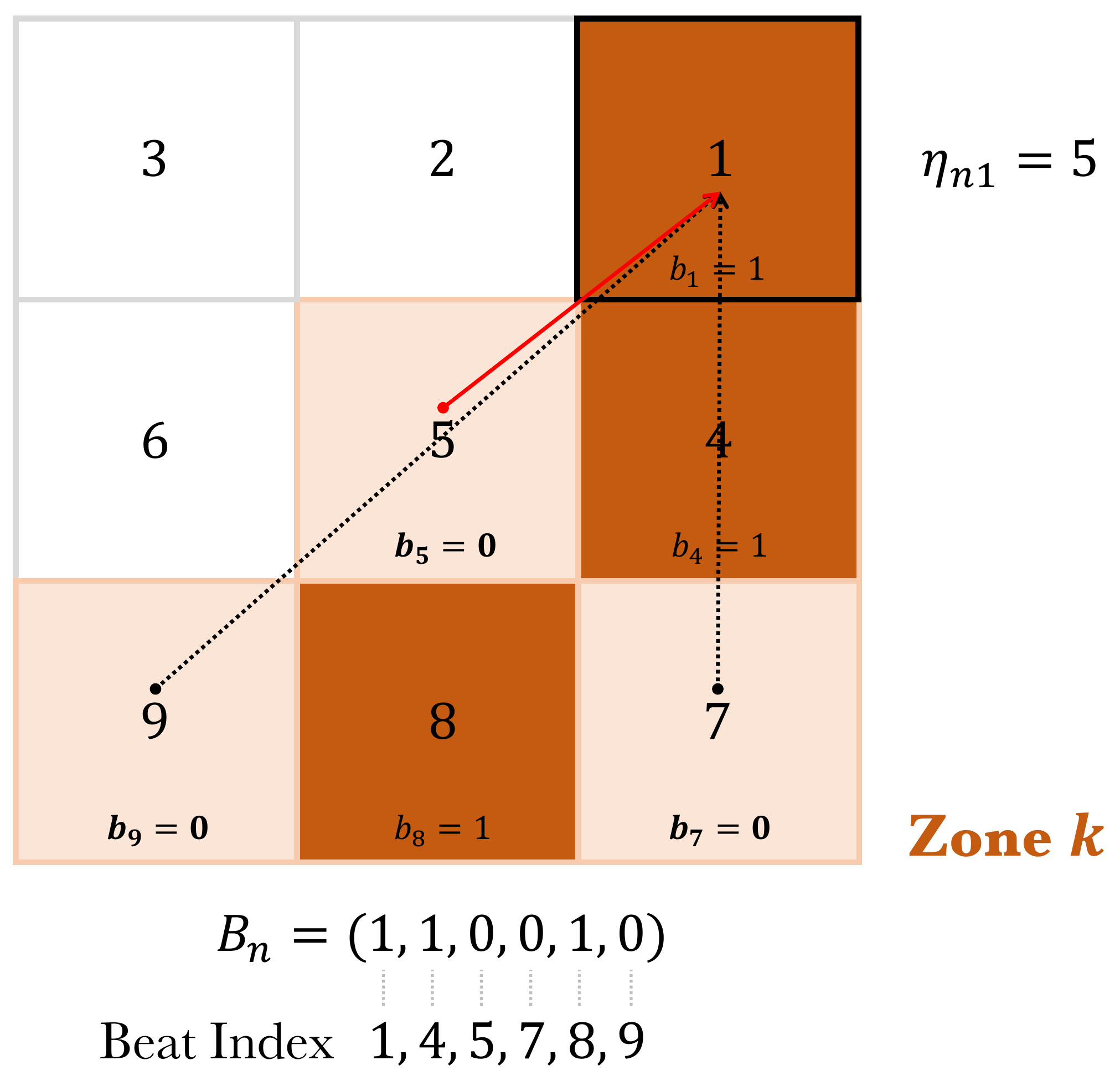}
\caption{$\eta_{n1}=5$ given $\mathscr{I}^{(k)} = \{1,4,5,7,8,9\}$.}
\label{fig:dispatch-exp-1}
\end{subfigure}
\hspace{.5in}
\begin{subfigure}[h]{0.38\linewidth}
\includegraphics[width=\linewidth]{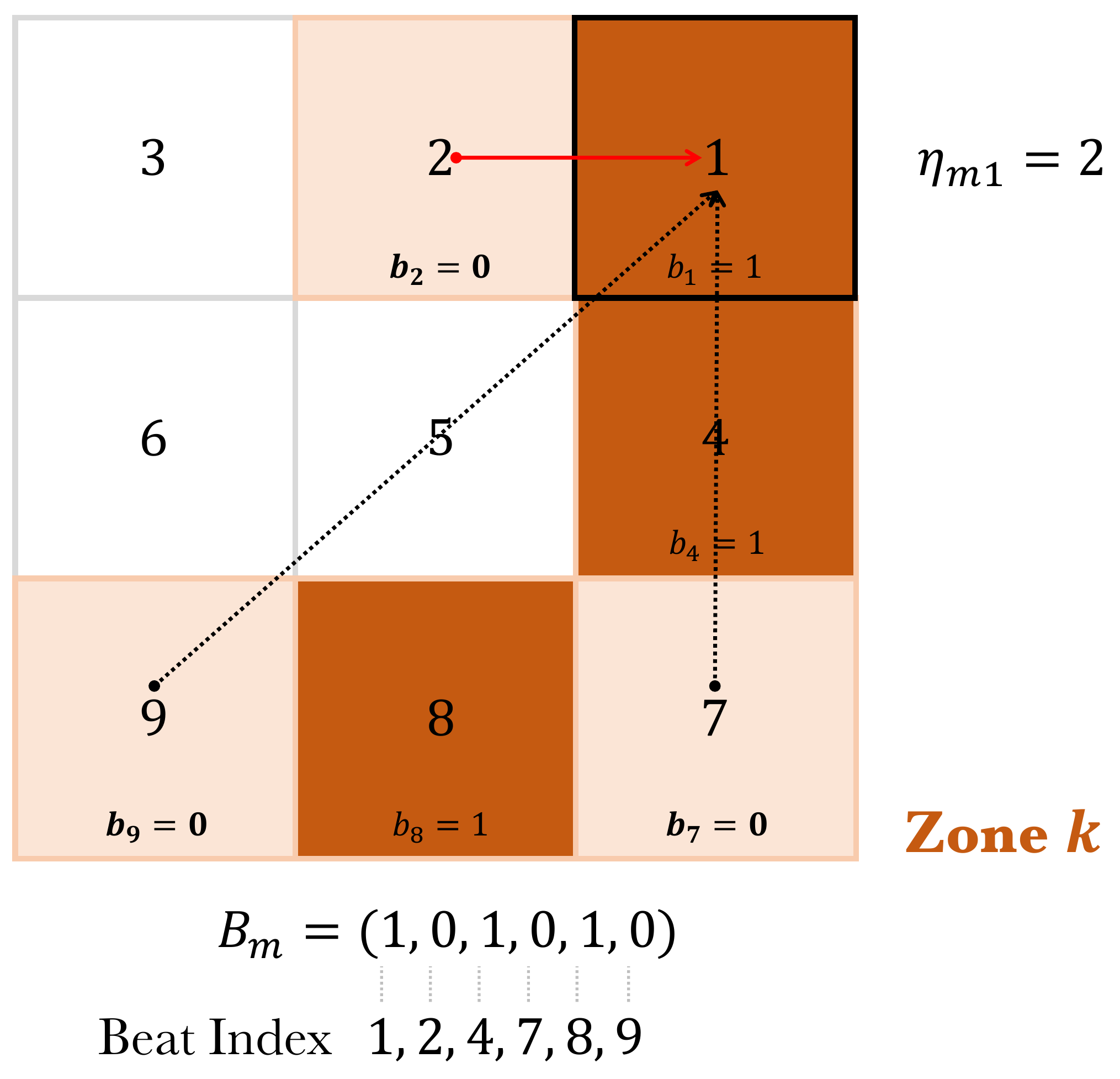}
\caption{$\eta_{m1}=2$ given $\mathscr{I}^{(k)} = \{1,2,4,7,8,9\}$.}
\label{fig:dispatch-exp-2}
\end{subfigure}
\vspace{6pt} 
\caption{Two examples show how districting decision affects the dispatch rule. Each numbered block represents a beat. The colored area represents a zone comprised of six beats. The darker blocks represent the corresponding patrol unit $i$ is busy ($b_i=1$). Euclidean distance is used as a proxy for travel time between beats.
The arrows indicate all possible dispatches, and the red line indicates the optimal dispatch with the shortest travel time.}
\label{fig:dispatch-exp}
\end{figure}

There are two classes of transitions on the hypercube:
\emph{upward transitions} that change a unit's status from idle to busy; 
\emph{downward transitions} that do the reverse. Let $d_{nm}$ denote the Hamming distance between two vertices $B_n$ and $B_m$ on the hypercube.  We define the ``upward'' Hamming distance as $d^+_{nm} = |B_n^C \cap B_m|$ and the ``downward'' Hamming distance as $d^-_{nm} = |B_n \cap B_m^C|$, where $|B|$ denote the number of binary ``one'' in $B$ (i.e., the number of busy units), and $B^C$ denote the complement of $B$. The downward transition rate from state $B_n$ to the adjacent states $B_m$ with $d_{nm}^- = 1$ is always $q^{(k)}_{nm} = \mu$. 
The upward transition rates, however, will depend on the dispatch rule (denoted by $H^{(k)}$) when a call is received.
As mentioned in the background,
we assume that when a call is ready to be processed, the nearest available unit with the minimum mean travel time to the location of the call is dispatched. That is, if the call is in beat $i$ and the system state is $B^{(k)}_n$, the index of the unit dispatched is $\eta^{(k)}_{ni} := \arg\min_{j: b_j = 0, b_j \in B^{(k)}_n} \tau_{ji}$ (we assume there are no ties in the mean travel time matrix).
The upward transition rate from state $B_n$ to an adjacent states $B_m$ with $d_{nm}^+ = 1$, can be obtained by $q^{(k)}_{nm} = \sum_{i:\eta_{ni}^{(k)} = j, \forall i} \lambda_i$, where $j$ is the index of the unit dispatched and satisfies $m = n + 2^{j-1}$.
To illustrate how this dispatch rule works, we give a detailed example in Figure~\ref{fig:dispatch-exp}.

The steady-state probabilities of this system are determined from the balance equations. 
Let $\mathbb{P}\{B^{(k)}_m\}$ and $\mathbb{P}\{S^{(k)}_n\}$
denote 
the probability that system is occupying state $B^{(k)}_m$ and $\mathbb{P}\{S^{(k)}_n\}$ under steady-state conditions.
For unsaturated states, we have
\begin{align}
  \left( \lambda^{(k)} + \mu |B^{(k)}_m| \right) \mathbb{P}\{B^{(k)}_m\} 
  ~=& \sum_{\{B^{(k)}_n: d^+_{nm} = 1\}} q_{nm} \mathbb{P}\{B^{(k)}_n\} 
  + \sum_{\{B^{(k)}_n: d^-_{nm} = 1\}} \mu \mathbb{P}\{B^{(k)}_n\},\quad \forall m = 1, \dots, 2^{N^{(k)}}-2,
 \nonumber \\
\left( \lambda^{(k)} + \mu N^{(k)} \right) \mathbb{P}\{B^{(k)}_{m}\} 
  ~=& \sum_{\{B^{(k)}_n: d^+_{nm} = 1\}} q_{nm} \mathbb{P}\{B^{(k)}_n\} 
  + \mu N^{(k)} \mathbb{P}\{S^{(k)}_1\},\quad \text{for } m = 2^{N^{(k)}}-1,  \label{eq:detailed-balance} 
\end{align}
where $\lambda^{(k)} =\sum_{i\in\mathscr{I}^{(k)}} \lambda_i$ is the aggregate call arrival rate in zone $k$.
For saturated states, we have
\[
 \left( \lambda^{(k)} + \mu N^{(k)} \right ) \mathbb{P}\{S^{(k)}_n\} 
  ~= \lambda^{(k)} \mathbb{P}\{S^{(k)}_{n-1}\} +  \mu N^{(k)} \mathbb{P}\{S^{(k)}_{n+1}\}, \quad \forall n \in \mathbb{Z}_+.
\]
In addition, we require that the probability sum to one, namely
\[
  \sum_{m=0}^{2^{N^{(k)}}-1} \mathbb{P}\{B^{(k)}_m\} + \sum_{n=1}^{\infty} \mathbb{P}\{S^{(k)}_n\} = 1.\]

We note that although the number of saturated states is infinite, it can be expressed in closed-form as in an M/M/c queue.
The main challenge for computing the steady-state distribution is that the number of \emph{unsaturated} states on the hypercube grows exponentially with the number of beats. 
In the following, 
we will propose an efficient computation method to approximate the steady-state distribution, which exploits the fact that the hypercube queue model has a sparse transition matrix.

\subsection{Measures of Performance}
\label{sec:measures}


\noindent\emph{Fraction of Dispatches.} 
To understand how often a patrol unit travels to different beats in the zone, it is critical to know how many calls in each beat have been assigned to a given unit. The fraction of dispatches in zone $k$ that send unit $i \in \mathscr{I}^{(k)}$ to beat $j \in \mathscr{I}^{(k)}$, denoted by $ \rho^{(k)}_{ij}$, can be divided into two terms: (1) $\rho^{(k,1)}_{ij}$ the fraction of dispatches that incurs no queue delay, and (2) $\rho^{(k,2)}_{ij}$ the fraction of dispatches that incurs a positive queue delay. The fraction of dispatches $\rho^{(k)}_{ij}$ can be expressed as
\begin{equation}
  \rho^{(k)}_{ij} = 
  \underbrace{\mystrut{4.ex} \sum_{B^{(k)}_n \in E^{(k)}_{ij}} \frac{\lambda_{j}}{\lambda^{(k)}} \mathbb{P}\{B^{(k)}_n\}}_{\rho^{(k,1)}_{ij}} + 
  \underbrace{\mystrut{4.ex}\frac{\lambda_{j} P'_Q}{\lambda^{(k)} N^{(k)}}}_{\rho^{(k,2)}_{ij}},
  \label{eq:fraction}
\end{equation}
where $P_Q' = \sum_{m=1}^{\infty} \mathbb{P}\{S^{(k)}_m\} + \mathbb{P}\{B^{(k)}_{2^{N^{(k)}}-1}\}$ specifies the probability that a randomly arriving call incurs a queue delay; the set $E^{(k)}_{ij} = \{B^{(k)}_n~|~\eta^{(k)}_{nj} = i, n=1,\cdots,2^{N^{(k)}}-1\}$ contains all the unsaturated states in which unit $i$ is the nearest available unit to assign to a call from beat $j$, which is determined by the dispatch rule $H^{(k)} = (\eta^{(k)}_{ni})$. 
Equation~\eqref{eq:fraction} can be found in \cite{Larson1974}, so we omit the proof.
Note that in Equation~\eqref{eq:fraction}, the number of beats in zone $k$ ($N^{(k)}$), the dispatch rule ($H^{(k)}$), and the steady-state probability ($\mathbb{P}\{B^{(k)}_n\}, \mathbb{P}\{S^{(k)}_m\}$) depends on the districting decision. 

\vspace{.1in}
\noindent\emph{Expected Travel Time.} 
Another key measure of performance is the travel time that a unit takes for a single dispatch. The travel time is affected by both the beat-to-beat travel time $\tau_{ij}$ and the fraction of dispatches $\rho^{(k)}_{ij}$. The beat-to-beat time $T=(\tau_{ij})$ can be estimated from the real traffic data, which will be discussed later.
The fraction of dispatches critically depends on the utilization factor of the hypercube queue system. If the utilization is low, most calls can be processed by units close to their home beats, resulting in short travel time. However, if the system is congested, calls must wait in the queue, which will be processed in an FCFS basis. Each patrol unit is equally likely to process calls that incur a queue delay, regardless of their locations. This will lead to longer travel time.

To derive the equation for the expected travel time for zone $k\in\mathscr{K}$, we denote the expected travel time of unit $i$ for a single dispatch as $\xi^{(k)}_i$, the expected travel time for a call in beat $j$ as $T^{(k)}_j$,
and the zone-wide mean travel time as $\bar{T}^{(k)}$. Using the result in the fraction of dispatches,
$\xi^{(k)}_i$ can be expressed as the weighted sum of mean travel times $\tau_{ij}$ for all beats $j$ in zone $k$ given the dispatch fraction $\rho^{(k)}_{ij}$. 
Specifically, the average travel time of unit $i$ in zone $k$ can be expressed \citep[see][]{Larson1974,larson1981urban} as
\begin{equation}
  \xi_i^{(k)} = \frac{\sum_{j \in \mathscr{I}^{(k)}} \rho^{(k,1)}_{ij} \tau_{ij} + \bar{T}^{(k)}_Q P'_Q/N^{(k)}}{\sum_{j \in \mathscr{I}^{(k)}} \rho^{(k,1)}_{ij} + P'_Q/N^{(k)}},
  \label{eq:travel-time}
\end{equation}
where $\bar{T}^{(k)}_Q = \sum_{i \in \mathscr{I}^{(k)}} \sum_{j \in \mathscr{I}^{(k)}} \lambda_{i} \lambda_{j} \tau_{ij} / (\lambda^{(k)})^2$ is the mean travel time for calls that incur a positive queue delay and $P'_Q$ is defined in Equation~\eqref{eq:fraction}. 
Similarly, the expected travel time for calls in beat $j$ is
\[
    T^{(k)}_j = \frac{\sum_{i=1}^{N^{(k)}} \rho^{(k,1)}_{ij}\tau_{ij}}{\sum_{i=1}^{N^{(k)}} \rho^{(k,1)}_{ij}} (1-P'_Q) + \sum_{i=1}^{N^{(k)}} \frac{\lambda_i \tau_{ij}}{\lambda^{(k)}}P'_Q.
\]
The zone-wide mean travel time is
\[
    \bar{T}^{(k)} = \sum_{i=1}^{N^{(k)}} \sum_{j \in \mathscr{I}^{(k)}}\rho^{(k,1)}_{ij}\tau_{ij}
    + P'_Q \bar{T}^{(k)}_Q.
\]

\vspace{.1in}
\noindent\emph{Response Time.}
The response time for a call includes both the travel time for the patrol unit to arrive at the scene and the waiting time if the call incurs a positive queue delay. Therefore, the average response time for calls in beat $j$ is
\[
     T^{(k)}_j + \sum_{m=1}^{\infty} \frac{m+1}{\mu N^{(k)}}\mathbb{P}\{S^{(k)}_m\} + \frac{1}{\mu N^{(k)}}\mathbb{P}\{B^{(k)}_{2^{N^{(k)}}-1}\}.
\]
As mentioned earlier, the infinite sum in the second term can be solved in closed form as in an M/M/c queue. The zone-wide mean response time is
\[
 \bar{T}^{(k)} + \sum_{m=1}^{\infty} \frac{m+1}{\mu N^{(k)}}\mathbb{P}\{S^{(k)}_m\} + \frac{1}{\mu N^{(k)}}\mathbb{P}\{B^{(k)}_{2^{N^{(k)}}-1}\}.
\]

\vspace{.1in}
\noindent\emph{Patrol Unit Workload.}
Finally, we consider the workload of the patrol units in each zone. In the following optimization model, 
we will consider districting designs that 
balance the workload across zones, so this metric will play a vital role in defining the optimization objective. The workload of patrol unit $i$ in zone $k$ is
\[
    w^{(k)}_i = \sum_{n: (B^{(k)}_n)_i = 1} \mathbb{P}\{B^{(k)}_n\} + \sum_{m=1}^{\infty}  \mathbb{P}\{S^{(k)}_m\}.
\]
Here, we introduce the districting decision $D$ in defining the zone workload and it will be used in the objective function of our redistricting optimization.
The total workload for all patrol units in zone $k$ is simply
\begin{equation}
    w^{(k)}(D) = \sum_{i\in \mathscr{I}^{(k)}} w^{(k)}_i.
    \label{eq:zone-workload}
\end{equation}
Note that the workload for zone $k$ depends on the decision variable $D$ through $\mathscr{I}^{(k)}$. A detailed explanation is given below.

\begin{figure}[!hbt]
\centering
\includegraphics[width=.95\linewidth]{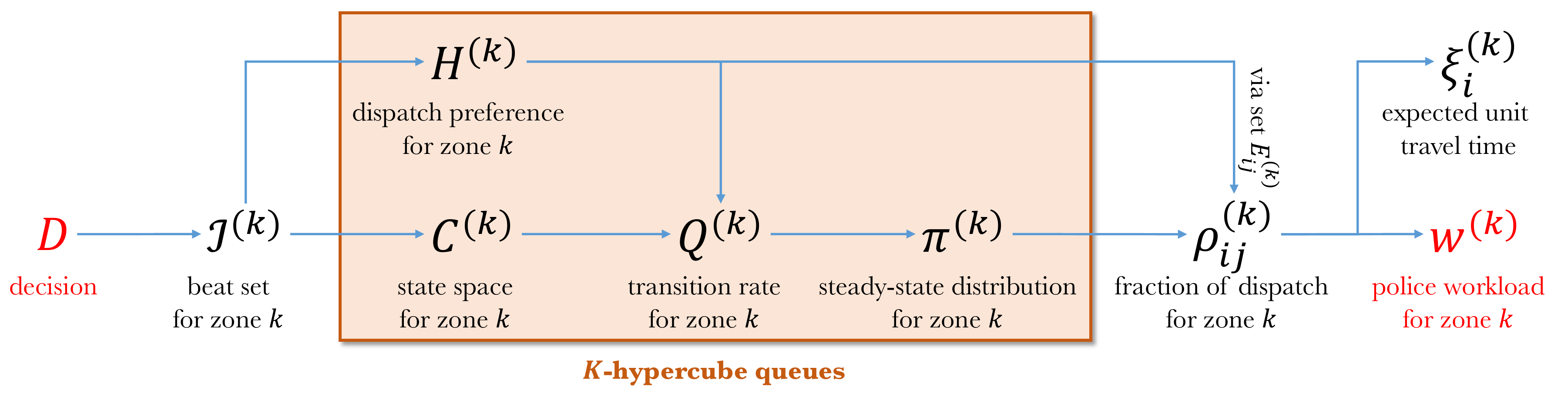}
\caption{The dependence diagram for the key variables in the $K$-hypercube model.}
\label{fig:dependence}
\end{figure}

\vspace{.1in}
\noindent\emph{Dependence of Performance Measures on the Districting Decision.}
A zone districting decision $D$ divides the entire geographical region into $K$ police zone, or $K$ hypercube queue models. The dynamics between the above performance measures and the beat-to-zone allocations are quite complex, as the parameters of $K$ hypercube queues are endogenously determined by the decision of beat allocation. As shown in Figure~\ref{fig:dependence}, the police zone workload $w^{(k)}(D)$ depends on decision $D$ through a long dependence chain. 
Specifically, given a particular decision $D$, the beats will be divided into $K$ zones, where zone $k \in \mathscr{K}$ corresponds to a set of beats $\mathscr{I}^{(k)}(D)$. This leads to a unique combination of $K$ hypercube queues, where their state spaces $\{C^{(k)}\}_{k\in\mathscr{K}}$ and dispatch preferences $\{H^{(k)}\}_{k\in\mathscr{K}}$ are jointly controlled by how decision $D$ partition the beats. By solving the linear equation system formed by the generator matrix $Q^{(k)}$ for each hypercube, we can obtain $K$ steady-state distributions $\{\pi^{(k)}\}_{k\in\mathscr{K}}$. Accordingly, we can compute the performance measures $\xi_i^{(k)}(D), \rho_{ij}^{(k)}(D)$ for each hypercube queue based on their steady-state distribution. 
The unknown police zone workload $w^{(k)}(D)$ for $K$ zones given a particular decision can be evaluated through the above process.

\subsection{Efficient Computation Method}
\label{sec:workload-estimation}

Performing workload estimation efficiently is critical to search for the optimal zone design in the solution space.  Take Atlanta as an example; there are more than 2,000 possible ways to shift the existing design even with less than only four changes that are allowed to make. 
Each zone in one particular design corresponds to a unique hypercube queue, where its zone workload needs to be re-evaluated separately. 

A critical step for obtaining police zone workload estimation given a decision is to compute the steady-state distribution. 

Let $\lambda^{(k)} = \sum_{i \in \mathscr{I}^{(k)}} \lambda_i$ be the total call arrival rate in zone $k$.
It is clear that for zone $k$, if the system is in some state $B^{(k)}_n$ with $|B^{(k)}_n|$ servers busy, the total transition rate ``downward'' to states with one less busy server is $|B^{(k)}_n| \mu$ and the total transition rate ``upward'' to states with one more busy server is $\lambda$ 
(except for $|B^{(k)}_n| = N^{(k)}$, in which case the upward transition rate is zero). 
Thus, if we aggregate states according to the numbers of servers busy, we obtain the well-known M/M/c multi-server queuing system. 
Then the steady-state probabilities in the saturated state set $\{S^{(k)}_n\}_{n \geq 1}$ is given by closed-form expression:
\[
  \mathbb{P}\{S^{(k)}_n\} = \frac{\frac{(\lambda^{(k)}/\mu)^{N^{(k)}}}{N^{(k)}!} \left(\frac{\lambda^{(k)}}{\mu N^{(k)}}\right)^n}{\sum_{s=0}^{N^{(k)}-1} \frac{(\lambda^{(k)}/\mu)^s}{s!} + \frac{(\lambda^{(k)}/\mu)^{N^{(k)}}}{N^{(k)}!} \frac{1}{1 - \lambda^{(k)} / (\mu N^{(k)})}}, \quad \text{for all }n \in \mathbb{Z}_+.
\]
Let the row vector $\pi^{(k)} = (\mathbb{P}\{B^{(k)}_m\})$ denote the steady-state probabilities in the unsaturated state set $C^{(k)}$. Then the balance equations \eqref{eq:detailed-balance} can be  written compactly as
\[
  \pi^{(k)} = \pi^{(k)} (I^{(k)} + Q^{(k)}/ \gamma^{(k)}),
\]
where $I^{(k)}$ is the identity matrix (of the same size as $Q^{(k)}$) and $\gamma^{(k)} = \max |q_{ii}^{(k)}|$ is the largest absolute value on the diagonal of $Q^{(k)}$. The matrix $(I^{(k)} + Q^{(k)}/ \gamma^{(k)})$ is the probability transition matrix of a discrete-time Markov chain obtained through uniformization of the original continuous-time Markov chain.  

Theoretically, the solution to this set of equations requires a matrix inversion. However, the size of the matrix is equal to $2^{2N^{(k)}}$, thus becoming huge for even moderate values of $N^{(k)}$. 
Fortunately, due to the sparsity of the transition rate matrix, we can find the stationary distribution ${\pi}^{(k)}$ iteratively by the power method. More specifically, given the transition rate matrix $Q^{(k)}$ and an initial vector $\pi^{(k)}_{0}$, compute $\pi^{(k)}_{t} = \pi^{(k)}_{t-1}(I^{(k)}+Q^{(k)}/\gamma^{(k)}) $ for $t=1,2,\dots$, until the distance (e.g.\ the sup norm) between $\pi^{(k)}_{t}$ and $\pi^{(k)}_{t-1}$ is small enough.

\vspace{.2in}
\section{Model Estimation and Workload Prediction}
\label{sec:model_estimation}

\subsection{Workload Estimation} 

We obtained the three data sets described above for five years from 2013 to 2017. This enabled us to fit the hypercube queue model 
to predict the future workload.

\vspace{.1in}
\noindent\emph{Call Arrival Rate Prediction.}
Predicting the call arrival rate is particularly challenging.
Although our operations model 
assumes the call arrival rates are time-homogeneous,
we observe in the actual data that the call arrival rates have a significant seasonality pattern and yearly trend (show in Figure~\ref{fig:zone-workload-validation} (b)), as well as correlation over adjacent geographical areas.
Therefore, we propose a spatio-temporal model to predict future call arrival rates.

We index the five years from 2013 to 2017 by $\ell \in \mathscr{L} = \{1, \dots, L\},~L=5$.
Let $\lambda_{i\ell} \in \mathbb{R}_+$ represent 911 call arrival rates of beat $i$ in year $\ell$.
Let vector $\boldsymbol{c}_{i\ell} \in \mathbb{R}^{M}$ represent the values of $M$ demographic factors of beat $i$ in year $\ell$.
The graph $\mathcal{G}$ is given by associating a node with every beat and connecting two nodes by an edge whenever the corresponding beats are geographically adjacent.
The set of adjacency pairs is defined by $\mathscr{A} = \{(i,j)\in \mathscr{I}: i,j\text{ are adjacent in }\mathcal{G}\}$. 
Then, we use the spatially lagged endogenous regressors \citep{Rosen1974} defined as
\begin{equation}
  \lambda_{i\ell} = \underbrace{\mystrut{3.ex} \sum_{(i,j)\in\mathscr{A}} \alpha_{ij} \lambda_{j\ell}}_{(\dagger)} + \underbrace{\mystrut{3.ex} \beta_0 \lambda_{i,\ell-1}}_{(\ddagger)} + \underbrace{\mystrut{3.ex} \sum_{t=0}^{p} \boldsymbol{c}_{i,\ell-t} \boldsymbol{\beta}_t}_{(\dagger\ddagger)} + \epsilon_{i}, \quad \ell = p + 1, \dots, L,
  \label{eq:lr-lam}
\end{equation}
where {$(\dagger)$ represents the influence of neighboring beats $j: (i, j) \in \mathscr{A}$,  $(\ddagger)$ represents the influence of the arrival rate in the last year, and $(\dagger \ddagger)$ represents the influence of demographic factors in the past $p$ years.} The $p$ is the total number of past years of data that we consider for fitting the regressor (in our case, $p=1$). The spatial coefficient matrix $A = (\alpha_{ij}) \in \mathbb{R}^{I \times I}$ specifies the spatial correlation between arrival rates of two beats. The temporal coefficient $\beta_0 \in \mathbb{R}$ specifies the influence of the arrival rate in the last year. The coefficient $\boldsymbol{\beta}_t \in \mathbb{R}^{M}, \forall 1 \le t \le p$ specifies the correlation of demographic factors.
Let $X_i$ denote a subset of the data $\{\{ \boldsymbol{c}_{i\ell} \}_{\ell\in\mathscr{L}}, \{\lambda_{j\ell}\}_{(i,j)\in\mathscr{A}, \ell\in\mathscr{L}}\}$.
The error $\epsilon_{i} \in \mathbb{R},~i\in\mathscr{I}$ is a random noise, 
$\mathbb{E}[\epsilon_i | X_{i}] = 0$ and $\text{Cov}[\epsilon_i, \epsilon_j | X_i, X_j] = \Sigma_\theta(i, j)$ for all $i, j\in \mathscr{I}$. 
The spatial covariance of the noise between two beats $i, j$ is determined by a correlation function $\Sigma_\theta$, which is a function of their spatial distance $s_{ij}$, and is parameterized by $\theta$. Some commonly used spatial models, including:
Gaussian model \citep{lee2014spatiotemporal},
Exponential model \citep{gaetan2010spatial}, and 
Mat\'{e}rn model \citep{gaetan2010spatial}.
Here we adopt the exponential model where $\Sigma_\theta(i, j) = \theta_1 \exp \{ - \theta s_{ij} \}
$, where $\theta$ is a pre-specified parameter. 
The other parameters of the model including $\{A, \beta_0, \{\boldsymbol{\beta}_t\}_{1\le t\le p}\}$ can be fitted by {maximum likelihood estimation} using historical data $\{X_i\}_{i \in \mathscr{I}}$. 

\vspace{.1in}
\noindent\emph{On-scene Time Estimation.}
As we assume the on-scene time of a response unit is independent of its service region, we investigate the distribution of the on-scene time of each unit in the same zone. By looking into the data, we have seen that the on-scene time of response units follows an exponential distribution shown in Figure~\ref{fig:on-scene}. We obtain the mean on-scene time $1/\mu = 31.2$ (minutes) by averaging the on-scene time of all dispatches in the past. 


\vspace{.1in}
\noindent\emph{Travel Time Estimation.}
We estimated the beat-to-beat travel time $T=(\tau_{ij})$ 
from the 911 call dispatch data. The estimation result is shown in Figure~\ref{fig:tau} 
by averaging the travel time of all 911 call dispatches that started in beat $i$ and ended in beat $j$.
In this figure, beats are ordered according to their relative location in the city. Thus, the beats in the same zone are consecutively numbered. Clearly, the figure reveals that local traveling within a neighborhood region usually taking a longer time than traveling across zones. 
There are two major reasons accounting for this phenomenon:
first, cross-zone calls usually have higher priority, which takes less time to response, and requires more patrol units to assist the call (that is why sometimes they need polices from other zones); 
second, response units usually resort to the highway system of the city when responding to the non-local calls, which may help reduce the travel time.

\vspace{.2in}
\section{Literature Review on Police Workload}

In our study, we consider balancing police workload as one of the design objectives.
Various definitions for police workload have been adopted in the literature.
\cite{mitchell1972optimal} defines police workload as the sum of service time and travel time. 
\cite{curtin2005integrating} uses the number or frequency of 911 calls occurring at each district as a proxy for the workload.
In \cite{bodily1978police} and \cite{Steven2002}, the workload is defined as the fraction of working time that an agent spends attending to calls. An equivalent measure is considered by \cite{benveniste1985solving}.
Response time is also an important performance measure representing the time between the arrival of a call for service and the arrival of a unit at the incident location. According to \cite{bodily1978police}, the reduction of the response time results in a number of beneficial effects such as (a) increased likelihood of intercepting a crime in progress; (b) deterrent effect on criminals; and (c) increased confidence of citizens in the police.

\vspace{.2in}
\section{Zone Design Optimization}

In this section, we develop an optimization model for zone redesign. We assume that all the beat areas are given, so the zone districting decision is equivalent to a graph partition where we allocate the beats to a fixed number of zones. Recall that the matrix $D = (d_{ik})_{I\times K}$ represents the beat allocation decisions, and $w^{(k)}(D)$ represents the mean police workload in zone $k\in\mathscr{K}$ given a districting design $D$.

The goal is to minimize the workload variance subject to some shape constraints (e.g.,\ contiguity and compactness) for each zone. The zone redesign problem can be expressed as
\begin{equation}
\begin{split}
\underset{D \in \{0,1\}^{I\times K}}{\text{minimize}} &\quad  f(D) \coloneqq \sum_{k=1}^{K} \left( w^{(k)}(D) - \frac{1}{K} \sum_{k^\prime=1}^{K} w^{(k^\prime)}(D) \right)^2\\
\mbox{subject to} 
&\quad \sum_{k = 1}^{K} d_{ik} = 1,\quad \forall i \in \mathscr{I},\\
&\quad \mbox{\it contiguity and compactness constraints for zone $k \in \mathscr{K}$}.
\end{split}
\label{eq:opt-objective-quadratic-1}
\end{equation}

Notice that the zone-level workload variance is chosen as the objective function in problem \eqref{eq:opt-objective-quadratic-1}, as APD suggested this metric for their zone redesign. 
However, one may also consider other performance metrics, such as mean response time or workload variance at the beat level (we provided expressions for these metrics in Appendix~B).
Our framework can be similarly adapted to solve these problems. 


\subsection{Zone Workload Approximation}
\label{sec:approx-workload}


We consider a local approximation of the workload function, which allows the objective function of the optimization problem \eqref{eq:opt-objective-quadratic-1} to be expressed explicitly.
Consider a given zone with an existing design. 
If we make a small change to this zone, e.g., by adding a few beats or removing a few beats, the workload of this zone can be approximated by its first-order Taylor expansion as
\begin{equation}
  w^{(k)}(D) = \theta_{0k} + \sum_{i\in\mathscr{I}} \theta_{ik} d_{ik} + \tilde{\epsilon}^{(k)}, \quad \forall k\in \mathscr{K}
  \label{eq:linear-approximation}
\end{equation}
where $\boldsymbol{\theta} = \{\{\theta_{0k} \in \mathbb{R}\}_{k\in\mathscr{K}}, \{\theta_{ik} \in \mathbb{R}\}_{i\in\mathscr{I}, k\in\mathscr{K}}\}$ denotes the coefficients of this approximation model, and $\tilde{\epsilon}^{(k)}, \forall k\in\mathscr{K}$ represent the approximation errors.

\begin{figure}[!hb]
\centering
\begin{subfigure}[h]{0.23\linewidth}
\includegraphics[width=\linewidth]{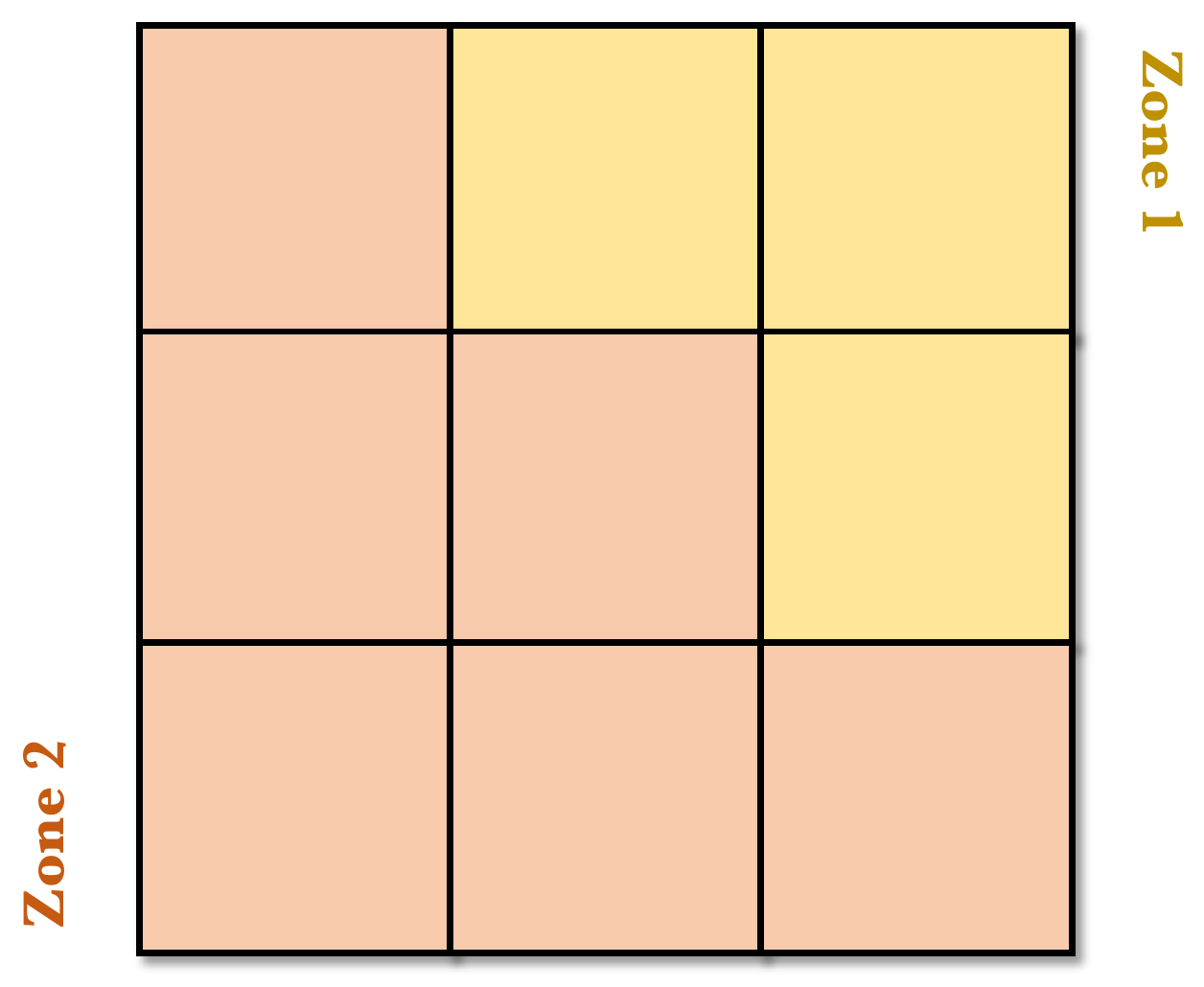}
\caption{$D_0$ without shift}
\label{fig:local-decision-no-swap}
\end{subfigure}
\begin{subfigure}[h]{0.23\linewidth}
\includegraphics[width=\linewidth]{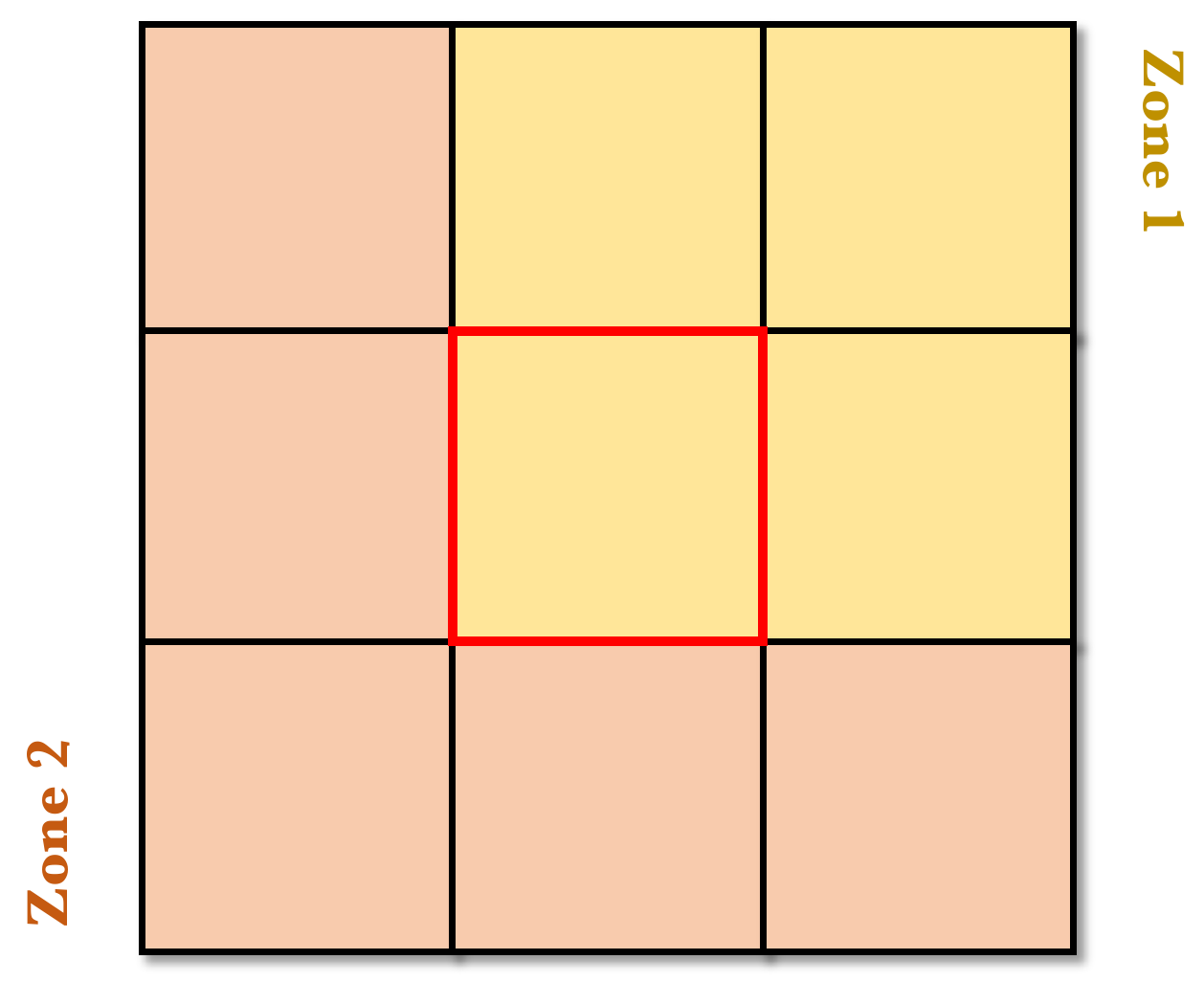}
\caption{$D_1$ with 1 shift}
\label{fig:local-decision-1-swap}
\end{subfigure}
\begin{subfigure}[h]{0.23\linewidth}
\includegraphics[width=\linewidth]{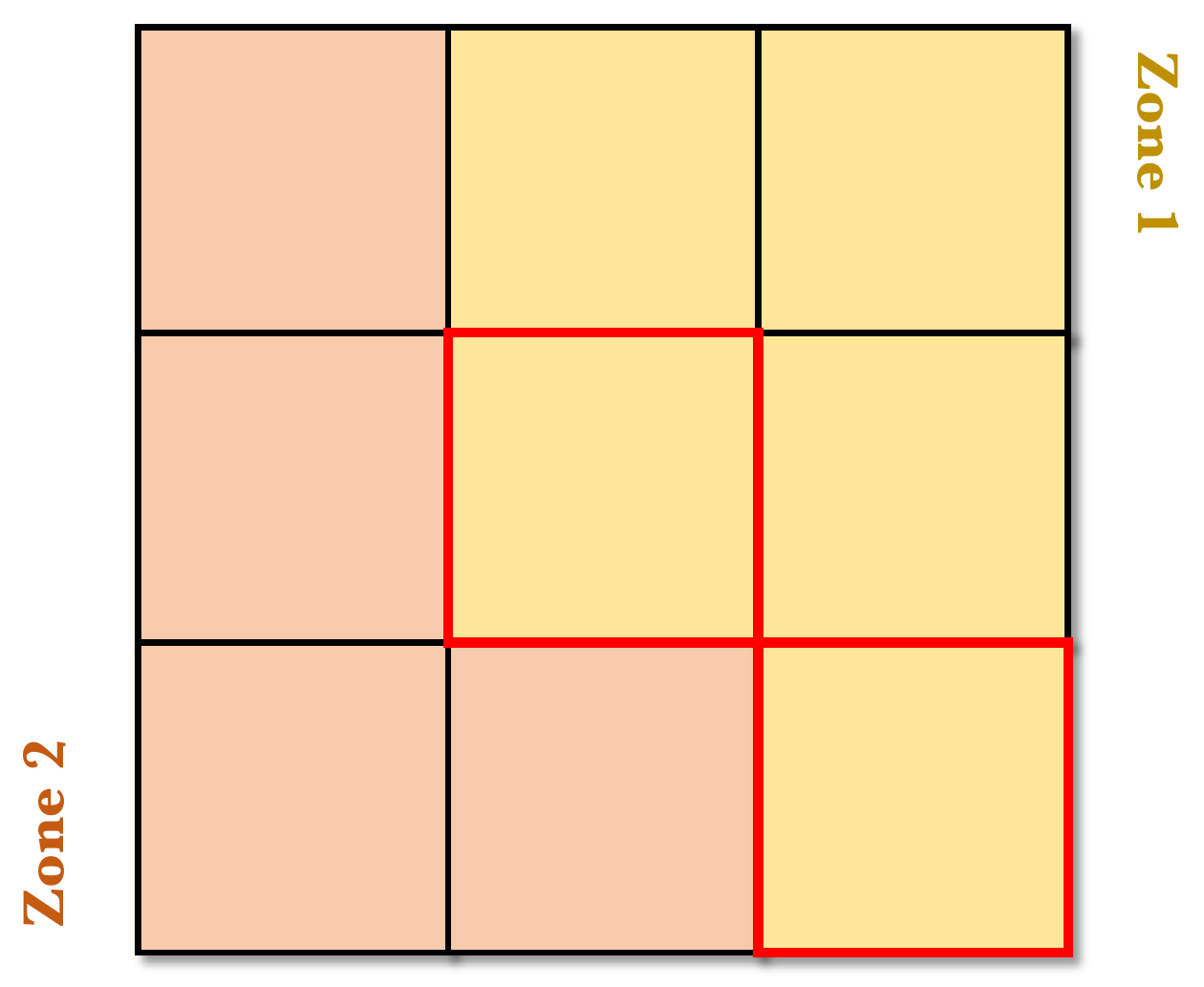}
\caption{$D_2$ with 2 shifts}
\label{fig:local-decision-2-swap-a}
\end{subfigure}
\begin{subfigure}[h]{0.23\linewidth}
\includegraphics[width=\linewidth]{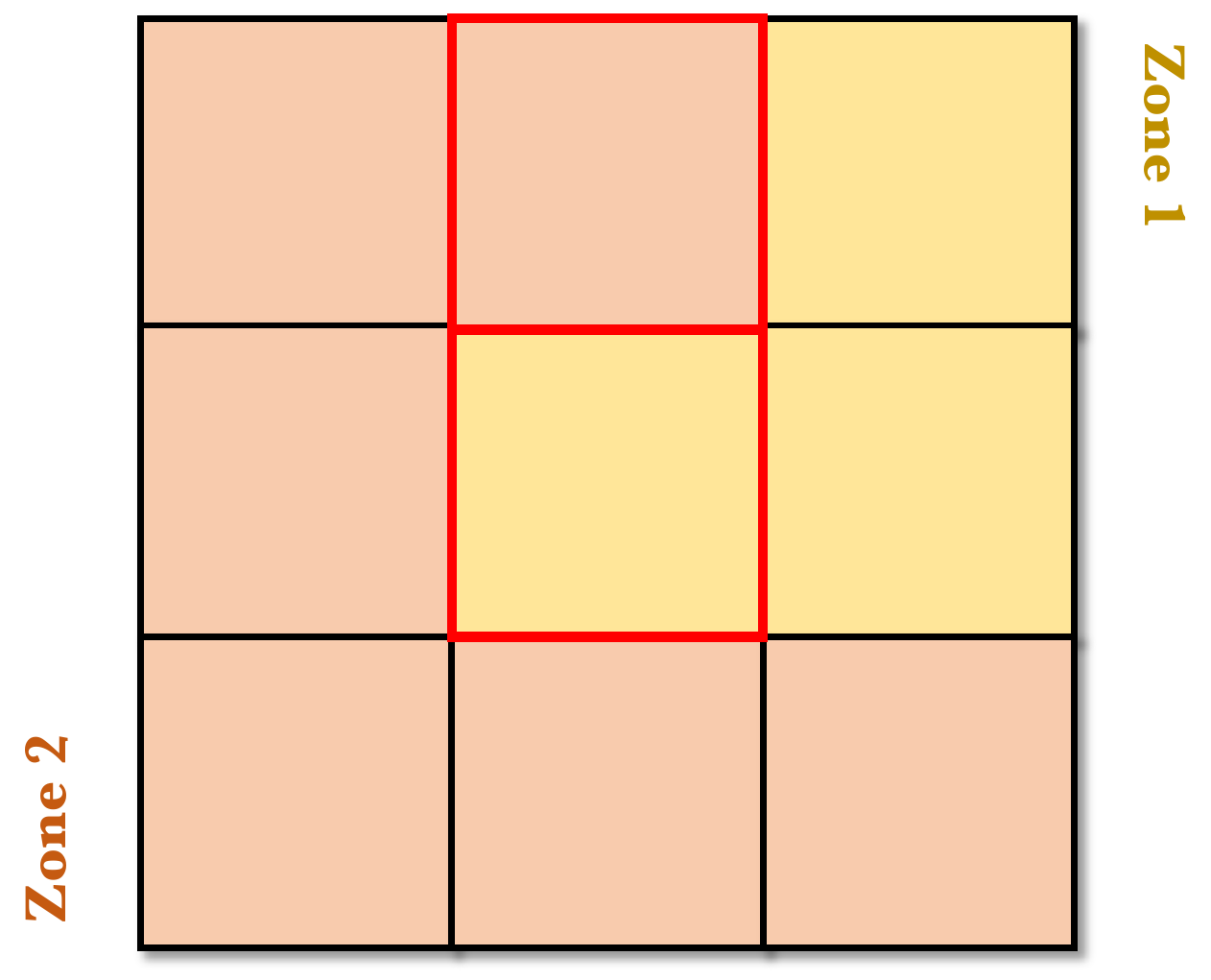}
\caption{$D_3$ with 2 shifts}
\label{fig:local-decision-2-swap-b}
\end{subfigure}
\vspace{10pt} \\
\caption{Examples show three local decisions by making 1, 2, 2 shifts, respectively. Each square block represents a single beat. Blocks with the same color are in the same zone. The red boxes indicate the shifts based on the original decision. (\ref{fig:local-decision-no-swap}) shows the original decision $D_0$ without any shift; (\ref{fig:local-decision-1-swap}) shows the decision $D_1$ with only 1 shift on the boundary; (\ref{fig:local-decision-2-swap-a}, \ref{fig:local-decision-2-swap-b}) show two decisions $D_2, D_3$ with different 2 shifts on the boundary.
}
\label{fig:local-decision}
\end{figure}

We estimate the parameters in Eq~\eqref{eq:linear-approximation} as follows.
First, for each of the six zones in Atlanta, we sample 1,000 designs that are created using small perturbations from the existing design. In particular, these samples are generated
by randomly changing a few beats on the boundary of the existing zone design, see Figure~\ref{fig:local-decision}. 
Obtaining a single sample involves procedures as follows: (a) generate a random design by changing a limited number of beat assignments; (b) compute the steady-state probabilities of six zones by solving \eqref{eq:detailed-balance}; (c) obtain simulated zone workload by computing \eqref{eq:zone-workload}. The computational bottleneck of getting one sample lies on step (b), which might take 2 minutes or so, depending on the number of the beat shifts and the computer-processing power. The entire generating process takes around 33 hours. 
Next, we estimate the parameters of the linear function in Eq~\eqref{eq:linear-approximation} using generated data samples.
More specifically, we apply 
the least squares method to obtain the estimated values
$\hat{\boldsymbol{\theta}}$, where the residual is $R^2 = 0.987$. 
Finally, we replace the exact workload function $w^{(k)}(D)$ with its approximation $\hat{w}^{(k)}(D|\hat{\boldsymbol{\theta}}) = \hat{\theta}_{0k} + \sum_{i\in\mathscr{I}} \hat{\theta}_{ik} d_{ik}$. The objective function $f(D)$ of the optimization problem~\eqref{eq:opt-objective-quadratic-1} is thus approximated by
\begin{equation}
\tilde{f}(D | \hat{\boldsymbol{\theta}}) = \sum_{k=1}^{K} \left( \hat{w}^{(k)}(D | \hat{\boldsymbol{\theta}}) - \frac{1}{K} \sum_{k^\prime=1}^{K} \hat{w}^{(k^\prime)}(D | \hat{\boldsymbol{\theta}}) \right)^2. 
\label{eq:opt-objective-quadratic-2}
\end{equation}

\begin{figure}[!htb]
\centering
\includegraphics[width=.5\linewidth]{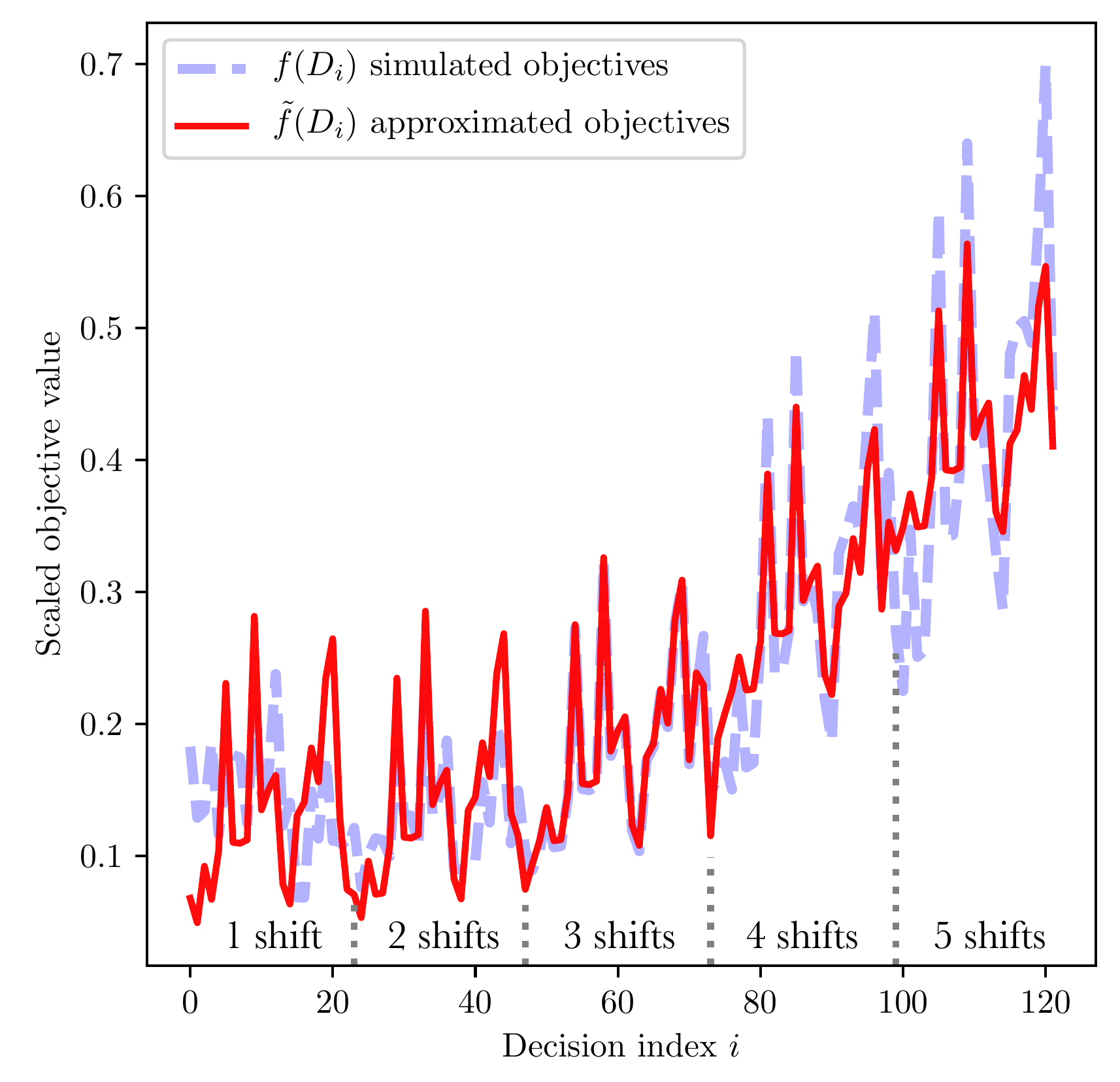}
\caption{The approximated and simulated objective values versus a subset of $\{D_i\}$ that differs to each other by shuffling a small number of beats. Decisions indexed by $0 \le i < 23$, $23 \le i < 47$, $47 \le i < 73$, $73 \le i < 99$, $99 \le i < 122$ correspond to randomly generated designs with 1, 2, 3, 4, and 5 shifts, respectively, comparing to the existing zone plan regarding the beat assignment.}
\label{fig:sim-approx-res}
\end{figure}

In Figure~\ref{fig:sim-approx-res}, we evaluate the effectiveness of our approximation by plotting the values of the objective function (i.e., workload variance) and corresponding approximations for 122 randomly selected \emph{out-of-sample} designs. These designs are also created based on the exiting design in Atlanta: 23 of them have one beat shift, 24 of them have two beats shift, 26 of them have three and four beats shift, and the other have five-beat shifts. Note that these designs are uniformly sampled from all possibilities and not necessarily optimal for the objective. 
The blue dash line plots the exact objective value $\{f(D_i)\}$, calculated by simulation. The red curve plots the approximated workload value $\{\tilde{f}(D_i|\hat{\boldsymbol{\theta}})\}$. The result shows that the approximated objective values match well with the true objective values for these randomly generated designs. 

\subsection{Linearizing Objective Function}
\label{sec:linear-obj}

With the parametric approximation \eqref{eq:opt-objective-quadratic-2}, the objective function measuring the workload variance becomes a \emph{quadratic} function of the binary decision variables $D$.
Using a standard technique for quadratic binary variables,
this quadratic function can be expressed as a linear form by 
by introducing auxiliary variables $e_{ijkk'}$ with additional constraints
\[
    e_{ijkk'} \leq d_{ik},\; e_{ijkk'} \leq d_{jk'},\; e_{ijkk'}\geq d_{ik} + d_{jk'} -1, e_{ijkk'}\geq 0 \quad \text{for all }i,j\in\mathscr{I}, k,k' \in \mathscr{K}.\]
These constraints ensure that $e_{ijkk'}= d_{ik}d_{jk'}$.
Therefore, the objective function can be written as
\begin{equation} \label{eq:obj-approx}
\begin{aligned}
\tilde{f}(D|\boldsymbol{\hat{\theta}}) = &~\sum_{k=1}^{K}\biggl[
c^{(k)} + \sum_{i=1}^I \sum_{j=1}^I \hat\theta_{ik} \hat\theta_{jk} e_{ijkk} + 2 \hat\theta_{0k} \sum_{i=1}^I \hat\theta_{ik} d_{ik} - \\
&\qquad \frac{2}{K} (\sum_{k'=1}^{K} \hat\theta_{0k'}) \sum_{i = 1}^I \hat\theta_{ik} d_{ik} - \frac{2}{K} \hat\theta_{0k} \sum_{k'=1}^{k} \sum_{i=1}^I \hat\theta_{ik'} d_{ik'} - \frac{2}{K} \sum_{i=1}^I \sum_{j=1}^I \sum_{k'=1}^{K} \hat\theta_{ik} \hat\theta_{jk'} e_{ijkk'} + \\
&\qquad \frac{2}{K^2} (\sum_{k'=1}^{K} \hat\theta_{0k'}) \sum_{k'=1}^{K} \sum_{i=1}^I \hat\theta_{ik'} d_{ik'} + \sum_{k'=1}^{K} \sum_{k''=1}^{K} \sum_{i=1} \sum_{j=1} \hat\theta_{ik'} \hat\theta_{jk''} e_{ijk'k''}
\biggr],
\end{aligned}
\end{equation}
where $c^{(k)} = \hat\theta_{0k}^2 + 2 \hat\theta_{0k} \sum_{k'=1}^{K} \hat\theta_{0k'}/ K + (\sum_{k'=1}^{K} \hat\theta_{0k'})^2 / K^2$ is a constant term for all $k\in\mathscr{K}$. 
The reformulated problem itself contains more than 240,000 variables and is characterized by a large search space with more than $2.4\times 10^{59}$ possible solutions. 
To tackle the computational complexity of solving this large-scale MILP, consider the implementation feasibility, and avoid drastic changes from the existing zone design, we use a fast local search method based on simulated annealing.
The optimization process takes about 30 to 50 minutes depending on where the problem converged given a single laptop's computational power with quad-core processors, which speed up to 4.7 GHz.

\subsection{Contiguity and Compactness Constraints}
\label{sec:constraints}


\noindent\emph{Contiguity Constraints.} 
We formulate the contiguity constraints using a network flow approach \citep[see][]{Shirabe2009}.  
We define a flow network, where each beat is represented by a node.
An arc from node $i$ to node $j$ exists if the beats associated with these two nodes are geographically adjacent.
Let $\mathscr{A}$ represent the set of arcs.

We now define flow constraints for each zone $k\in\mathscr{K}$ to ensure contiguity.
Let $n_{\max}$ be the maximum number of beats allowed in a zone.
If beat $i$ belongs to zone $k$, i.e., $d_{ik}=1$, we send one unit of flow into node $i$. Let $\upsilon_{ijk}$ be the amount of flow from beat $i$ to beat $j$ if $(i,j)\in\mathscr{A}$ and both beats belong to zone $k$.
In addition, for each zone $k$, we specify a special ``sink node.'' Let $h_{ik}$ be 1 if beat $i\in\mathscr{I}$ is selected as the sink node in zone $k\in\mathscr{K}$; otherwise, let $h_{ik}$ be 0. The sink node receives one unit of flow from every other node in the same zone. 
This ensures that there is a path of positive flow from any node in the zone to the sink node, implying contiguity.
In sum, we formulate the contiguity constraints as
\begin{subequations}\label{eq:opt-con-1}
\begin{align}
  \sum_{j: (i, j) \in \mathscr{A}} \upsilon_{ijk} - \sum_{l: (l, i) \in \mathscr{A}} \upsilon_{lik} & \ge d_{ik} - n_{\max} h_{ik}, & \forall i \in \mathscr{I}, k \in \mathscr{K}, \label{eq:opt-con2}\\
    \upsilon_{ijk} + \upsilon_{jik} & \le (n_{\max}-1) d_{ik}, & \forall (i,j)\in \mathscr{A}, k\in\mathscr{K}, \label{eq:opt-con7}\\
  \sum_{i=1}^{N} h_{ik} & = 1, & \forall k \in \mathscr{K}, \label{eq:opt-con4}\\
  h_{ik} - d_{ik} & \le 0, & \forall i \in \mathscr{I}, k \in \mathscr{K} \label{eq:opt-con6}\\
   d_{ik}, h_{ik} & \in \{0, 1\}, & \forall i\in\mathscr{I}, k\in\mathscr{K}, \label{eq:opt-con9}\\
  \upsilon_{ijk} & \ge 0, & \forall (i, j)\in\mathscr{A}, k\in\mathscr{K}. \label{eq:opt-con10}
\end{align}
\end{subequations}
Specifically, constraints \eqref{eq:opt-con2} represent the net outflow from all nodes. The two terms on the left indicate, respectively, the total outflow and total inflow of beat $i$. If beat $i$ is included in zone $k$ and it is not a sink, then we have $d_{ik}=1$, $h_{ik}=0$, and thus beat $i$ must have an outflow of one unit. If beat $i$ is included in zone $k$ and is a sink, then we have $d_{ik}=1$, $h_{ik}=1$, and thus beat $i$ has a negative outflow lower bounded by $1-n_{\max}$. 
Constraints \eqref{eq:opt-con7} ensure that if beat $i$ is not in zone $k$, i.e.\ $d_{ik}=0$,
there must be no flow into or out of beat $i$ with any nodes in zone $k$; if beat $i$ is in zone $k$ (i.e.\ $d_{ik}=1$), the flow between two beats does not exceed $n_{\max}-1$. 
Constraints \eqref{eq:opt-con4} ensure that each zone must have a unique sink node. 
Constraints \eqref{eq:opt-con6} ensure that beat $i$ must belong to zone $k$ if it is the sink for that zone. These constraints together ensure that each zone is contiguous.
  
\vspace{.1in}
\noindent\emph{Compactness Constraints.} 
We consider two ways to enforce the compactness of each zone: distance compactness \citep{Ni1990, Yo1988} and shape compactness \citep{GaNe1970}.
For distance compactness, we require the distance between any two beats in a zone to be bounded by some parameters.
For shape compactness, we require 
the square of the diameter divided by the zone's area to be upper bounded by another parameter. 
These constraints can be formulated as follows.

For each beat $i \in \mathscr{I}$, let $A_i$ be its area. For each pair of beats $(i,j) \in \mathscr{I}\times \mathscr{I}$, let $l_{ij}$ be the square of the distance between the center of these beats. 
We choose two parameters $\zeta_1, \zeta_2>0$ to control the level of distance compactness and shape constraints, respectively.  
Then we have
\begin{subequations}\label{eq:opt-con-2}
\begin{align}
  l_{ij}e_{ijkk} & \leq \zeta_1, & \forall i\in \mathscr{I}, j\in \mathscr{I}, k\in \mathscr{K}, \label{eq:opt-con11}\\
  l_{ij}e_{ijkk} & \leq \zeta_2 \sum_{i=1}^{K} d_{ik} A_i,  & \forall i\in \mathscr{I}, j \in \mathscr{I}, k\in \mathscr{K}. \label{eq:opt-con12}
\end{align}
\end{subequations}
Recall that $e_{ijkk}=d_{ik}d_{jk}$ is equal to 1 if both beat $i$ and beat $j$ are in zone $k$, and is equal to 0 otherwise. Constraints \ref{eq:opt-con11} represent the distance compactness criterion. Constraints \ref{eq:opt-con12} represent the shape compactness criterion, where $\sum_{i=1}^{K} d_{ik} A_i$ is the area of zone $k\in\mathscr{K}$.







\end{document}